\documentclass[a5paper, 10pt]{article}

\usepackage{cite, amsmath, amssymb, booktabs, lastpage, graphicx}
\usepackage[hidelinks]{hyperref}

\usepackage{geometry}
\usepackage{setspace}

\usepackage{amsthm}

\theoremstyle{plain}
\newtheorem{Lemma}{Lemma}
\newtheorem{Theorem}{Theorem}

\newtheorem{Corollary}{Corollary}

\newtheorem{Assumption}{Assumption}

\theoremstyle{remark}

\theoremstyle{definition}
\newtheorem{Definition}{Definition}

\usepackage{graphicx}
\makeatletter
\renewcommand{\maketitle}{
	\begin{center}
    \rule[.2em]{\textwidth}{0.353mm}
		\begin{minipage}[m]{0.35\textwidth}
			{\scriptsize
				\begin{center}
					\textsf{\textbf{\huge MATCH}\\
						\textit{Communications in Mathematical\\
							and in Computer Chemistry}
					}
			\end{center}}
		\end{minipage}\hfill
		\begin{minipage}[m]{0.65\textwidth}
			\begin{flushright}
				\baselineskip=10px
				{\scriptsize\sffamily{\itshape  MATCH Commun. Math. Comput. Chem.}
				\vol{}
				(\pubyear)
					\thepage--\pageref{LastPage}}\\
				{\scriptsize\sffamily {\bfseries ISSN:} 0340--6253}\\
				{\scriptsize\sffamily \textbf{doi:} \doi{10.46793/match.xx-x.xxxX}}
			\end{flushright}
		\end{minipage}
		\rule[1em]{\textwidth}{.353mm}
		\baselineskip=0.30in
		{\Large\bfseries \@title} \par
		\vspace{5mm}
		\baselineskip=0.2in
		{\large\bfseries \@author}\par
		\vspace{1mm}
		{\it \@address} \par
		{\small\tt \@email} \par
		\vspace{3mm}
		{\small (Received \@date)} \par
	\end{center}
	\vspace{3mm}
}
\newcommand{\address}[1]{\def\@address{#1}}
\newcommand{\email}[1]{\def\@email{#1}}
\newcommand{\doi}[1]{\href{https://doi.org/#1}{#1}}

\newcommand{\vol}{\textbf{00}}
\newcommand{\pubyear}{0000}

\makeatother

\newcommand{\acknowledgment}[1]{\vspace{5mm}\singlespacing
	{\noindent\textbf{\textit{Acknowledgment\/}:} #1}
}

\usepackage{fancyhdr}

\usepackage[margin=1cm,%
			font=footnotesize,%
			format=hang,%
			labelsep=period,%
			labelfont=bf]{caption}

\title{Asymptotic Distribution of Degree--Based Topological Indices}
\author{Mingao Yuan 
 }

\address{$^a$Department of Statistics, North Dakota State University,  USA\\

    }

\email{mingao.yuan@ndsu.edu 
}

\date{May 30, 2023}

\allowdisplaybreaks

\begin{document}

\maketitle

\begin{abstract}
Topological indices play a significant role in mathematical chemistry. Given a graph $\mathcal{G}$ with vertex set $\mathcal{V}=\{1,2,\dots,n\}$ and edge set $\mathcal{E}$, let $d_i$ be the degree of node $i$. The degree-based topological index is defined as $\mathcal{I}_n=$ $\sum_{\{i,j\}\in \mathcal{E}}f(d_i,d_j)$, where $f(x,y)$ is a symmetric function.
In this paper, 
we investigate the asymptotic distribution of
the degree-based topological indices of a heterogeneous  Erd\H{o}s-R\'{e}nyi random graph. We show that after suitably centered and scaled, the  topological indices converges in distribution to the standard normal distribution. Interestingly, we find that the general Randi\'{c} index with $f(x,y)=(xy)^{\tau}$ for a constant  $\tau$ exhibits a phase change at $\tau=-\frac{1}{2}$.
\end{abstract}

\onehalfspacing

\section{Introduction}
\label{S:1}

Topological index is a numerical parameter of a graph. It is graph invariant and characterizes the topology of a graph. Topological indices are used to model many physicochemical properties in QSAR \cite{AZG19,  BBGG00,G13}. One of the most important types of topological indices is the  degree-based topological indices, which are defined as a function of the degrees of nodes in a graph \cite{G13}. 

The first degree-based topological index is the Randi\'{c} index \cite{R75}. It measures the branching extent of a graph.  The Randi\'{c} index is the most popular and most
studied index among all
topological indices. It plays a central role in understanding quantitative structure-property and structure-activity relations in chemistry and pharmocology \cite{R08,RNP16}. Moreover, the Randi\'{c} index possesses a wealth of non-trivial and interesting mathematical properties \cite{BE98,BES99,CFK10,DSG17,LS08}. In addition,  the Randi\'{c} index finds countless applications in network (graph) data analysis. For instance, it was used to quantify the similarity between different networks or subgraphs of the same network \cite{FT13}, it serves as a quantitative characterization of network heterogeneity \cite{E10}, and \cite{DMRSV18,DMMRSF21} used the Randi\'{c} index to measure graph robustness.

 Motivated by the Randi\'{c} index, various degree-based topological indices have been introduced and attracted great interest in the past years \cite{MCSGDF18}. For example, \cite{BE98} proposed the general Randi\'{c} index, which includes the Randi\'{c} index and the second Zagreb index as special cases. \cite{ZT09,ZT10} defined the general sum-connectivity index, which includes the harmonic index as a special case \cite{Z12}. \cite{VG10,SV15} introduced the inverse sum indeg index to predict the total surface area of octane isomers. The reader is referred to \cite{G13} for more references.

One of the interesting research topics in the study of topological indices is to investigate their properties of random graphs. Recently, \cite{MMRS20,MMRS21} performed numeric and analytic  analyses of the Randi\'{c} index and the harmonic index in  the Erd\H{o}s-R\'{e}nyi random graph. Their simulation studies show that the indices are approximately equal to one half of the number of nodes, and the distributions of the indices are symmetric around their expectations. \cite{DMRSV18, DHIR20, LSG21} calculated the expectations of the Randi\'{c} index, the generalized Zagreb indices  and two modified Zagreb indices of the Erd\H{o}s-R\'{e}nyi random graph, respectively. \cite{Y22} derived the asymptotic limits of the Randi\'{c} index and the harmonic index of a heterogeneous random graph.

In this paper, we are interested in the asymptotic distribution of the degree-based topological index of a heterogeneous  Erd\H{o}s-R\'{e}nyi random graph. We show that after suitably centered and scaled, the degree-based topological index converges in distribution to the standard normal distribution. We apply our results to several well-known topological indices and observe that the general Randi\'{c} index exhibits a phase change phenomenon.

This paper is organized as follows. In Section \ref{smres}, we present the asymptotic distribution of the topological index of a heterogeneous random graph. In Section \ref{exmin}, we provide several examples. The proof is deferred to Section \ref{proof}.

\medskip

\noindent
{\bf Notation:} We adopt the  Bachmann–Landau notation throughout this paper. Let $a_n$  and $b_n$ be two positive sequences. Denote $a_n=\Theta(b_n)$ if $c_1b_n\leq a_n\leq c_2 b_n$ for some positive constants $c_1,c_2$. Denote  $a_n=\omega(b_n)$ if $\lim_{n\rightarrow\infty}\frac{a_n}{b_n}=\infty$. Denote $a_n=O(b_n)$ if $a_n\leq cb_n$ for some positive constants $c$. Denote $a_n=o(b_n)$ if $\lim_{n\rightarrow\infty}\frac{a_n}{b_n}=0$. Let $\mathcal{N}(0,1)$ be the standard normal distribution and $X_n$ be a sequence of random variables. Then $X_n\Rightarrow\mathcal{N}(0,1)$ means $X_n$ converges in distribution to the standard normal distribution as $n$ goes to infinity. Denote $X_n=O_P(a_n)$ if $\frac{X_n}{a_n}$ is bounded in probability. Denote $X_n=o_P(a_n)$ if $\frac{X_n}{a_n}$ converges to zero in probability as $n$ goes to infinity. Let $\mathbb{E}[X_n]$ and $Var(X_n)$ denote the expectation and variance of a random variable $X_n$ respectively. $\mathbb{P}[E]$ denote the probability of an event $E$. Let $f=f(x,y)$ be a function. For non-negative integers $s,t$, $f^{(s,t)}=f^{(s,t)}(x,y)$ denote the partial derivative $\frac{\partial^{s+t}f(x,y)}{\partial x^s\partial y^t}$. For convenience, sometimes we write $f_x=f^{(1,0)}$, $f_y=f^{(0,1)}$, $f_{xx}=f^{(2,0)}$,$f_{yy}=f^{(0,2)}$ and $f_{xy}=f^{(1,1)}$. $\exp(x)$ denote the exponential function $e^x$.
For positive integer $n$, denote $[n]=\{1,2,\dots,n\}$. Given a finite set $E$, $|E|$ represents the number of elements in $E$.

\section{Main results}\label{smres}

A graph consists of a set of nodes (vertices) and a set of edges. Given a positive integer $n$, an \textit{undirected} graph on $\mathcal{V}=[n]$ is a pair $\mathcal{G}=(\mathcal{V},\mathcal{E})$, where $\mathcal{E}$ is a collection of subsets of $\mathcal{V}$ such that $|e|=2$ for every $e\in\mathcal{E}$. Elements in $\mathcal{E}$ are called edges. A graph can be conveniently represented as an adjacency matrix $A$, where $A_{ij}=1$ if $\{i,j\}$ is an edge, $A_{ij}=0$ otherwise and $A_{ii}=0$. Since $\mathcal{G}$ is undirected, the adjacency matrix $A$ is symmetric. The degree $d_i$ of node $i$ is the number of edges connecting it, that is, $d_i=\sum_{j\neq i}A_{ij}$.  A graph is said to be random if $A_{ij} (1\leq i<j\leq n)$ are random. 

The degree-based topological index of a graph is defined as follows \cite{G13,AZG19}.

\begin{Definition}\label{def2}
The degree-based topological index of a graph $\mathcal{G}=(\mathcal{V},\mathcal{E})$  is defined as
\begin{equation}\label{topindex}
\mathcal{I}_n=\sum_{\{i,j\}\in \mathcal{E}}f(d_i,d_j),
\end{equation}
where $f(x,y)$ is a real function satisfying $f(x,y)=f(y,x)$. 
\end{Definition}
Many well-known topological indices can be expressed as (\ref{topindex}). For example, the Randi\'{c} index corresponds to $f(x,y)=(xy)^{-\frac{1}{2}}$  and the hyper-Zagreb index corresponds to  $f(x,y)=(x+y)^2$.  

\begin{Definition}\label{def1}
Let $\beta$ be a constant between zero and one, $n$ be an positive integer, and 
\[W=\{w_{ij}\in[\beta,1]| 1\leq i<j\leq n, w_{ji}=w_{ij} ,w_{ii}=0\}.\]
Define a heterogeneous random graph  $\mathcal{G}_n(\beta, W)$ as  
\[\mathbb{P}(A_{ij}=1)=p_nw_{ij},\]
where $A_{ij}$ $(1\leq i<j\leq n)$ are independent, $A_{ij}=A_{ji}$ and $p_n\in(0,1)$.
\end{Definition}

The expected degree of node $i$ in $\mathcal{G}_n(\beta, W)$ is  $\mathbb{E}[d_i]=\sum_{k\neq i}p_nw_{ik}$. In general, $\mathbb{E}[d_i]\neq \mathbb{E}[d_j]$ if $i\neq j$, that is, the expected degrees of nodes are not the same. Hence $\mathcal{G}_n(\beta, W)$ is a heterogeneous random graph.
When $w_{ij}=c$ $(1\leq i<j\leq n)$ for a constant $c\in(0,1)$, $\mathcal{G}_n(\beta, W)$ is the Erd\H{o}s-R\'{e}nyi random graph. It is homogeneous in the sense that nodes in it share the same expected degree.

Several recent works have studied the expectations of some special topological indices  of the Erd\H{o}s-R\'{e}nyi random graph \cite{MMRS20,MMRS21,DMRSV18, DHIR20, LSG21}. In this paper,  we derive the asymptotic distribution of the topological index $\mathcal{I}_n$  of  the heterogeneous random graph $\mathcal{G}_n(\beta, W)$. Our results can be applied to all the topological indices studied in \cite{MMRS20,MMRS21,DMRSV18, DHIR20, LSG21}.

Before presenting our results, we introduce several notations and assumptions. Let 
$w_{i(k)}=1+\sum_{l\notin\{ i,k\}}p_nw_{il}$
and 
\begin{equation}\label{sigman}
\sigma_n^2=\sum_{i< j}(a_{ij}+a_{ji})^2p_nw_{ij}(1-p_nw_{ij}),
\end{equation}
where
\[ a_{ij}=\frac{1}{2}f(w_{i(j)},w_{j(i)})+\frac{1}{2}\sum_{l\notin\{i,j\}}p_nw_{il}\big[f_x(w_{i(l)},w_{l(i)})+f_y(w_{i(l)},w_{l(i)})\big].\]

\begin{Assumption}\label{assump}
Let $k_0 (k_0\geq3)$, $s,t$ be non-negative integers. Suppose $np_n=\omega(\log(n))$ and the following conditions hold.\\
(C1)
\[
\sum_{i<j}(a_{ij}+a_{ji})^4p_n=o(\sigma_n^4).
\]
(C2). For all non-negative integers $s,t$ satisfying $s+t\leq k_0$,  there is some positive constant $C$  such that 
\[|f^{(s,t)}(x,y)|\leq (xy)^C.\]
(C3). Given $s,t$ satisfying $s+t=k_0$, $|f^{(s,t)}(x,y)|$ is monotone in $x$ and $y$.\\
(C4). For a large positive constant $M$ and positive sequences $a_n,b_n\in[(\log(np_n))^{-2},M]$, the following holds. For $s+t=k_0$,
\[n(np_n)^{\frac{k_0}{2}+1}|f^{(s,t)}(a_nnp_n,b_nnp_n)|=o\left(\sigma_n\right).\]
(C5). For $1\leq s+t\leq k_0-1$,
\[n(np_n)^{2(s+t)-1}|f^{(s,t)}(np_n,np_n)|^2=o\left(\sigma_n^2\right).\]
\end{Assumption}

Assumption \ref{assump} is not restrictive and many common degree-based topological indices satisfy this assumption as shown later. Under Assumption \ref{assump}, we derive the asymptotic distribution of the topological index $\mathcal{I}_n$ of $\mathcal{G}_n(\beta, W)$ as follows.

\begin{Theorem}\label{mainthm} Let $\mathcal{I}_n$ be the topological index defined in (\ref{topindex}) of the random graph
$\mathcal{G}_n(\beta, W)$  and $\sigma_n^2$ be defined in  (\ref{sigman}). Suppose Assumption \ref{assump} holds. 
Then
\begin{equation}
\frac{\mathcal{I}_n-\mathbb{E}[\mathcal{I}_n]}{\sigma_n}\Rightarrow \mathcal{N}(0,1),
\end{equation}
as $n$ goes to infinity. 
In addition, the expectation $\mathbb{E}[\mathcal{I}_n]$ has the following asymptotic expression
\begin{equation}\label{expecma}
\mathbb{E}[\mathcal{I}_n]=\left(1+O\left(\frac{1}{np_n}\right)\right)\sum_{1\leq i<j\leq n}p_nw_{ij}f(w_{i(j)},w_{j(i)}),
\end{equation}
where the error rate $\frac{1}{np_n}$ cannot be improved.
\end{Theorem}

Based on Theorem \ref{mainthm}, the degree-based topological index $\mathcal{I}_n$ (suitably centered and scaled) of the heterogeneous random graph
$\mathcal{G}_n(\beta, W)$ 
 converges in distribution to the standard normal distribution. As far as we know, this is the first theoretical result on limiting distribution of topological indices. Moreover, Theorem \ref{mainthm} provides the best approximation of the expectation of $\mathcal{I}_n$, in the sense that the error rate $\frac{1}{np_n}$ cannot be improved. For some special topological indices of the Erd\H{o}s-R\'{e}nyi random graph, it is possible to get an exact and compact expression of $\mathbb{E}[\mathcal{I}_n]$. For instance, \cite{LSG21} and \cite{DHIR20} obtained the exact expressions of the expectation of the hyper-Zagreb index and the forgotten topological
index of the Erd\H{o}s-R\'{e}nyi random graph respectively. However, for most topological indices, it seems impossible to get exact and closed-form expressions of the expectations \cite{DHIR20}. Our result (\ref{expecma}) provides an approximation of the expectations.

 The proof of Theorem \ref{mainthm} proceeds by decomposing $\mathcal{I}_n$ as a sum of leading term and remainder term, followed by finding the limiting distribution of the leading term and showing the remainder term is negligible. The condition $(C1)$ of Assumption \ref{assump} is used to prove the leading term converges in distribution to the standard normal distribution. The conditions $(C2)$--$(C5)$ are needed to bound the remainder term. The condition
$np_n=\omega(\log(n))$ requires the random graph to be relatively dense. This condition is common in theoretical network analysis. Assumption \ref{assump} is weak. We shall provide several examples of degree--based topological indices that satisfy this assumption in the subsequent section.

\section{Application to several topological indices}\label{exmin}

In this section, we apply Theorem \ref{mainthm} to several well-known topological indices of a special heterogeneous random graph.

Let $p_n=n^{-\alpha}$ for a constant $\alpha\in(0,1)$ and $w_{ij}=e^{-\frac{\kappa i}{n}}e^{-\frac{\kappa j}{n}}, (i\neq j)$ with non-negative constant $\kappa$. Then $e^{-2\kappa}\leq w_{ij}\leq 1$. Denote the corresponding random graph as $\mathcal{G}_n(\alpha,\kappa)$. When $\kappa=0$, $\mathcal{G}_n(\alpha,0)$ is the Erd\H{o}s-R\'{e}nyi random graph. In this case, we denote it as $\mathcal{G}_n(\alpha)$ for convenience. For $\kappa>0$, $\mathcal{G}_n(\alpha,\kappa)$ is heterogeneous. 

Denote $c(\kappa)=\frac{1-e^{-\kappa}}{\kappa}$ for $\kappa>0$ and  $c(0)=1$.  
Note that
\[\sum_{i=1}^ne^{-\frac{\kappa i}{n}}=\frac{e^{-\frac{\kappa}{n}}(1-e^{-\kappa})}{1-e^{-\frac{\kappa}{n}}}=nc(\kappa)+O(1),\ \ \ \ \ \kappa>0.\]
Then
\begin{eqnarray*}
w_{i(k)}&=&1+\sum_{l\notin\{ i,k\}}p_nw_{il}\\
&=&1+p_ne^{-\frac{\kappa i}{n}}\left(nc(\kappa)+O(1)-e^{-\frac{\kappa i}{n}}-e^{-\frac{\kappa k}{n}}\right)\\
&=&1+np_nc(\kappa)e^{-\frac{\kappa i}{n}}+O\left(p_n\right)\\
&=&np_nc(\kappa)e^{-\frac{\kappa i}{n}}+O\left(1\right), \ \   \kappa>0.
\end{eqnarray*}
When  $\kappa=0$, 
\[w_{i(k)}=1+\sum_{l\notin\{ i,k\}}p_n=1+(n-2)p_n.\]

\subsection{The general Randi\'{c} index}\label{geR}

The general Randi\'{c} index is a generalization of the well-known Randi\'{c} index and has been widely studied in literature \cite{MMRS20,MMRS21,DMRSV18}. Let $f(x,y)=(xy)^{\tau}$ for a non-zero constant $\tau$. The general Randi\'{c} index $\mathcal{I}_n$ is defined as
\begin{equation}\label{geRandi}
\mathcal{I}_n=\sum_{\{i,j\}\in\mathcal{E}}(d_id_j)^{\tau}.    
\end{equation}
When $\tau=-\frac{1}{2}$, $\mathcal{I}_n$ is the Randi\'{c} index.

Given non-negative integers $s,t$, straightforward computation yields
\[
f^{(s,t)}(x,y)=\left(\prod_{k=0}^{s-1}(\tau-k)\right)\left(\prod_{k=0}^{t-1}(\tau-k)\right)x^{\tau-s}y^{\tau-t}.
\]
Then
\begin{eqnarray*}
&&a_{ij}\\
&=&\frac{1}{2}\left[(np_n)^2c(\kappa)^2e^{-\frac{\kappa i}{n}}e^{-\frac{\kappa j}{n}}+O\left(np_n\right)\right]^{\tau}\\
&&+\frac{1}{2}\sum_{l\notin\{i,j\}}p_ne^{-\frac{\kappa i}{n}}e^{-\frac{\kappa l}{n}}\tau\left(np_nc(\kappa)e^{-\frac{\kappa i}{n}}+O\left(1\right)\right)^{\tau-1}\\
&&\times\left(np_nc(\kappa)e^{-\frac{\kappa l}{n}}+O\left(1\right)\right)^{\tau}\\
&&+\frac{1}{2}\sum_{l\notin\{i,j\}}p_ne^{-\frac{\kappa i}{n}}e^{-\frac{\kappa l}{n}}\tau\left(np_nc(\kappa)e^{-\frac{\kappa i}{n}}+O\left(1\right)\right)^{\tau}\\
&&\times\left(np_nc(\kappa)e^{-\frac{\kappa l}{n}}+O\left(1\right)\right)^{\tau-1}.
\end{eqnarray*}
Note that $e^{-\kappa}\leq e^{-\frac{\kappa i}{n}}\leq 1$.
If $\tau>0$, then $a_{ij}=\Theta((np_n)^{2\tau})$. In this case,
\begin{eqnarray*}
    \sigma_n^2=\sum_{i< j}(a_{ij}+a_{ji})^2p_n=\Theta\left(n(np_n)^{4\tau+1}\right).
\end{eqnarray*}
Hence
\begin{eqnarray*}
\frac{\sum_{i<j}(a_{ij}+a_{ji})^4p_n}{\sigma_n^4}
&=&O\left(\frac{n(np_n)^{8\tau+1}}{n^2(np_n)^{8\tau+2}}\right)=o(1).
\end{eqnarray*}

For $1\leq s+t$,
\begin{eqnarray*}
n(np_n)^{2(s+t)}|f^{(s,t)}(np_n,np_n)|^2&=&\Theta\left(n(np_n)^{4\tau}\right)=o\left(\sigma_n^2\right)
\end{eqnarray*}
Let $k_0=\max\left\{\lfloor1+\frac{1}{1-\alpha}\rfloor+1,3\right\}$. Then $k_0>1+\frac{1}{1-\alpha}$. For $s+t=k_0$, it is easy to verify that
\begin{eqnarray*}
n(np_n)^{\frac{k_0}{2}+1}|f^{(s,t)}(np_n,np_n)|=\Theta\left(n(np_n)^{1+2\tau-\frac{k_0}{2}}\right)
=o\left(\sigma_n\right).
\end{eqnarray*}
Then Assumption \ref{assump} holds and Theorem \ref{mainthm} applies. 

\begin{Corollary}
Let $\mathcal{I}_n$ be the general Randi\'{c} index defined in (\ref{geRandi}) of the random graph
$\mathcal{G}_n(\alpha,\kappa)$ and $\sigma_n^2$ be defined in  (\ref{sigman}). If $\tau>0$,  
then
\begin{equation}
\frac{\mathcal{I}_n-\mathbb{E}[\mathcal{I}_n]}{\sigma_n}\Rightarrow \mathcal{N}(0,1),
\end{equation}
as $n$ goes to infinity. 
In addition, the expectation $\mathbb{E}[\mathcal{I}_n]$ has the following asymptotic expression
\begin{equation}\label{expec}
\mathbb{E}[\mathcal{I}_n]=\left(1+O\left(\frac{1}{np_n}\right)\right)\sum_{1\leq i<j\leq n}p_nw_{ij}(w_{i(j)}w_{j(i)})^{\tau},
\end{equation}
where the error rate $\frac{1}{np_n}$ cannot be improved.
\end{Corollary}

When $\tau<0$, Assumption \ref{assump} may not hold. To see this, let $\kappa=0$. Then
\begin{eqnarray*}\nonumber
a_{ij}&=&\frac{(1+(n-2)p_n)^{2\tau}}{2}+\tau (1+(n-2)p_n)^{2\tau-1}(n-2)p_n\\
&=&(1+(n-2)p_n)^{2\tau}\left(\frac{1+2\tau}{2}-\frac{\tau}{1+(n-2)p_n}\right).
\end{eqnarray*}

If $\tau\neq-\frac{1}{2}$, then $a_{ij}=\Theta\left((np_n)^{2\tau}\right)$. Similar to the case $\tau>0$, Assumption \ref{assump} holds and Theorem \ref{mainthm} applies. 

If $\tau=-\frac{1}{2}$, then $a_{ij}=\Theta\left(\frac{1}{(np_n)^{2}}\right)$. In this case,
$\sigma_n^2=\Theta\left(\frac{n}{(np_n)^{3}}\right)$ and
\begin{eqnarray*}
n(np_n)^{2(s+t)-1}|f^{(s,t)}(np_n,np_n)|^2&=&\Theta\left(\frac{n}{(np_n)^{3}}\right).
\end{eqnarray*}
Clearly, $\frac{n}{(np_n)^{3}}\neq o(\sigma_n^2)$. Hence  Assumption \ref{assump} does not hold. We have to study this case separately.

When $\tau=-\frac{1}{2}$,  $\mathcal{I}_n$ is the well-known Randi\'{c} index.
The Randi\'{c} index is perhaps the first degree-based topological index introduced in \cite{R75}. Recently, \cite{MMRS20} performed simulation studies of the expectation and distribution of the Randi\'{c} index in Erd\H{o}s-R\'{e}nyi random graph. Its asymptotic limit was given in \cite{Y22}. Here we derive the asymptotic distribution of the Randi\'{c} index of the Erd\H{o}s-R\'{e}nyi random graph as follows.

\begin{Theorem}\label{geRandic}
Let $\mathcal{I}_n$ be the general Randi\'{c} index of the Erd\H{o}s-R\'{e}nyi random graph
$\mathcal{G}_n(\alpha)$ with constant $\alpha\in(0,1)$. Then
\begin{equation}
\frac{\mathcal{I}_n-\mathbb{E}[\mathcal{I}_n]}{\sigma_n}\Rightarrow \mathcal{N}(0,1),
\end{equation}
where 
\[\mathbb{E}[\mathcal{I}_n]=\frac{n(np_n)^{2\tau+1}}{2}\left(1+O\left(\frac{1}{np_n}\right)\right), \]
\[\sigma_n^2=\frac{n(n-1)(n-2)p_n^2(1-p_n)^2}{32(1+(n-2)p_n))^4}=\Theta\left(\frac{n}{(np_n)^2}\right),\]
 for $\tau=-\frac{1}{2}$ and
\[\sigma_n^2=\frac{(1+2\tau)^2}{2}n(n-1)p_n(1+(n-2)p_n)^{4\tau}(1+o(1))=\Theta\left(n(np_n)^{4\tau+1}\right),\]
 for $\tau\neq-\frac{1}{2}$.
\end{Theorem}

Based on Theorem \ref{geRandic}, we have an interesting finding. For $\tau\neq-\frac{1}{2}$, 
the order of $\sigma_n^2$ is $n(np_n)^{4\tau+1}$, while for $\tau=-\frac{1}{2}$, the order of $\sigma_n^2$ is $\frac{n}{(np_n)^{2}}$. Note that
\[n(np_n)^{4\left(-\frac{1}{2}\right)+1}=\frac{n}{np_n}=\omega\left(\frac{n}{(np_n)^{2}}\right). \]
Therefore, as a function of $\tau$, the order of $\sigma_n^2$  is continuous at $\tau\neq-\frac{1}{2}$, but discontinuous at $\tau=-\frac{1}{2}$. In this sense, the general Randi\'{c} index exhibits a phase change at $\tau=-\frac{1}{2}$.

\subsection{Hyper-Zagreb index}
 The Zagreb indices and its variants
 are frequently used to measure physical-chemical properties
of compounds. Recently, \cite{LSG21} studied the 
expectation of the hyper-Zagreb index of  the Erd\H{o}s-R\'{e}nyi random graph.
 Let $f(x,y)=(x+y)^2$.  The hyper-Zagreb index  $\mathcal{I}_n$ is defined as
\begin{equation}\label{hyperzagreb}
\mathcal{I}_n=\sum_{\{i,j\}\in\mathcal{E}}(d_i+d_j)^2.
\end{equation}
Clearly, $f^{(s,t)}(x,y)=0$ for $s+t\geq3$ and
\[f_x(x,y)=f_y(x,y)=2(x+y),\ \ f_{xx}(x,y)=f_{yy}(x,y)=f_{xy}(x,y)=2.\]
Straightforward calculation yields
\begin{eqnarray*}
&&\sum_{l\notin\{i,j\}}w_{il}\big[f_x(w_{i(l)},w_{l(i)})+f_y(w_{i(l)},w_{l(i)})\big]\\
&=&4e^{-\frac{\kappa i}{n}}\sum_{l\notin\{i,j\}}e^{-\frac{\kappa l}{n}}\left(np_nc(\kappa)e^{-\frac{\kappa i}{n}}+np_nc(\kappa)e^{-\frac{\kappa l}{n}}+O\left(1\right)\right)\\
&=&4e^{-\frac{\kappa i}{n}}\left[n^2p_nc(\kappa)^2e^{-\frac{\kappa i}{n}}+n^2p_nc(\kappa)c(2\kappa)+O(n)\right]\\
&=&4n^2p_nc(\kappa)^2e^{-\frac{2\kappa i}{n}}+4n^2p_nc(\kappa)c(2\kappa)e^{-\frac{\kappa i}{n}}+O(n).
\end{eqnarray*}
Since $e^{-\kappa}\leq e^{-\frac{\kappa i}{n}}\leq 1$ and $e^{-2\kappa}\leq e^{-\frac{2\kappa i}{n}}\leq 1$, then
\begin{eqnarray*}
a_{ij}&=&\frac{(np_n)^2c(\kappa)^2}{2}\left(e^{-\frac{\kappa i}{n}}+e^{-\frac{\kappa j}{n}}\right)^2\\
&&+2(np_n)^2c(\kappa)^2e^{-\frac{2\kappa i}{n}}+2(np_n)^2c(\kappa)c(2\kappa)e^{-\frac{\kappa i}{n}}+O(np_n)\\
&=&(np_n)^2\Bigg[\frac{5}{2}c(\kappa)^2e^{-\frac{2\kappa i}{n}}+\frac{1}{2}c(\kappa)^2e^{-\frac{2\kappa j}{n}}\\
&&+c(\kappa)^2e^{-\frac{\kappa i}{n}}e^{-\frac{\kappa j}{n}}+2c(\kappa)c(2\kappa)e^{-\frac{\kappa i}{n}}\Bigg]+O(np_n)\\
&=&\Theta\left((np_n)^2\right).
\end{eqnarray*}
Consequently, we have $\sigma_n^2=\Theta\left(n(np_n)^5\right)$ and
\begin{eqnarray*}
\frac{\sum_{i<j}(a_{ij}+a_{ji})^4p_n}{\sigma_n^4}
&=&O\left(\frac{n(np_n)^{9}}{n^2(np_n)^{10}}\right)=o(1).
\end{eqnarray*}
For $s+t\geq 3$, $f^{(s,t)}(x,y)=0$. For $s+t=2$, we have
\begin{eqnarray*}
n(np_n)^{4}=o\left(\sigma_n^2\right).
\end{eqnarray*}
Then Assumption \ref{assump} holds and Theorem \ref{mainthm} applies.

\begin{Corollary}\label{corhzage}
Let $\mathcal{I}_n$ be the hyper-Zagreb index defined in (\ref{hyperzagreb}) of the random graph
$\mathcal{G}_n(\alpha,\kappa)$  and $\sigma_n^2$ be defined in  (\ref{sigman}). Then
\begin{equation}
\frac{\mathcal{I}_n-\mathbb{E}[\mathcal{I}_n]}{\sigma_n}\Rightarrow \mathcal{N}(0,1),
\end{equation}
as $n$ goes to infinity. 
In addition, the expectation $\mathbb{E}[\mathcal{I}_n]$ has the following asymptotic expression
\begin{equation}\label{hexpec}
\mathbb{E}[\mathcal{I}_n]=\left(1+O\left(\frac{1}{np_n}\right)\right)\sum_{1\leq i<j\leq n}p_nw_{ij}(w_{i(j)}+w_{j(i)})^{2},
\end{equation}
where the error rate $\frac{1}{np_n}$ cannot be improved.
\end{Corollary}

By Theorem 1 in \cite{LSG21},
the 
expectation of the hyper-Zagreb index of  the Erd\H{o}s-R\'{e}nyi random graph $\mathcal{G}_n(\alpha)$ is equal to
\begin{equation}\label{hyperzagee}
    \mathbb{E}[\mathcal{I}_n]=n(n-1)(n-2)(2n-5)p_n^3+5n(n-1)(n-2)p_n^2+2n(n-1)p_n.
\end{equation}

By (\ref{hexpec}), we have
\begin{equation}\label{hypereq1}
\mathbb{E}[\mathcal{I}_n]=2n(n-1)(n-2)^2p_n^3\left(1+O\left(\frac{1}{np_n}\right)\right).
\end{equation}
Then our result (\ref{hypereq1}) is consistent with (\ref{hyperzagee}). In addition, by (\ref{hyperzagee}), the error rate $\frac{1}{np_n}$ cannot be improved, as stated in Corollary \ref{corhzage}.

\subsection{Forgotten topological index}
The forgotten topological index is another chemical index. \cite{DHIR20} studied the 
expectation of the forgotten topological index of  the Erd\H{o}s-R\'{e}nyi random graph.
 Let $f(x,y)=x^2+y^2$.  The  forgotten topological index $\mathcal{I}_n$ is defined as
\begin{equation}\label{fortopindex}
\mathcal{I}_n=\sum_{\{i,j\}\in\mathcal{E}}(d_i^2+d_j^2).
\end{equation}
Similar to the hyper-Zagreb index, it is easy to verify that Assumption \ref{assump} holds. Then Theorem \ref{mainthm} applies.

\begin{Corollary}\label{corfort}
Let $\mathcal{I}_n$ be the forgotten topological index defined in (\ref{fortopindex}) of the random graph
$\mathcal{G}_n(\alpha,\kappa)$  and $\sigma_n^2$ be defined in  (\ref{sigman}). Then
\begin{equation}
\frac{\mathcal{I}_n-\mathbb{E}[\mathcal{I}_n]}{\sigma_n}\Rightarrow \mathcal{N}(0,1),
\end{equation}
as $n$ goes to infinity. 
In addition, the expectation $\mathbb{E}[\mathcal{I}_n]$ has the following asymptotic expression
\begin{equation}\label{finexpec}
\mathbb{E}[\mathcal{I}_n]=\left(1+O\left(\frac{1}{np_n}\right)\right)\sum_{1\leq i<j\leq n}p_nw_{ij}(w_{i(j)}^{2}+w_{j(i)}^{2}),
\end{equation}
where the error rate $\frac{1}{np_n}$ cannot be improved.
\end{Corollary}

For  the Erd\H{o}s-R\'{e}nyi random graph $\mathcal{G}_n(\alpha)$, the expectation of  the forgotten topological index \cite{DHIR20} is equal to 
\begin{equation}\label{exaforin}
\mathbb{E}[\mathcal{I}_n]= n(n-1)(n-2)(n-3)p_n^3+3n(n-1)(n-2)p_n^2+n(n-1)p_n.
\end{equation}
By (\ref{finexpec}), we have
\begin{eqnarray*}
&&\mathbb{E}[\mathcal{I}_n]\\
&=&\left(1+O\left(\frac{1}{np_n}\right)\right)\frac{n(n-1)p_n}{2}\Big[(1+(n-2)p_n)^{2}+(1+(n-2)p_n)^{2}\Big]\\
&=&n(n-1)(n-2)^2p_n^3\left(1+O\left(\frac{1}{np_n}\right)\right).
\end{eqnarray*}
Hence our approximation (\ref{finexpec}) is consistent with (\ref{exaforin}). Moreover, by (\ref{exaforin}), the error rate $\frac{1}{np_n}$ cannot be improved as in Corollary \ref{corfort}.

\subsection{The inverse sum indeg index}
 The inverse sum indeg index is a significant predictor of total surface area of
octane isomers \cite{SV15,P18,VG10}.
Let $f(x,y)=\frac{xy}{x+y}$.  The inverse sum indeg index $\mathcal{I}_n$ is defined as
\begin{equation}\label{invsumin}
\mathcal{I}_n=\sum_{\{i,j\}\in\mathcal{E}}\frac{d_id_j}{d_i+d_j}.
\end{equation}
As far as we know, the inverse sum indeg index of random graph has not been studied in literature. Here we provide its asymptotic distribution and an approximation of its expectation.

Let $g(x,y)=xy$ and $h(x,y)=\frac{1}{x+y}$. For simplicity, we denote $f(x,y)$ as $f$. Then $f=gh$. Given positive integers $s,t$,  $h^{(s,t)}=\frac{c_{s,t}}{(x+y)^{1+s+t}}$, where $c_{s,t}$ is a constant dependent on $s,t$. Straightforward calculation yields
\begin{eqnarray*}
f^{(s,0)}&=&\sum_{r=0}^s\binom{s}{r}g^{(r,0)}h^{(s-r,0)}=gh^{(s,0)}+sg^{(1,0)}h^{(s-1,0)}\\
&=&\frac{c_{s,0}xy}{(x+y)^{1+s}}+\frac{sc_{s-1,0}y}{(x+y)^{s}},\\
f^{(0,t)}&=&\sum_{r=0}^t\binom{t}{r}g^{(0,r)}h^{(0,t-r)}=gh^{(0,t)}+tg^{(0,1)}h^{(0,t-1)}\\
&=&\frac{c_{0,t}xy}{(x+y)^{1+t}}+\frac{tc_{0,t-1}x}{(x+y)^{t}}.
\end{eqnarray*}
Further, for $s\geq1$ or $t\geq1$, we have
\begin{eqnarray*}
f^{(s,t)}&=&\sum_{r=0}^t\binom{t}{r}g^{(0,r)}h^{(s,t-r)}+s\sum_{r=0}^t\binom{t}{r}g^{(1,r)}h^{(s-1,t-r)}\\
&=&gh^{(s,t)}+tg^{(0,1)}h^{(s,t-1)}+sg^{(1,0)}h^{(s-1,t)}+stg^{(1,1)}h^{(s-1,t-1)}\\
&=&\frac{c_{s,t}xy}{(x+y)^{1+s+t}}+\frac{tc_{s,t-1}x}{(x+y)^{s+t}}+\frac{sc_{s-1,t}y}{(x+y)^{s+t}}+\frac{stc_{s-1,t-1}}{(x+y)^{s+t-1}}.
\end{eqnarray*}

Hence, for $s+t\geq1$, $|f^{(s,t)}(np_n,np_n)|$ can be bounded as follows
\[|f^{(s,t)}(np_n,np_n)|=O\left(\frac{1}{(np_n)^{s+t-1}}\right).\]
Note that
\[f_x(np_n,np_n)=f_y(np_n,np_n)=\frac{(np_n)^2}{(np_n+np_n)^2}=\frac{1}{4}.\]
Then $a_{ij}=\Theta\left(np_n\right)$, $\sigma_n^2=\Theta\left(n(np_n)^3\right)$ and
\begin{eqnarray*}
\frac{\sum_{i<j}(a_{ij}+a_{ji})^4p_n}{\sigma_n^4}
&=&O\left(\frac{n(np_n)^{5}}{n^2(np_n)^{6}}\right)=o(1).
\end{eqnarray*}

When $1\leq s+t$, we have
\begin{eqnarray*}
n(np_n)^{2(s+t)}|f^{(s,t)}(np_n,np_n)|^2&=&O\left(n(np_n)^{2}\right)=o\left(\sigma_n^2\right).
\end{eqnarray*}

Let $k_0=\max\{\lfloor1+\frac{1}{1-\alpha}\rfloor+1,3\}$. Then $k_0>1+\frac{1}{1-\alpha}$ and
\begin{eqnarray*}
n(np_n)^{\frac{k_0}{2}+1}|f^{(s,t)}(np_n,np_n)|=O\left(n(np_n)^{2-\frac{k_0}{2}}\right)
=o\left(\sigma_n\right).
\end{eqnarray*}
Assumption \ref{assump} holds. Then Theorem \ref{mainthm} applies.

\begin{Corollary}
Let $\mathcal{I}_n$ be the inverse sum indeg index defined in (\ref{invsumin}) of the random graph
$\mathcal{G}_n(\alpha,\kappa)$  and $\sigma_n^2$ be defined in  (\ref{sigman}). Then
\begin{equation}
\frac{\mathcal{I}_n-\mathbb{E}[\mathcal{I}_n]}{\sigma_n}\Rightarrow \mathcal{N}(0,1),
\end{equation}
as $n$ goes to infinity. 
In addition, the expectation $\mathbb{E}[\mathcal{I}_n]$ has the following asymptotic expression
\begin{equation}\label{expec}
\mathbb{E}[\mathcal{I}_n]=\left(1+O\left(\frac{1}{np_n}\right)\right)\sum_{1\leq i<j\leq n}p_nw_{ij}\frac{w_{i(j)}w_{j(i)}}{w_{i(j)}+w_{j(i)}},
\end{equation}
where the error rate $\frac{1}{np_n}$ cannot be improved.
\end{Corollary}

For  the Erd\H{o}s-R\'{e}nyi random graph $\mathcal{G}_n(\alpha)$, the expectation of the inverse sum indeg index can be expressed as 
\begin{equation*}
\mathbb{E}[\mathcal{I}_n]=\left(1+O\left(\frac{1}{np_n}\right)\right)\frac{n(n-1)(n-2)p_n^2}{4}.
\end{equation*}

\section{Proof of main results} \label{proof}

In this section, we provide detailed proofs of Theorem \ref{mainthm} and Theorem \ref{geRandic}. It is not easy to work on $\mathcal{I}_n$ defined in (\ref{topindex}) directly. In stead, 
we provide an alternative expression of $\mathcal{I}_n$ as follows
\begin{equation}\label{freq1}
\mathcal{I}_n=\sum_{1\leq i<j\leq n}A_{ij}f(d_{i(j)},d_{j(i)})=\frac{1}{2}\sum_{i\neq j}A_{ij}f(d_{i(j)},d_{j(i)}),
\end{equation}
where $d_{i(j)}=1+\sum_{l\notin\{ i,j\}}A_{il}$ and $d_{j(i)}=1+\sum_{l\notin\{ i,j\}}A_{jl}$. Note that $A_{ij}$, $d_{i(j)}$ and $d_{j(i)}$ are independent , $\mathbb{E}[d_{i(j)}]=w_{i(j)}$ and $\mathbb{E}[d_{j(i)}]=w_{j(i)}$. We will use these facts frequently in the proof.

\subsection{Lemmas}

Before proving Theorem \ref{mainthm} and Theorem \ref{geRandic},
we present two lemmas.
\begin{Lemma}\label{lem1}
 Let $\mathcal{G}_n(\beta, W)$ be defined in Definition \ref{def1}, $\delta_n=\left(\log(np_n)\right)^{-2}$ and $M$ be a constant greater than $\frac{e^2}{1-p_n\beta}$. 
For any $i\in[n]$, we have
\[\mathbb{P}(d_{i(j)}-1=k)\leq \exp(-np_n\beta(1+o(1))),\ \ \ \ k\leq \delta_nnp_n,\]
\[\mathbb{P}(d_{i(j)}-1=k)\leq \exp(-np_n\beta(1+o(1))),\ \ \ \ k\geq Mnp_n.\]
\end{Lemma}
\noindent
{\bf Proof of Lemma \ref{lem1}: } Given distinct indices $i,j$, let $\theta_{ij}=\{p_nw_{il}|l\in[n]\setminus\{i,j\}\}$. Then
$d_{i(j)}-1$ follows the Poisson-Binomial distribution $PB(\theta_{ij})$. Recall that $\beta\leq w_{ij}\leq1$. Then 
\begin{eqnarray}\nonumber
\mathbb{P}(d_{i(j)}-1=k)&=&\sum_{S\subset[n]\setminus\{i,j\},|S|=k}\prod_{l\in S}p_nw_{il}\prod_{l\in S^C\setminus\{i,j\}}(1-p_nw_{il})\\ \nonumber
&\leq&\sum_{S\subset[n]\setminus\{i,j\},|S|=k}\prod_{l\in S}p_n\prod_{l\in S^C\setminus\{i,j\}}(1-p_n\beta)\\ \label{rheq2}
&=&\binom{n-2}{k}p_n^k(1-p_n\beta)^{n-2-k}.
\end{eqnarray}
Note that $\binom{n-2}{k}\leq e^{k\log n-k\log k+k}$ and 
$(1-p_n\beta)^{n-2-k}=e^{(n-2-k)\log(1-p_n\beta)}$.
Then by (\ref{rheq2}) we get
\begin{eqnarray}\nonumber
&&\mathbb{P}(d_{i(j)}-1=k)\\ \label{rheq3}
&\leq&\exp\left(k\log (np_n)-k\log k+k+(n-2-k)\log(1-p_n\beta)\right).
\end{eqnarray}
Let $g(k)=k\log (np_n)-k\log k+k+(n-2-k)\log(1-p_n\beta)$. Considering $k$ as continuous variable, the derivative of $g(k)$ with respect to $k$ is equal to 
\[g^{\prime}(k)=\log\left(\frac{np_n}{1-p_n\beta}\right)-\log k.\]
Clearly, $g^{\prime}(k)>0$ for $k<\frac{np_n}{1-p_n\beta}$ and $g^{\prime}(k)<0$ for $k>\frac{np_n}{1-p_n\beta}$. Then $g(k)$ achieves its maximum at $k=\frac{np_n}{1-p_n\beta}$. For $k\leq \delta_nnp_n$,  $g(k)\leq g(\delta_nnp_n)$. Hence
\begin{eqnarray*}
&&\mathbb{P}(d_{i(j)}-1=k)\\
&\leq&\exp\left(\delta_nnp_n\log\frac{1}{\delta_n(1-p_n\beta)}+\delta_nnp_n+n\log(1-p_n\beta)\right)\\
&\leq& \exp\left(-np_n\beta(1+o(1))\right).
\end{eqnarray*}

Since $M>\frac{e^2}{1-p_n\beta}$, then $ Mnp_n\geq\frac{np_n}{1-p_n\beta}$.
 For $k\geq Mnp_n$,  $g(k)\leq g(Mnp_n)$. Hence
\begin{eqnarray}\nonumber
&&\mathbb{P}(d_{i(j)}-1=k)\\ \nonumber
&\leq&\exp\left(-M\log(M)np_n+Mnp_n+(n-2-Mnp_n)\log(1-p_n\beta)\right)\\ \nonumber
&=&\exp\left(-M(\log M-1)np_n+|n-2-Mnp_n|p_n\beta\right)\\ \nonumber
&\leq& \exp\left(-np_n\beta(1+o(1))\right).
\end{eqnarray}

\qed

\begin{Lemma}\label{lem2}
    Suppose ($C1$) of Assumption \ref{assump} holds. For the random graph $\mathcal{G}_n(\beta, W)$, we have
    \begin{equation*}
\frac{\sum_{i\neq j}a_{ij}(A_{ij}-p_nw_{ij})}{\sqrt{\sum_{i< j}(a_{ij}+a_{ji})^2p_nw_{ij}(1-p_nw_{ij})}}\Rightarrow\mathcal{N}(0,1),
    \end{equation*}
where 
 \[ a_{ij}=\frac{1}{2}f(w_{i(j)},w_{j(i)})+\frac{1}{2}\sum_{l\notin\{i,j\}}p_nw_{il}\big[f_x(w_{i(l)},w_{l(i)})+f_y(w_{i(l)},w_{l(i)})\big].\]
\end{Lemma}

\noindent
{\bf Proof of Lemma \ref{lem2}: } 
Let
\begin{equation*}
\mathcal{Z}_n=\frac{\sum_{i\neq j}a_{ij}(A_{ij}-p_nw_{ij})}{\sqrt{\sum_{i< j}(a_{ij}+a_{ji})^2p_nw_{ij}(1-p_nw_{ij})}}.
\end{equation*}
Then $\mathbb{E}[\mathcal{Z}_n]=0$.
Note that 
\[\sum_{i\neq j}a_{ij}(A_{ij}-p_nw_{ij})=\sum_{i<j}(a_{ij}+a_{ji})(A_{ij}-p_nw_{ij}).\]
Then 
    \begin{equation*}
\mathcal{Z}_n=\frac{\sum_{i<j}(a_{ij}+a_{ji})(A_{ij}-p_nw_{ij})}{\sqrt{\sum_{i< j}(a_{ij}+a_{ji})^2p_nw_{ij}(1-p_nw_{ij})}}.
    \end{equation*}
Recall that $A_{ij} (1\leq i<j\leq n)$ are independent. It is easy to verify that
\begin{eqnarray*}\nonumber
&&\mathbb{E}\left[\left(\sum_{i<j}(a_{ij}+a_{ji})(A_{ij}-p_nw_{ij})\right)^2\right]\\ \nonumber
&=&\sum_{i<j}(a_{ij}+a_{ji})^2\mathbb{E}[(A_{ij}-p_nw_{ij})^2]\\ \label{thueq4}
&=&\sum_{i<j}(a_{ij}+a_{ji})^2p_nw_{ij}(1-p_nw_{ij}).
\end{eqnarray*}
Hence $Var[\mathcal{Z}_n]=1$. 
Straightforward calculation yields
\[\mathbb{E}[(A_{ij}-p_nw_{ij})^4]=p_nw_{ij}[(1-p_nw_{ij})^4+p_n^3w_{ij}^3(1-p_nw_{ij})]\leq2p_nw_{ij}.\]
By ($C1$) of Assumption \ref{assump}, we have
\begin{eqnarray*}
    &&\frac{\sum_{i<j}(a_{ij}+a_{ji})^4\mathbb{E}[(A_{ij}-p_nw_{ij})^4]}{\left(\sum_{i< j}(a_{ij}+a_{ji})^2p_nw_{ij}(1-p_nw_{ij})\right)^2}\\
    &\leq&    \frac{2\sum_{i<j}(a_{ij}+a_{ji})^4p_nw_{ij}}{\left(\sum_{i< j}(a_{ij}+a_{ji})^2p_nw_{ij}(1-p_nw_{ij})\right)^2}\\
    &=&o(1).
\end{eqnarray*}
According to the  Lyapunov Central Limit Theorem, $\mathcal{Z}_n$ converges in distribution to the standard normal distribution.

\qed

\subsection{Proof of Theorem \ref{mainthm}}

Let $k_0$ be the integer in Assumption \ref{assump}. The proof strategy is as follows: firstly we use the Taylor expansion to expand the function $f(d_{i(j)},d_{j(i)})$ at $(w_{i(j)},w_{j(i)})$ to $k_0$-th order; then we write $\mathcal{I}_n$ as the sum of the leading term and remainder terms; finally we show the leading term (after suitable scaling) converges in distribution to the standard normal distribution and the remainder terms are negligible.

By the Taylor expansion, $f(d_{i(j)},d_{j(i)})$ can be decomposed as
\begin{eqnarray}\label{apeq1}
f(d_{i(j)},d_{j(i)})&=&M_{ij}+S_{ij}+T_{ij}+R_{ij},
\end{eqnarray}
where 
\begin{eqnarray}\nonumber
M_{ij}&=&f(w_{i(j)},w_{j(i)})+f_x(w_{i(j)},w_{j(i)})(d_{i(j)}-w_{i(j)})\\ \nonumber
&&+f_y(w_{i(j)},w_{j(i)})(d_{j(i)}-w_{j(i)}),\\ \nonumber
S_{ij}&=&\frac{1}{2}f_{xx}(w_{i(j)},w_{j(i)})(d_{i(j)}-w_{i(j)})^2\\ \nonumber
&&+\frac{1}{2}f_{yy}(w_{i(j)},w_{j(i)})(d_{j(i)}-w_{j(i)})^2\\ \nonumber
&&+f_{xy}(w_{i(j)},w_{j(i)})(d_{i(j)}-w_{i(j)})(d_{j(i)}-w_{j(i)}),\\ \nonumber
T_{ij}
&=&\sum_{k=3}^{k_0-1}\sum_{s+t=k}\frac{f^{(s,t)}(w_{i(j)},w_{j(i)})}{s!t!}(d_{i(j)}-w_{i(j)})^s(d_{j(i)}-w_{j(i)})^t,\\  \nonumber
R_{ij}&=&\sum_{s+t=k_0}\frac{f^{(s,t)}(X_{i(j)},X_{j(i)})}{s!t!}(d_{i(j)}-w_{i(j)})^s(d_{j(i)}-w_{j(i)})^t,
\end{eqnarray}
 $X_{i(j)}$ is between $d_{i(j)}$ and $w_{i(j)}$, and $X_{j(i)}$ is between $d_{j(i)}$ and $w_{j(i)}$.
By (\ref{freq1}), the topological index $\mathcal{I}_n$ is equal to
\begin{equation}\label{iineq1}
\mathcal{I}_n 
=\frac{1}{2}\sum_{ i\neq j}M_{ij}A_{ij}+\frac{1}{2}\sum_{ i\neq j}S_{ij}A_{ij}+\frac{1}{2}\sum_{ i\neq j}T_{ij}A_{ij}+\frac{1}{2}\sum_{ i\neq j}R_{ij}A_{ij}. 
\end{equation}
Then
\begin{eqnarray}\nonumber
&&\frac{\mathcal{I}_n-\mathbb{E}[\mathcal{I}_n]}{\sigma_n}\\ \nonumber
&=&\frac{\frac{1}{2}\sum_{ i\neq j}(M_{ij}A_{ij}-\mathbb{E}[M_{ij}A_{ij}])}{\sigma_n}+\frac{\frac{1}{2}\sum_{ i\neq j}(S_{ij}A_{ij}-\mathbb{E}[S_{ij}A_{ij}])}{\sigma_n}\\ \label{iiineq1}
&+&\frac{\frac{1}{2}\sum_{ i\neq j}(T_{ij}A_{ij}-\mathbb{E}[T_{ij}A_{ij}])}{\sigma_n}+\frac{\frac{1}{2}\sum_{ i\neq j}(R_{ij}A_{ij}-\mathbb{E}[R_{ij}A_{ij}])}{\sigma_n},
\end{eqnarray}
where $\sigma_n^2$ is defined in  (\ref{sigman}).

Next we show the first term in (\ref{iiineq1}) is leading term and the last three terms are negligible.

\subsubsection{Asymptotic normality of the first term in (\ref{iiineq1})}

To begin with, we study the first term of (\ref{iiineq1}).
Note that 
\begin{eqnarray}\nonumber
&&\sum_{ i\neq j}M_{ij}A_{ij}\\ \nonumber
&=&\sum_{ i\neq j}f(w_{i(j)},w_{j(i)})A_{ij}+
\sum_{ i\neq j}f_x(w_{i(j)},w_{j(i)})(d_{i(j)}-w_{i(j)})A_{ij}\\ \label{apeq2}
&&+\sum_{ i\neq j}f_y(w_{i(j)},w_{j(i)})(d_{j(i)}-w_{j(i)})A_{ij}.
\end{eqnarray}
Since $d_{i(j)}-w_{i(j)}=\sum_{l\notin\{ i,j\}}(A_{il}-p_nw_{il})$, then $(d_{i(j)}-w_{i(j)})$ does not contain $A_{ij}$. Hence, $(d_{i(j)}-w_{i(j)})$ and  $A_{ij}$ are independent. Similarly, $d_{j(i)}-w_{j(i)}$  and  $A_{ij}$ are independent. In addition, $\mathbb{E}[d_{i(j)}]=w_{i(j)}$ and $\mathbb{E}[d_{j(i)}]=w_{j(i)}$.
Then the expectation of $\sum_{ i\neq j}M_{ij}A_{ij}$ is equal to
\begin{equation}\label{WWeq2}
\mathbb{E}\left[\sum_{ i\neq j}M_{ij}A_{ij}\right]=\sum_{ i\neq j}f(w_{i(j)},w_{j(i)})p_nw_{ij}.
\end{equation}

The first term of (\ref{apeq2}) can be written as
\begin{eqnarray}\nonumber
\sum_{ i\neq j}f(w_{i(j)},w_{j(i)})A_{ij}&=&\sum_{ i\neq j}f(w_{i(j)},w_{j(i)})(A_{ij}-p_nw_{ij})\\ \label{apeq5}
&&+\sum_{ i\neq j}f(w_{i(j)},w_{j(i)})p_nw_{ij}.
\end{eqnarray}
Similarly, the second term of (\ref{apeq2}) is written as
\begin{eqnarray}\nonumber
&&\sum_{ i\neq j}f_x(w_{i(j)},w_{j(i)})(d_{i(j)}-w_{i(j)})A_{ij}\\ \nonumber
&=&\sum_{i\neq j}f_x(w_{i(j)},w_{j(i)})(d_{i(j)}-w_{i(j)})(A_{ij}-p_nw_{ij})\\ \nonumber
&&+\sum_{i\neq j}p_nw_{ij}f_x(w_{i(j)},w_{j(i)})(d_{i(j)}-w_{i(j)})\\ \nonumber
&=&\sum_{i\neq j\neq l}f_x(w_{i(j)},w_{j(i)})(A_{il}-p_nw_{il})(A_{ij}-p_nw_{ij})\\ \label{apeq3}
&&+\sum_{i\neq j}\left(\sum_{l\notin\{i,j\}}p_nw_{il}f_x(w_{i(l)},w_{l(i)})\right)(A_{ij}-p_nw_{ij}).
\end{eqnarray}
We will show the first term of (\ref{apeq3}) is of smaller order than the second term. To this end, we find the second moment of it. Recall that
if $\{i,j\}\neq\{s,t\}$,  $A_{ij}$ and $A_{st}$ are independent. Let $i,j,l$ be three arbitrary distinct indices and $i_1,j_1,l_1$ be another three arbitrary distinct indices. If $\{i,j,l\}\neq \{i_1,j_1,l_1\}$, then
\[\mathbb{E}[(A_{ij}-p_nw_{ij})(A_{il}-p_nw_{il})(A_{i_1j_1}-p_nw_{i_1j_1})(A_{i_1l_1}-p_nw_{i_1l_1})]=0.\]
When $\{i,j,l\}=\{i_1,j_1,l_1\}$, it is easy to get
\begin{eqnarray*}\nonumber
&&\mathbb{E}[(A_{ij}-p_nw_{ij})(A_{il}-p_nw_{il})(A_{i_1j_1}-p_nw_{i_1j_1})(A_{i_1l_1}-p_nw_{i_1l_1})]\\ \nonumber
&=&\mathbb{E}[(A_{ij}-p_nw_{ij})^2(A_{il}-p_nw_{il})^2]\\ \label{ineq2}
&=&p_nw_{ij}(1-p_nw_{ij})p_nw_{il}(1-p_nw_{il}).
\end{eqnarray*}

Then the second moment of the first term of (\ref{apeq3}) can be calculated as follows.
\begin{eqnarray*}\nonumber
&&\mathbb{E}\left[\left(\sum_{i\neq j\neq l}f_x(w_{i(j)},w_{j(i)})(A_{il}-p_nw_{il})(A_{ij}-p_nw_{ij})\right)^2\right]\\ \nonumber
&=&\sum_{\substack{i\neq j\neq l\\ i_1\neq j_1\neq l_1}}f_x(w_{i(j)},w_{j(i)})f_x(w_{i_1(j_1)},w_{j_1(i_1)})\\ \nonumber
&&\times\mathbb{E}[(A_{il}-p_nw_{il})(A_{ij}-p_nw_{ij})(A_{i_1l_1}-p_nw_{i_1l_1})(A_{i_1j_1}-p_nw_{i_1j_1})]\\ \nonumber
&=&\sum_{\substack{i\neq j\neq l}}f_x(w_{i(j)},w_{j(i)})^2\mathbb{E}[(A_{il}-p_nw_{il})^2(A_{ij}-p_nw_{ij})^2]\\ \nonumber
&=&\sum_{\substack{i\neq j\neq l}}f_x(w_{i(j)},w_{j(i)})^2p_nw_{il}(1-p_nw_{il})p_nw_{ij}(1-p_nw_{ij})\\
&=&\Theta\left(n^3p_n^2f_x(np_n,np_n)^2\right).
\end{eqnarray*}
By Markov's inequality, it follows that
\begin{eqnarray}\nonumber
   &&\sum_{i\neq j\neq l}f_x(w_{i(j)},w_{j(i)})(A_{il}-p_nw_{il})(A_{ij}-p_nw_{ij})\\ \label{ineq3}
&=&O_P\left(\sqrt{n^3p_n^2f_x(np_n,np_n)^2}\right).
\end{eqnarray}
Similarly, one has
\begin{eqnarray}\nonumber
&&\sum_{i\neq j\neq l}f_y(w_{i(j)},w_{j(i)})(A_{il}-p_nw_{il})(A_{ij}-p_nw_{ij})\\ \label{fyineq3}
&=&O_P\left(\sqrt{n^3p_n^2f_y(np_n,np_n)^2}\right).
\end{eqnarray}

Denote
  \[ a_{ij}=\frac{1}{2}f(w_{i(j)},w_{j(i)})+\frac{1}{2}\sum_{l\notin\{i,j\}}p_nw_{il}\big[f_x(w_{i(l)},w_{l(i)})+f_y(w_{i(l)},w_{l(i)})\big].\]
Then combining (\ref{apeq2})- (\ref{fyineq3}) yields
\begin{eqnarray}\nonumber
&&\frac{\frac{1}{2}\sum_{ i\neq j}(M_{ij}A_{ij}-\mathbb{E}[M_{ij}A_{ij}])}{\sigma_n}\\ \nonumber
&=&\frac{\sum_{ i\neq j}a_{ij}(A_{ij}-p_nw_{ij})}{\sigma_n}+O_P\left(\sqrt{\frac{n^3p_n^2f_x(np_n,np_n)^2}{\sigma_n^2}}\right).
\end{eqnarray}
By Lemma \ref{lem2} and $(C5)$ of Assumption \ref{assump} ( let $s+t=1$), we conclude that
\[\frac{\frac{1}{2}\sum_{ i\neq j}(M_{ij}A_{ij}-\mathbb{E}[M_{ij}A_{ij}])}{\sigma_n}\Rightarrow\mathcal{N}(0,1). \]
Then the proof is complete if the second, third and last term of (\ref{iiineq1}) converge to zero in probability.

\subsubsection{Bound the second term of (\ref{iiineq1})}

We prove the second term of (\ref{iiineq1})  is equal to $o_P(1)$. By the definition of $S_{ij}$, we have
\begin{eqnarray}\nonumber
   &&\sum_{ i\neq j}S_{ij}A_{ij}\\ \nonumber
   &=&\frac{1}{2}\sum_{ i\neq j}f_{xx}(w_{i(j)},w_{j(i)})(d_{i(j)}-w_{i(j)})^2A_{ij}\\ \nonumber
   &&+\frac{1}{2}\sum_{ i\neq j}f_{yy}(w_{i(j)},w_{j(i)})(d_{j(i)}-w_{j(i)})^2A_{ij}\\ \label{apeq10}
&&+\sum_{ i\neq j}f_{xy}(w_{i(j)},w_{j(i)})(d_{i(j)}-w_{i(j)})(d_{j(i)}-w_{j(i)})A_{ij}.
\end{eqnarray}
Then
\begin{eqnarray}\nonumber
&&\mathbb{E}\left[ \sum_{ i\neq j}S_{ij}A_{ij}\right]\\ \nonumber
  &=&\frac{1}{2}\sum_{ i\neq j}f_{xx}(w_{i(j)},w_{j(i)})\mathbb{E}\left[(d_{i(j)}-w_{i(j)})^2\right]p_nw_{ij}\\ \nonumber
  &&+\frac{1}{2}\sum_{ i\neq j}f_{yy}(w_{i(j)},w_{j(i)})\mathbb{E}\left[(d_{j(i)}-w_{j(i)})^2\right]p_nw_{ij}\\ \nonumber
  &=&\frac{1}{2}\sum_{ i\neq j\neq l}f_{xx}(w_{i(j)},w_{j(i)})p_nw_{il}(1-p_nw_{il})p_nw_{ij}\\ \label{apeq14} 
  &&+\frac{1}{2}\sum_{ i\neq j\neq l}f_{yy}(w_{i(j)},w_{j(i)})p_nw_{jl}(1-p_nw_{jl})p_nw_{ij}.
\end{eqnarray}

The first term of (\ref{apeq10}) can be expressed as
\begin{eqnarray}\nonumber
&&\sum_{ i\neq j}f_{xx}(w_{i(j)},w_{j(i)})(d_{i(j)}-w_{i(j)})^2A_{ij}\\ \nonumber
&=&\sum_{ i\neq j}f_{xx}(w_{i(j)},w_{j(i)})(d_{i(j)}-w_{i(j)})^2(A_{ij}-p_nw_{ij})\\ \label{apeq13}
&&+\sum_{ i\neq j}p_nw_{ij}f_{xx}(w_{i(j)},w_{j(i)})(d_{i(j)}-w_{i(j)})^2.
\end{eqnarray}
We will find an upper bound of (\ref{apeq13}).
Note that
\begin{eqnarray}\nonumber
&&\sum_{ i\neq j}f_{xx}(w_{i(j)},w_{j(i)})(d_{i(j)}-w_{i(j)})^2(A_{ij}-p_nw_{ij})\\ \nonumber
&=&\sum_{i\neq j}f_{xx}(w_{i(j)},w_{j(i)})(A_{ij}-p_nw_{ij})\sum_{\substack{s\neq t\\ s,t\notin\{i,j\}}}(A_{is}-p_nw_{is})(A_{it}-p_nw_{it})\\ \label{apeq11}
&&+\sum_{i\neq j}f_{xx}(w_{i(j)},w_{j(i)})(A_{ij}-p_nw_{ij})\sum_{s\notin\{i,j\}}(A_{is}-p_nw_{is})^2.
\end{eqnarray}
The second moment of the first term of (\ref{apeq11}) is equal to 
\begin{eqnarray}\nonumber
&&\mathbb{E}\Bigg[\sum_{i\neq j}f_{xx}(w_{i(j)},w_{j(i)})(A_{ij}-p_nw_{ij})\\ \nonumber
&&\times\sum_{\substack{s\neq t\\ s,t\notin\{i,j\}}}(A_{is}-p_nw_{is})(A_{it}-p_nw_{it})\Bigg]^2\\ \nonumber
&=&\sum_{\substack{i\neq j,s\neq t\\ s,t\notin\{i,j\}}}f_{xx}(w_{i(j)},w_{j(i)})^2\mathbb{E}\big[(A_{ij}-p_nw_{ij})^2\\ \nonumber
&&\times(A_{is}-p_nw_{is})^2(A_{it}-p_nw_{it})^2\big]\\ \nonumber
&=&\sum_{\substack{i\neq j,s\neq t\\ s,t\notin\{i,j\}}}f_{xx}(w_{i(j)},w_{j(i)})^2p_nw_{ij}(1-p_nw_{ij})\\  \nonumber
&&\times p_nw_{is}(1-p_nw_{is})p_nw_{it}(1-p_nw_{it})\\ \label{moneq7}
&=&\Theta\left(n^4p_n^3f_{xx}(np_n,np_n)^2\right).
\end{eqnarray}

The second moment of the second term of (\ref{apeq11}) is equal to 
\begin{eqnarray}\nonumber
&&\mathbb{E}\left[\sum_{i\neq j}f_{xx}(w_{i(j)},w_{j(i)})(A_{ij}-p_nw_{ij})\sum_{s\notin\{i,j\}}(A_{is}-p_nw_{is})^2\right]^2\\ \nonumber
&=&\sum_{i\neq j\neq s}f_{xx}(w_{i(j)},w_{j(i)})^2\mathbb{E}\left[(A_{ij}-p_nw_{ij})^2(A_{is}-p_nw_{is})^4\right]\\ \nonumber
&&+\sum_{i\neq j\neq s\neq t}f_{xx}(w_{i(j)},w_{j(i)})^2\mathbb{E}\big[(A_{ij}-p_nw_{ij})^2(A_{is}-p_nw_{is})^2\\  \nonumber
&&\times (A_{it}-p_nw_{it})^2\big]\\ \nonumber
&&+\sum_{i\neq j\neq s}f_{xx}(w_{i(j)},w_{j(i)})f_{xx}(w_{i(s)},w_{s(i)})\\  \nonumber
&&\times\mathbb{E}\left[(A_{ij}-p_nw_{ij})^3(A_{is}-p_nw_{is})^3\right]\\ \nonumber
&\leq&\sum_{i\neq j\neq s}f_{xx}(w_{i(j)},w_{j(i)})^2p_n^2+\sum_{i\neq j\neq s\neq t}f_{xx}(w_{i(j)},w_{j(i)})^2p_n^3\\ \nonumber
&&+\sum_{i\neq j\neq s}f_{xx}(w_{i(j)},w_{j(i)})f_{xx}(w_{i(s)},w_{s(i)})p_n^2\\ \label{moneq6}
&=&\Theta\left(n^4p_n^3f_{xx}(np_n,np_n)^2\right).
\end{eqnarray}

Now we consider the second term of (\ref{apeq13}). Note that
\begin{eqnarray}\nonumber
&&\sum_{ i\neq j}p_nw_{ij}f_{xx}(w_{i(j)},w_{j(i)})(d_{i(j)}-w_{i(j)})^2\\ \nonumber
&=&\sum_{ i\neq j}p_nw_{ij}f_{xx}(w_{i(j)},w_{j(i)})\sum_{\substack{s\neq t\\ s,t\notin\{i,j\}}}(A_{is}-p_nw_{is})(A_{it}-p_nw_{it})\\ \nonumber
&&+\sum_{ i\neq j}p_nw_{ij}f_{xx}(w_{i(j)},w_{j(i)})\\  \nonumber
&&\times\sum_{\substack{ s\notin\{i,j\}}}\big[(A_{is}-p_nw_{is})^2-\mathbb{E}[(A_{is}-p_nw_{is})^2]\big]\\ \label{mmoneq1}
&&+\sum_{ i\neq j}p_nw_{ij}f_{xx}(w_{i(j)},w_{j(i)})\sum_{\substack{ s\notin\{i,j\}}}\mathbb{E}\left[(A_{is}-p_nw_{is})^2\right].
\end{eqnarray}

The second moment of the first term of (\ref{mmoneq1}) is equal to

\begin{eqnarray}\nonumber
&&\mathbb{E}\left[\sum_{ i\neq j}p_nw_{ij}f_{xx}(w_{i(j)},w_{j(i)})\sum_{\substack{s\neq t\\ s,t\notin\{i,j\}}}(A_{is}-p_nw_{is})(A_{it}-p_nw_{it})\right]^2\\ \nonumber
&=&
\sum_{\substack{i\neq j,i\neq j_1, s\neq t\\ s,t\notin\{i,j,j_1\}}}p_n^2w_{ij}w_{ij_1}f_{xx}(w_{i(j)},w_{j(i)})f_{xx}(w_{i(j_1)},w_{j_1(i)})\\ \nonumber
&&\times\mathbb{E}\left[(A_{is}-p_nw_{is})^2(A_{it}-p_nw_{it})^2\right]\\ \nonumber
&&+\sum_{\substack{i\neq j, s\neq t\\ s,t\notin\{i,j\}}}p_n^2w_{ij}^2f_{xx}(w_{i(j)},w_{j(i)})^2\mathbb{E}\left[(A_{is}-p_nw_{is})^2(A_{it}-p_nw_{it})^2\right]\\ \nonumber
&=&\sum_{\substack{i\neq j,i\neq j_1, s\neq t\\ s,t\notin\{i,j,j_1\}}}p_n^2w_{ij}w_{ij_1}f_{xx}(w_{i(j)},w_{j(i)})f_{xx}(w_{i(j_1)},w_{j_1(i)})\\ \nonumber
&&\times p_nw_{is}(1-p_nw_{is})p_nw_{it}(1-p_nw_{it})\\ \nonumber
&&+\sum_{\substack{i\neq j,s\neq t\\ s,t\notin\{i,j\}}}p_n^2w_{ij}^2f_{xx}(w_{i(j)},w_{j(i)})^2p_nw_{is}(1-p_nw_{is})p_nw_{it}(1-p_nw_{it})\\ \label{mmoneq4}
&=&\Theta\left(n^5p_n^4f_{xx}(np_n,np_n)^2\right)
\end{eqnarray}

The second moment of the second term of (\ref{mmoneq1}) is equal to
\begin{eqnarray}\nonumber
&&\mathbb{E}\Bigg[\sum_{ i\neq j}p_nw_{ij}f_{xx}(w_{i(j)},w_{j(i)})\\  \nonumber
&&\times\sum_{\substack{ s\notin\{i,j\}}}\left[(A_{is}-p_nw_{is})^2-\mathbb{E}\left[(A_{is}-p_nw_{is})^2\right]\right]\Bigg]^2\\ \nonumber
&=&\sum_{\substack{i\neq j\neq s\\ i\neq j_1\neq s}}p_nw_{ij}f_{xx}(w_{i(j)},w_{j(i)})p_nw_{ij_1}f_{xx}(w_{i(j_1)},w_{j_1(i)})\\ \nonumber
&&\times\mathbb{E}\left[\left((A_{is}-p_nw_{is})^2-\mathbb{E}\left[(A_{is}-p_nw_{is})^2\right]\right)^2\right]\\ \nonumber
&&+\sum_{\substack{i\neq j\neq s}}p_n^2w_{ij}^2f_{xx}(w_{i(j)},w_{j(i)})^2\\  \nonumber
&&\times\mathbb{E}\left[\left((A_{is}-p_nw_{is})^2-\mathbb{E}\left[(A_{is}-p_nw_{is})^2\right]\right)^2\right]\\ \label{mmoneq3}
&=&\Theta\left(n^4p_n^3f_{xx}(np_n,np_n)^2\right).
\end{eqnarray}

By Markov's inequality and equations  (\ref{apeq13})- (\ref{mmoneq3}), we have
\begin{eqnarray}\nonumber
&&\sum_{ i\neq j}f_{xx}(w_{i(j)},w_{j(i)})\left[(d_{i(j)}-w_{i(j)})^2A_{ij}-\mathbb{E}\left[(d_{i(j)}-w_{i(j)})^2\right]A_{ij}\right]\\ \label{mmoneq2}
&=&O_P\left(\sqrt{n^5p_n^4f_{xx}(np_n,np_n)^2}\right).
\end{eqnarray}

Similarly, one has
\begin{eqnarray}\nonumber
&&\sum_{ i\neq j}f_{yy}(w_{i(j)},w_{j(i)})\left[(d_{i(j)}-w_{i(j)})^2A_{ij}-\mathbb{E}\left[(d_{i(j)}-w_{i(j)})^2\right]A_{ij}\right]\\ \label{mmoneq8}
&=&O_P\left(\sqrt{n^5p_n^4f_{yy}(np_n,np_n)^2}\right).
\end{eqnarray}

Now we bound the third term of (\ref{apeq10}). Note that
\begin{eqnarray}\nonumber
&&\sum_{i\neq j}f_{xy}(w_{i(j)},w_{j(i)})(d_{i(j)}-w_{i(j)})(d_{j(i)}-w_{j(i)})A_{ij}\\ \nonumber
&=&\sum_{i\neq j}f_{xy}(w_{i(j)},w_{j(i)})(d_{i(j)}-w_{i(j)})(d_{j(i)}-w_{j(i)})(A_{ij}-p_nw_{ij})\\ \nonumber
&&+\sum_{i\neq j}f_{xy}(w_{i(j)},w_{j(i)})(d_{i(j)}-w_{i(j)})(d_{j(i)}-w_{j(i)})p_nw_{ij}\\ \nonumber
&=&\sum_{\substack{i\neq j,s\neq j\\ t\neq i}}f_{xy}(w_{i(j)},w_{j(i)})(A_{ij}-p_nw_{ij})(A_{is}-p_nw_{is})(A_{jt}-p_nw_{jt})\\ \label{mmoneq9}
&&+\sum_{\substack{i\neq j,s\neq j\\ t\neq i}}f_{xy}(w_{i(j)},w_{j(i)})p_nw_{ij}(A_{is}-p_nw_{is})(A_{jt}-p_nw_{jt}).
\end{eqnarray}

The second moment of the first term of (\ref{mmoneq9}) is equal to 
\begin{eqnarray}\nonumber
&&\mathbb{E}\Bigg[\sum_{i\neq j,s\neq j,t\neq i}f_{xy}(w_{i(j)},w_{j(i)})(A_{ij}-p_nw_{ij})\\  \nonumber
&&\times(A_{is}-p_nw_{is})(A_{jt}-p_nw_{jt})\Bigg]^2\\ \nonumber
&=&\sum_{\substack{i\neq j,s\neq j,\\t\neq i, s\neq t}}f_{xy}(w_{i(j)},w_{j(i)})^2\mathbb{E}\big[(A_{ij}-p_nw_{ij})^2\\  \nonumber
&&\times(A_{is}-p_nw_{is})^2(A_{jt}-p_nw_{jt})^2\big]\\ \nonumber
&&+\sum_{i\neq j,s\neq j,}f_{xy}(w_{i(j)},w_{j(i)})^2\mathbb{E}\big[(A_{ij}-p_nw_{ij})^2\\  \nonumber
&&\times(A_{is}-p_nw_{is})^2(A_{js}-p_nw_{js})^2\big]\\ \nonumber
&=&\sum_{\substack{i\neq j,s\neq j,\\t\neq i,s\neq t}}f_{xy}(w_{i(j)},w_{j(i)})^2p_nw_{ij}(1-p_nw_{ij})\\  \nonumber
&&\times p_nw_{is}(1-p_nw_{is})p_nw_{jt}(1-p_nw_{jt})\\ \nonumber
&&+\sum_{\substack{i\neq j,\\ s\neq j}}f_{xy}(w_{i(j)},w_{j(i)})^2p_nw_{ij}(1-p_nw_{ij})\\  \nonumber
&&\times p_nw_{is}(1-p_nw_{is})p_nw_{js}(1-p_nw_{js})\\ \label{caneq2}
&=&\Theta\left(n^4p_n^3f_{xy}(np_n,np_n)^2\right).
\end{eqnarray}

The second moment of the second term of (\ref{mmoneq9}) is equal to 
\begin{eqnarray}\nonumber
  &&\mathbb{E} \left[\sum_{i\neq j,s\neq j,t\neq i}f_{xy}(w_{i(j)},w_{j(i)})p_nw_{ij}(A_{is}-p_nw_{is})(A_{jt}-p_nw_{jt})\right]^2\\ \nonumber
  &=&\sum_{\substack{i\neq j,s\neq j\\ t\neq i,s\neq t}}f_{xy}(w_{i(j)},w_{j(i)})^2p_n^2w_{ij}^2\mathbb{E} \left[(A_{is}-p_nw_{is})^2(A_{jt}-p_nw_{jt})^2\right]\\ \nonumber
  &&+\sum_{i\neq j,s\neq j}f_{xy}(w_{i(j)},w_{j(i)})^2p_n^2w_{ij}^2\mathbb{E} \left[(A_{is}-p_nw_{is})^2(A_{js}-p_nw_{js})^2\right]\\ \nonumber
  &=&\sum_{\substack{i\neq j,s\neq j\\ t\neq i,s\neq t}}f_{xy}(w_{i(j)},w_{j(i)})^2p_n^2w_{ij}^2p_nw_{is}(1-p_nw_{is})p_nw_{jt}(1-p_nw_{jt})\\ \nonumber
  &&+\sum_{i\neq j,s}f_{xy}(w_{i(j)},w_{j(i)})^2p_n^2w_{ij}^2p_nw_{is}(1-p_nw_{is})p_nw_{js}(1-p_nw_{js})\\ \label{caneq1}
  &=&\Theta\left(n^4p_n^4f_{xy}(np_n,np_n)^2\right).
\end{eqnarray}

Combining (\ref{mmoneq9}), (\ref{caneq2}) and (\ref{caneq1}) yields
\begin{eqnarray}\nonumber
&&\sum_{ i\neq j}f_{xy}(w_{i(j)},w_{j(i)})(d_{i(j)}-w_{i(j)})(d_{j(i)}-w_{j(i)})A_{ij}\\ \label{caneq3}
&=&O_P\left(\sqrt{n^4p_n^3f_{xy}(np_n,np_n)^2}\right).
\end{eqnarray}

By (\ref{mmoneq2}), (\ref{mmoneq8}), (\ref{caneq3}) and $(C5)$ of Assumption \ref{assump} (let $s+t=2$), we get 
\begin{eqnarray}\nonumber
&&\sum_{i\neq j}(S_{ij}A_{ij}-\mathbb{E}\left[S_{ij}A_{ij}\right])=\\ \nonumber
&&O_P\left(\sqrt{n^5p_n^4[f_{xx}(np_n,np_n)^2+f_{yy}(np_n,np_n)^2]+n^4p_n^3f_{xy}(np_n,np_n)^2}\right)\\ \label{caneq4}
&&=o_P(\sigma_n).
\end{eqnarray}
Hence the second term of (\ref{iiineq1})  is equal to $o_P(1)$.

\subsubsection{Bound the third term of (\ref{iiineq1})}

Now we prove the third term of (\ref{iiineq1}) converges in probability to zero. This is the most complex part of the proof. Note that 
\begin{eqnarray}\nonumber
&&\sum_{i\neq j}T_{ij}A_{ij}\\ \nonumber
&=&\sum_{k=3}^{k_0-1}\sum_{s+t=k}\sum_{i\neq j}\frac{f^{(s,t)}(w_{i(j)},w_{j(i)})}{s!t!}(d_{i(j)}-w_{i(j)})^s(d_{j(i)}-w_{j(i)})^tA_{ij}\\ \nonumber
&=&\sum_{k=3}^{k_0-1}\sum_{s+t=k}\sum_{i\neq j}\frac{f^{(s,t)}(w_{i(j)},w_{j(i)})}{s!t!}(d_{i(j)}-w_{i(j)})^s(d_{j(i)}-w_{j(i)})^t\\  \nonumber
&&\times (A_{ij}-p_nw_{ij})\\   \nonumber
&&+\sum_{k=3}^{k_0-1}\sum_{s+t=k}\sum_{i\neq j}\frac{f^{(s,t)}(w_{i(j)},w_{j(i)})}{s!t!}(d_{i(j)}-w_{i(j)})^s\\ \label{tueeq1}
&&\times (d_{j(i)}-w_{j(i)})^tp_nw_{ij}.
\end{eqnarray}
Next we bound the second moment of the first term of (\ref{tueeq1}) and the variance of the second term. Since $k_0$ is a fixed finite integer, the quantities $s!,t!$ in (\ref{tueeq1}) are finite. We will ignore them in the subsequent analysis for simplicity. Given finite integer $k_0\geq 4$, there are finitely many non-negative integers $s,t$ such that $s+t=k$ for any $k=3,4,\dots,k_0-1$. Hence we only need to bound the second moment of 
\begin{eqnarray}\label{tueeq2}
\sum_{i\neq j}f^{(s,t)}(w_{i(j)},w_{j(i)})(d_{i(j)}-w_{i(j)})^s(d_{j(i)}-w_{j(i)})^t(A_{ij}-p_nw_{ij}),
\end{eqnarray}
and the variance of 
\begin{eqnarray}\label{tueeq3}
\sum_{i\neq j}f^{(s,t)}(w_{i(j)},w_{j(i)})(d_{i(j)}-w_{i(j)})^s(d_{j(i)}-w_{j(i)})^tp_nw_{ij},
\end{eqnarray}
where $s,t$ are given non-negative integers with $s+t=k$ for $k=3,4,\dots,k_0-1$.

We consider the variance of (\ref{tueeq3}) first. Fix integer $k\in\{3,4,\dots,k_0-1\}$ and integers $s,t\in\{0,1,2,\dots,k\}$ satisfying $s+t=k$. For positive integers $r\leq s$ and $v\leq t$, let $\lambda_{r;1},\lambda_{r;2},\dots,\lambda_{r;r}$, $\gamma_{v;1},\gamma_{v;2},\dots,\gamma_{v;v}$ be positive integers such that $\lambda_{r;1}+\lambda_{r;2}+\dots+\lambda_{r;r}=s$ and $\gamma_{v;1}+\gamma_{v;2}+\dots+\gamma_{v;v}=t$.
Given indices $i,j$, we have
\begin{eqnarray}\nonumber
(d_{i(j)}-w_{i(j)})^s&=&\sum_{j_1,j_2,\dots,j_s\notin\{i,j\}}\prod_{l=1}^s(A_{ij_l}-p_nw_{ij_l})\\ \nonumber
&=&\sum_{r=1}^s\sum_{\substack{j_1,j_2,\dots,j_r\notin\{i,j\}\\ j_1\neq j_2\neq \dots\neq j_r}}\prod_{l=1}^r(A_{ij_l}-p_nw_{ij_l})^{\lambda_{r;l}},\\ \nonumber
(d_{j(i)}-w_{j(i)})^t&=&\sum_{i_1,i_2,\dots,i_t\notin\{i,j\}}\prod_{l=1}^t(A_{ji_l}-p_nw_{ji_l})\\ \label{mmmeq1}
&=&\sum_{v=1}^t\sum_{\substack{i_1,i_2,\dots,i_v\notin\{i,j\}\\i_1\neq i_2\neq \dots\neq i_v}}\prod_{l=1}^v(A_{ji_l}-p_nw_{ji_l})^{\gamma_{v;l}}.
\end{eqnarray}
Then (\ref{tueeq3}) can be written as
\begin{eqnarray}\nonumber
   &&\sum_{i\neq j}p_nw_{ij}f^{(s,t)}(w_{i(j)},w_{j(i)})(d_{i(j)}-w_{i(j)})^s(d_{j(i)}-w_{j(i)})^t=\sum_{r=1}^s\sum_{v=1}^tV_{rv},
\end{eqnarray}
where
\begin{eqnarray}\nonumber
V_{rv}&=&\sum_{\substack{i\neq j\\ j_1,\dots,j_r\notin\{i,j\}\\ j_1\neq j_2\neq \dots\neq j_r\\
i_1,\dots,i_v\notin\{i,j\}\\i_1\neq i_2\neq \dots\neq i_v}}p_nw_{ij}f^{(s,t)}(w_{i(j)},w_{j(i)})\prod_{l=1}^r(A_{ij_l}-p_nw_{ij_l})^{\lambda_{r;l}}\\ \label{tueeq4}
&&\times\prod_{m=1}^v(A_{ji_m}-p_nw_{ji_m})^{\gamma_{v;m}}. 
\end{eqnarray}

Note that
\begin{eqnarray}  \nonumber
Var\left(\sum_{r=1}^s\sum_{v=1}^tV_{rv}\right)&=&\sum_{r=1}^s\sum_{v=1}^t\sum_{r_1=1}^s\sum_{v_1=1}^tCov(V_{rv},V_{r_1v_1})\\ \nonumber
&\leq& \sum_{r=1}^s\sum_{v=1}^t\sum_{r_1=1}^s\sum_{v_1=1}^t\Big(Var(V_{rv})+Var(V_{r_1v_1})\Big),
\end{eqnarray}
and $s,t$ are finite non-negative integers, we only need to bound $Var(V_{rv})$ for each given $r,v$. Fix $r\in\{1,2,\dots,s\}$ and $v\in\{1,2,\dots,t\}$. There are two cases: (I) there exists $l_0\in\{1,2,\dots,r\}$ or $m_0\in\{1,2,\dots,v\}$ such that  $\lambda_{r;l_0}=1$ or $\gamma_{v;m_0}=1$ ; (II) $\lambda_{r;l}\geq2$ for all $l\in\{1,2,\dots,r\}$ and $\gamma_{v;m}\geq2$ for all $m\in\{1,2,\dots,v\}$.

We study case (I) first. Suppose there are some $\lambda_{r;l}$ or $\gamma_{v;m}$ which are equal to one. Without loss of generality, let $\lambda_{r;1}=\lambda_{r;2}=\dots=\lambda_{r;r_0}=1$ and $\lambda_{r;l}\geq2$ for $l\in\{r_0+1,\dots,r\}$, $\gamma_{v;1}=\gamma_{v;2}=\dots=\gamma_{v;v_0}=1$ and $\gamma_{v;l}\geq2$ for $l\in\{v_0+1,\dots,v\}$. Here, either $r_0\geq1$ or $v_0\geq1$. Without loss of generality, let $r_0\geq1$. In this case,
\begin{eqnarray} \nonumber
&&\prod_{l=1}^r(A_{ij_l}-p_nw_{ij_l})^{\lambda_{r;l}}\\ \label{tueeq5}
&=&\left(\prod_{l=1}^{r_0}(A_{ij_l}-p_nw_{ij_l})\right)\left(\prod_{l=r_0+1}^r(A_{ij_l}-p_nw_{ij_l})^{\lambda_{r;l}}\right),
\end{eqnarray}
\begin{eqnarray} \nonumber
&&
\prod_{m=1}^v(A_{ji_m}-p_nw_{ji_m})^{\gamma_{v;m}}\\ \label{tueeq6}
&=&\left(\prod_{m=1}^{v_0}(A_{ji_m}-p_nw_{ji_m})\right)\left(\prod_{m=v_0+1}^{v}(A_{ji_m}-p_nw_{ji_m})^{\gamma_{v;m}}\right),
\end{eqnarray}
and $ \mathbb{E}[V_{rv}]=0$. Then $Var(V_{rv})=\mathbb{E}[V_{rv}^2]$. For convenience, denote $\bar{A}_{ij}=A_{ij}-p_nw_{ij}$. By (\ref{tueeq4}),  (\ref{tueeq5}) and  (\ref{tueeq6}), we have
\begin{eqnarray}\nonumber
&&\mathbb{E}[V_{rv}^2]\\ \nonumber
 &=&\sum_{\substack{i\neq j\\ j_1,\dots,j_r\notin\{i,j\}\\ j_1\neq j_2\neq \dots\neq j_r\\
i_1,\dots,i_v\notin\{i,j\}\\i_1\neq i_2\neq \dots\neq i_v}}\sum_{\substack{i^{\prime}\neq j^{\prime}\\j_1^{\prime},\dots,j_r^{\prime}\notin\{i^{\prime},j^{\prime}\}\\ j_1^{\prime}\neq j_2^{\prime}\neq \dots\neq j_r^{\prime}\\
i_1^{\prime},\dots,i_v^{\prime}\notin\{i^{\prime},j^{\prime}\}\\ i_1^{\prime}\neq i_2^{\prime}\neq \dots\neq i_v^{\prime}}}p_nw_{ij}f^{(s,t)}(w_{i(j)},w_{i(j)})\\ \nonumber
&&\times p_nw_{i^{\prime}j^{\prime}}f^{(s,t)}(w_{i^{\prime}(j^{\prime})},w_{j^{\prime}(i^{\prime})})\\ \nonumber
&&\times\mathbb{E}\Bigg[\left(\prod_{l=1}^{r_0}\bar{A}_{ij_l}\right)\left(\prod_{l=r_0+1}^r\bar{A}_{ij_l}^{\lambda_{r;l}}\right)\left(\prod_{l=1}^{r_0}\bar{A}_{i^{\prime}j_l^{\prime}}\right)\left(\prod_{l=r_0+1}^r\bar{A}_{i^{\prime}j_l^{\prime}}^{\lambda_{r;l}}\right)\\ \nonumber
&&\times\left(\prod_{m=1}^{v_0}\bar{A}_{ji_m}\right)\left(\prod_{m=v_0+1}^{v}\bar{A}_{ji_m}^{\gamma_{v;m}}\right)\left(\prod_{m=1}^{v_0}\bar{A}_{j^{\prime}i_m^{\prime}}\right)\left(\prod_{m=v_0+1}^{v}\bar{A}_{j^{\prime}i_m^{\prime}}^{\gamma_{v;m}}\right)\Bigg]. \\ \label{tueeq7} 
\end{eqnarray}

Next we find an upper bound of (\ref{tueeq7}). Recall that $i\neq j$ and  $i^{\prime}\neq j^{\prime}$. We shall decompose the summation in (\ref{tueeq7}) into six cases: $i\neq i^{\prime}$ and $j=j^{\prime}$; $i= i^{\prime}$ and $j\neq j^{\prime}$; $i\neq j^{\prime}$ and $j=i^{\prime}$; $i= j^{\prime}$ and $j\neq i^{\prime}$; $\{i,j\}=\{i^{\prime},j^{\prime}\}$; $\{i,j\}\cap\{i^{\prime},j^{\prime}\}=\emptyset$. For convenience, denote the expectation in (\ref{tueeq7}) as $E$.

Firstly, we consider the case $\{i,j\}=\{i^{\prime},j^{\prime}\}$. There are two scenarios:  (i) $i=i^{\prime}$ and $j=j^{\prime}$, (ii) $i=j^{\prime}$ and $j=i^{\prime}$.

Consider (i) first. In this case, $A_{ij_l}$, $A_{ij_l^{\prime}}$ are independent of $A_{ji_m}$, $A_{ji_m^{\prime}}$. Then
\begin{eqnarray}\nonumber
E&=&\mathbb{E}\Bigg[\left(\prod_{l=1}^{r_0}\bar{A}_{ij_l}\bar{A}_{ij_l^{\prime}}\right)\left(\prod_{l=r_0+1}^r\bar{A}_{ij_l}^{\lambda_{r;l}}\bar{A}_{ij_l^{\prime}}^{\lambda_{r;l}}\right)\Bigg]\\ \label{tueeq8}
&&\times\mathbb{E}\Bigg[\left(\prod_{m=1}^{v_0}\bar{A}_{ji_m}\bar{A}_{ji_m^{\prime}}\right)\left(\prod_{m=v_0+1}^{v}\bar{A}_{ji_m}^{\gamma_{v;m}}\bar{A}_{ji_m^{\prime}}^{\gamma_{v;m}}\right)\Bigg].
\end{eqnarray}
Recall that $j_1,j_2,\dots,j_r$ are mutually distinct and $j_1^{\prime},j_2^{\prime},\dots,j_r^{\prime}$ are mutually distinct. Moreover, $\mathbb{E}[\bar{A}_{ij_{l}}]=0$ for all $l=1,2,\dots,r$.
If there exists an index $j_{l_1}$ with $1\leq l_1\leq r_0$ such that $j_{l_1}\notin\{j_1^{\prime},j_2^{\prime},\dots,j_r^{\prime}\}$, then 
\begin{eqnarray}\nonumber
&&\mathbb{E}\Bigg[\left(\prod_{l=1}^{r_0}\bar{A}_{ij_l}\bar{A}_{ij_l^{\prime}}\right)\left(\prod_{l=r_0+1}^r\bar{A}_{ij_l}^{\lambda_{r;l}}\bar{A}_{ij_l^{\prime}}^{\lambda_{r;l}}\right)\Bigg]\\  \nonumber
&=&\mathbb{E}[\bar{A}_{ij_{l_1}}]\mathbb{E}\Bigg[\left(\prod_{l=1,j_l\neq j_{l_1}}^{r_0}\bar{A}_{ij_l}\right)\left(\prod_{l=1}^{r_0}\bar{A}_{ij_l^{\prime}}\right)\left(\prod_{l=r_0+1}^r\bar{A}_{ij_l}^{\lambda_{r;l}}\bar{A}_{ij_l^{\prime}}^{\lambda_{r;l}}\right)\Bigg]\\ \label{suneq1}
&=&0.
\end{eqnarray}
Hence $E=0$ by (\ref{tueeq8}). Similarly, 
if there exists an index $l_1$ with $1\leq l_1\leq r_0$ such that $j_{l_1}^{\prime}\notin\{j_1,j_2,\dots,j_r\}$, then  $E=0$. In addition, if there is an index $m_1$ with $1\leq m_1\leq v_0$ such that $i_{m_1}\notin\{i_1^{\prime},i_2^{\prime},\dots,i_v^{\prime}\}$ or $i_{m_1}^{\prime}\notin\{i_1,i_2,\dots,i_v\}$, then $E=0$.  Consequently,  $E\neq 0$ implies the following
\[\{j_1,j_2,\dots,j_{r_0}\}\subset\{j_1^{\prime},j_2^{\prime},\dots,j_r^{\prime}\}, \ \ \ \ \{j_1^{\prime},j_2^{\prime},\dots,j_{r_0}^{\prime}\}\subset\{j_1,j_2,\dots,j_r\}, \]
\[\{i_1,i_2,\dots,i_{v_0}\}\subset\{i_1^{\prime},i_2^{\prime},\dots,i_v^{\prime}\}, \ \ \ \ \{i_1^{\prime},i_2^{\prime},\dots,i_{v_0}^{\prime}\}\subset\{i_1,i_2,\dots,i_v\}. \]
Without loss of generality,
suppose 
\begin{equation}\label{tueeq11}\{j_1,j_2,\dots,j_{r_1}\}=\{j_1^{\prime},j_2^{\prime},\dots,j_{r_1}^{\prime}\},\ \ \ \ \ \ \ \{j_{r_1+1},\dots,j_r\}\cap\{j_{r_1+1}^{\prime},\dots,j_{r}^{\prime}\}=\emptyset,
\end{equation}
\begin{equation}\label{tueeq12}\{i_1,i_2,\dots,i_{v_1}\}=\{i_1^{\prime},i_2^{\prime},\dots,i_{v_1}^{\prime}\},\ \ \ \ \{i_{v_1+1},\dots,i_v\}\cap\{i_{v_1+1}^{\prime},\dots,i_{v}^{\prime}\}=\emptyset,
\end{equation}
for some $r_1$ ($r_0\leq r_1\leq r$) and $v_1$ ($v_0\leq v_1\leq v$). There are at most $n^{2+2r-r_1+2v-v_1}$ possible choices for the indices $i,j,i_1,\dots,i_v,j_1,\dots,j_r$, 
$i^{\prime}$, $j^{\prime}$, $i_1^{\prime}$, $\dots$, $i_v^{\prime}$, $j_1^{\prime}$,$\dots$,$j_r^{\prime}$ satisfying $i=i^{\prime}$, 
$j=j^{\prime}$, (\ref{tueeq11}) and (\ref{tueeq12}). Let $\sigma_1$ be a permutation  of $\{1,2,\dots,r_1\}$ such that $j_l=j_{\sigma_1(l)}^{\prime}$ and $\sigma_2$ be a permutation  of $\{1,2,\dots,v_1\}$ such that $i_m=i_{\sigma_2(m)}^{\prime}$. The numbers of the permutations $\sigma_1$ and $\sigma_2$ are $r_1!$ and $v_1!$ respectively. Then
\begin{eqnarray}\nonumber
&&\mathbb{E}\Bigg[\left(\prod_{l=1}^{r_0}\bar{A}_{ij_l}\bar{A}_{ij_l^{\prime}}\right)\left(\prod_{l=r_0+1}^r\bar{A}_{ij_l}^{\lambda_{r;l}}\bar{A}_{ij_l^{\prime}}^{\lambda_{r;l}}\right)\Bigg]\\ \nonumber
&=&\mathbb{E}\Bigg[\left(\prod_{l=1}^{r_1}\bar{A}_{ij_l}^{\lambda_{r;l}+\lambda_{r;\sigma_1(l)}}\right)\Bigg]\left(\prod_{l=r_1+1}^{r}\mathbb{E}[\bar{A}_{ij_l}^{\lambda_{r;l}}]\mathbb{E}[\bar{A}_{ij_l^{\prime}}^{\lambda_{r;l}}]\right)\\ \label{tueeq9}
&=&O\left(p_n^{2r-r_1}\right).
\end{eqnarray}
Similarly, we have
\begin{eqnarray}\nonumber
&&\mathbb{E}\Bigg[\left(\prod_{m=1}^{v_0}\bar{A}_{ji_m}\bar{A}_{ji_m^{\prime}}\right)\left(\prod_{m=v_0+1}^{v}\bar{A}_{ji_m}^{\gamma_{v;m}}\bar{A}_{ji_m^{\prime}}^{\gamma_{v;m}}\right)\Bigg]\\ \nonumber
&=&\mathbb{E}\Bigg[\left(\prod_{m=1}^{v_1}\bar{A}_{ji_l}^{\gamma_{v;m}+\gamma_{v;\sigma_2(m)}}\right)\Bigg]\left(\prod_{m=v_1+1}^{v}\mathbb{E}[\bar{A}_{ji_m}^{\gamma_{v;m}}]\mathbb{E}[\bar{A}_{ji_m^{\prime}}^{\gamma_{v;m}}]\right)\\ \label{tueeq10}
&=&O\left(p_n^{2v-v_1}\right).
\end{eqnarray}
Note that $2(r+v)-(r_1+v_1)\leq 2(s+t)-1$. By (\ref{tueeq7}), (\ref{tueeq8}), (\ref{tueeq9}), (\ref{tueeq10}) and $(C5)$ of Assumption \ref{assump}, the sum  in (\ref{tueeq7}) over the indices $i=i^{\prime}$, $j=j^{\prime}$, (\ref{tueeq11}) and (\ref{tueeq12}) is bounded by
\begin{eqnarray}\nonumber
(np_n)^{2+2r-r_1+2v-v_1}f^{(s,t)}(np_n,np_n)^2&\leq& (np_n)(np_n)^{2(s+t)}f^{(s,t)}(np_n,np_n)^2\\ \label{tueeq13}
&=&o\left(\sigma_n^2\right).
\end{eqnarray}

Consider case (ii) $i=j^{\prime}$ and $j=i^{\prime}$. If there is an index $l_1$ with $1\leq l_1\leq r_0$ such that $j_{l_1}\notin\{i_1^{\prime},\dots,i_v^{\prime}\}$, then $E=0$ (the same argument as in (\ref{suneq1})). If there is an index $m_1$ with $1\leq m_1\leq v_0$ such that $i_{m_1}^{\prime}\notin\{j_1 ,\dots,j_r\}$, then $E=0$ (the same argument as in (\ref{suneq1})). Then 
$E\neq 0$ implies the following
\[\{j_1,j_2,\dots,j_{r_0}\}\subset\{i_1^{\prime},i_2^{\prime},\dots,i_v^{\prime}\}, \ \ \ \ \{i_1^{\prime},i_2^{\prime},\dots,i_{v_0}^{\prime}\}\subset\{j_1,j_2,\dots,j_r\}, \]
\[\{i_1,i_2,\dots,i_{v_0}\}\subset\{j_1^{\prime},j_2^{\prime},\dots,j_v^{\prime}\}, \ \ \ \ \{j_1^{\prime},j_2^{\prime},\dots,j_{r_0}^{\prime}\}\subset\{i_1,i_2,\dots,i_v\}. \]
Without loss of generality,
suppose 
\begin{equation}\label{suneq4}\{j_1,j_2,\dots,j_{r_1}\}=\{i_1^{\prime},i_2^{\prime},\dots,i_{r_1}^{\prime}\},\ \ \ \ \ \ \ \{j_{r_1+1},\dots,j_r\}\cap\{i_{r_1+1}^{\prime},\dots,i_{v}^{\prime}\}=\emptyset,
\end{equation}
\begin{equation}\label{suneq5}\{i_1,i_2,\dots,i_{v_1}\}=\{j_1^{\prime},j_2^{\prime},\dots,j_{v_1}^{\prime}\},\ \ \ \ \{i_{v_1+1},\dots,i_v\}\cap\{j_{v_1+1}^{\prime},\dots,j_{v}^{\prime}\}=\emptyset,
\end{equation}
for some $r_1$ with $\max\{r_0,v_0\}\leq r_1\leq \min\{r,v\}$ and $v_1$ with $\max\{r_0,v_0\}\leq v_1\leq \min\{r,v\}$. There are at most $n^{2+2r-r_1+2v-v_1}$ possible choices for indices $i,j,i_1,\dots,i_v,j_1,\dots,j_r$, $i^{\prime},j^{\prime},i_1^{\prime},\dots,i_v^{\prime}$, $j_1^{\prime},\dots,j_r^{\prime}$ satisfying $i=j^{\prime}$, $j=i^{\prime}$, (\ref{suneq4}) and (\ref{suneq5}). Let $\sigma_1$ be a permutation  of $\{1,2,\dots,r_1\}$ such that $j_l=i_{\sigma_1(l)}^{\prime}$ and $\sigma_2$ be a permutation  of $\{1,2,\dots,v_1\}$ such that $i_m=j_{\sigma_2(m)}^{\prime}$. Then
\begin{eqnarray*}\nonumber
E&=&\mathbb{E}\Bigg[\left(\prod_{l=1}^{r_1}\bar{A}_{ij_l}^{\lambda_{r;l}+\gamma_{v;\sigma_1(l)}}\right)\left(\prod_{l=r_1+1}^r\bar{A}_{ij_l}^{\lambda_{r;l}}\right)\left(\prod_{m=r_1+1}^r\bar{A}_{ii_m^{\prime}}^{\gamma_{v;m}}\right)    \\ &&\times \left(\prod_{m=1}^{v_1}\bar{A}_{ji_m}^{\lambda_{r;m}+\gamma_{v;\sigma_2(m)}}\right)\left(\prod_{m=v_1+1}^v\bar{A}_{ji_m}^{\gamma_{v;m}}\right)\left(\prod_{m=v_1+1}^r\bar{A}_{jj_m^{\prime}}^{\lambda_{v;m}}\right)\Bigg]    \\
&=&\Theta\left(p_n^{2(r+v)-r_1-v_1}\right).
\end{eqnarray*}
Then the sum in (\ref{tueeq7}) over the indices $i=j^{\prime}$, $j=i^{\prime}$, (\ref{suneq4}) and (\ref{suneq5}) is bounded by
\begin{eqnarray}\nonumber
&&(np_n)^{2+2r-r_1+2v-v_1}f^{(s,t)}(np_n,np_n)^2\\ \label{suneq6}
&\leq& (np_n)(np_n)^{2(s+t)}f^{(s,t)}(np_n,np_n)^2=o\left(\sigma_n^2\right).
\end{eqnarray}

\medskip

Consider the case $i\neq i^{\prime}$ and $j=j^{\prime}$. For any $l_0\in\{1,2,\dots,r\}$, $\{i,j_{l_0}\}\neq \{j,i_l\}$ and $\{i,j_{l_0}\}\neq \{j,i_l^{\prime}\}$ for any $1\leq l\leq v$. If $r_0\geq2$, it is not possible that $\{i,j_1\}=\{i^{\prime},j_{l_1}^{\prime}\}$ and $\{i,j_2\}=\{i^{\prime},j_{l_2}^{\prime}\}$ for distinct $l_1$ and $l_2$. Then $E=0$ (the same argument as in (\ref{suneq1})). Let $r_0=1$. In this case, $i=j_1^{\prime}$ and $i^{\prime}=j_1$. Otherwise $E=0$. Suppose 
$\{i_1,\dots,i_{v_1}\}=\{i_1^{\prime},\dots,i_{v_1}^{\prime}\}$ and $\{i_{v_1+1},\dots,i_v\}\cap\{i_{v_1+1}^{\prime},\dots,i_v^{\prime}\}=\emptyset$ for $v_0\leq v_1\leq v$. Suppose 
$\{j_2,\dots,j_{r_1}\}=\{j_2^{\prime},\dots,j_{r_1}^{\prime}\}$ and $\{j_{r_1+1},\dots,j_r\}\cap\{j_{r_1+1}^{\prime},\dots,j_r^{\prime}\}=\emptyset$ for $1\leq r_1\leq r$. There are at most $n^{2+2(r+v)-v_1-r_1}$ choices for the indices $i,j,i_1,\dots,i_v,j_1,\dots,j_r,i^{\prime},j^{\prime},i_1^{\prime},\dots,i_v^{\prime},j_1^{\prime},\dots,j_r^{\prime}$ satisfying these conditions. 
Then
\begin{eqnarray}\nonumber
E&=&\mathbb{E}\Bigg[\bar{A}_{ij_1}^2\left(\prod_{l=2}^{r}\bar{A}_{ij_l}^{\lambda_{r;l}}\right)\left(\prod_{l=2}^{r_1}\bar{A}_{j_1j_l}^{\lambda_{r;l}}\right)\left(\prod_{l=r_1+1}^{r}\bar{A}_{j_1j_l^{\prime}}^{\lambda_{r;l}}\right)\Bigg]\\ \nonumber
&&\times\mathbb{E}\Bigg[\left(\prod_{m=1}^{v_1}\bar{A}_{ji_m}^{2\gamma_{v;m}}\right)\left(\prod_{m=v_1+1}^{v}\bar{A}_{ji_m}^{\gamma_{v;m}}\bar{A}_{ji_m^{\prime}}^{\gamma_{v;m}}\right)\Bigg]\\ \nonumber
&=&O\left(p_n^{1+2(r-1)+v_1+2(v-v_1)}\right),
\end{eqnarray}
and hence the sum over $i\neq i^{\prime}$ and $j=j^{\prime}$ in (\ref{tueeq7})  is bounded by
\begin{equation}\label{fneq1}
(np_n)^{2(r+v)-v_1+1}f^{(s,t)}(np_n,np_n)^2=o(\sigma_n^2).
\end{equation}

Consider the case $i= i^{\prime}$ and $j\neq j^{\prime}$. If $v_0\geq1$, the summation is similarly bounded by (\ref{fneq1}). Suppose $v_0=0$. Suppose $\{j_1,\dots,j_{r_1}\}=\{j_1^{\prime},\dots,j_{r_1}^{\prime}\}$ and $\{j_{r_1+1},\dots,j_{r}\}\cap\{j_{r_1+1}^{\prime},\dots,j_{r}^{\prime}\}=\emptyset$ for $r_0\leq r_1\leq r$. There are at most $n^{3+2r-r_1+2v}$ choices for the indices $i,j,i_1,\dots,i_v,j_1,\dots$, $j_r$, $i^{\prime},j^{\prime},i_1^{\prime},\dots,i_v^{\prime},j_1^{\prime},\dots,j_r^{\prime}$ satisfying these conditions. Then 
\begin{eqnarray*}
E
&=&\mathbb{E}\Bigg[\left(\prod_{l=1}^{r_1}\bar{A}_{ij_l}^{2\lambda_{r;l}}\right)\left(\prod_{l=r_1+1}^r\bar{A}_{ij_l}^{\lambda_{r;l}}\right)\left(\prod_{l=r_1+1}^r\bar{A}_{ij_l^{\prime}}^{\lambda_{r;l}}\right)\\
&&\times\left(\prod_{m=1}^{v}\bar{A}_{ji_m}^{\gamma_{v;m}}\right)\left(\prod_{m=1}^{v}\bar{A}_{j^{\prime}i_m^{\prime}}^{\gamma_{v;m}}\right)\Bigg]\\
&=&O\left(p_n^{2r-r_1+2v}\right).
\end{eqnarray*}
Then the sum over $i=i^{\prime}$ and $j\neq j^{\prime}$ in (\ref{tueeq7})  is bounded by
\begin{equation}\label{fneq2}
n(np_n)^{2+2r-r_1+2v}f^{(s,t)}(np_n,np_n)^2=o(\sigma_n^2).
\end{equation}

Suppose $i= j^{\prime}$ and $j\neq i^{\prime}$. If $r_0\geq2$ or $v_0\geq2$, then $E=0$ (similar to the argument in the case $i\neq i^{\prime}$ and $j=j^{\prime}$). Let $r_0=1$. If $v_0=1$, then  $j=j_1^{\prime}$ and $i_1=i^{\prime}$. If $v_0=0$, then $j_1=i_{m_1}^{\prime}$ for some $1\leq m_1\leq v$, $j=j_{l_1}^{\prime}$ for some $1\leq l_1\leq r$ and $i_1=i^{\prime}$. Without loss of generality, let $m_1=l_1=1$. Suppose
\[\{j_1,\dots,j_{r_1}\}=\{i_1^{\prime},\dots,i_{r_1}^{\prime}\},\ \ \ \{j_{r_1+1},\dots,j_r\}\cap\{i_{r_1+1}^{\prime},\dots,i_{v}^{\prime}\}=\emptyset.\]
\[\{i_2,\dots,i_{v_1}\}=\{j_2^{\prime},\dots,j_{v_1}^{\prime}\},\ \ \ \{i_{v_1+1},\dots,i_{v}\}\cap\{j_{v_1+1}^{\prime},\dots,j_{r}^{\prime}\}=\emptyset,\]
where $1\leq r_1,v_1\leq\min\{r,v\}$. There are at most $n^{2(r+v)-r_1-v_1+2}$ choices for the indices $i,j,i_1,\dots,i_v,j_1,\dots,j_r,i^{\prime},j^{\prime},i_1^{\prime},\dots,i_v^{\prime},j_1^{\prime},\dots,j_r^{\prime}$ satisfying these conditions.
In this case,
\begin{eqnarray}\nonumber
E&=&\mathbb{E}\Bigg[\bar{A}_{ji_1}^2\left(\prod_{l=1}^{r_1}\bar{A}_{ij_l}^{\lambda_{r;l}+\gamma_{v;l}}\right)\left(\prod_{l=r_1+1}^{r}\bar{A}_{ij_l}^{\lambda_{r;l}}\right)\left(\prod_{l=r_1+1}^{v}\bar{A}_{ii_l^{\prime}}^{\gamma_{v;l}}\right)\Bigg]\\ \nonumber
&&\times\mathbb{E}\Bigg[\left(\prod_{m=2}^{v}\bar{A}_{ji_m}^{\gamma_{v;m}}\right)\left(\prod_{l=2}^{v_1}\bar{A}_{i_1i_l}^{\lambda_{r;l}}\right)\left(\prod_{l=v_1+1}^{r}\bar{A}_{i_1j_l^{\prime}}^{\lambda_{r;l}}\right)\Bigg]\\ \nonumber
&=&O\left(p_n^{1+r+(r-v_1)+v-r_1+v-1+v_1-1}\right).
\end{eqnarray}
Then the sum over $i= j^{\prime}$ and $j\neq i^{\prime}$ in (\ref{tueeq7})  is bounded by
\begin{equation}\label{fneq3}
(np_n)^{2(r+v)-r_1-1}f^{(s,t)}(np_n,np_n)^2=o(\sigma_n^2).
\end{equation}

The case $i\neq j^{\prime}$ and $j= i^{\prime}$ can be similarly bounded as in (\ref{fneq3}).

Now we consider $\{i,j\}\cap\{i^{\prime},j^{\prime}\}=\emptyset$. If $r_0\geq3$, then at least one of $\{i,j_1\}$, $\{i,j_2\}$ and $\{i,j_3\}$ is not in $\{\{j^{\prime},i_{m_1}^{\prime}\},\{i^{\prime},j_l^{\prime}\}\}$ for any $m,m_1,l$. Hence $E=0$. Similarly, if $v_0\geq3$, $E=0$.

Suppose $r_0=v_0=2$. Then $\{i,j_1\}=\{i^{\prime},j_l^{\prime}\}$ for some $1\leq l\leq r$ or $\{i,j_1\}=\{j^{\prime},i_m^{\prime}\}$ for some $1\leq m\leq v$. Otherwise $E=0$. Without loss of generality, suppose
$\{i,j_1\}=\{i^{\prime},j_l^{\prime}\}$. In this case, $i=j_l^{\prime}$ and $j_1=i^{\prime}$. If $l\geq3$, then $\{i^{\prime},j_1^{\prime}\}=\{j,i_{m_1}\}$ and $\{i^{\prime},j_2^{\prime}\}=\{j^{\prime},i_{m_2}^{\prime}\}$ or $\{i^{\prime},j_1^{\prime}\}=\{j^{\prime},i_{m_2}^{\prime}\}$ and $\{i^{\prime},j_2^{\prime}\}=\{j,i_{m_1}\}$ (otherwise $E=0$). Either case is impossible, due to the fact that $j^{\prime}\neq i^{\prime}$, $j^{\prime}\neq j_{l_1}^{\prime}$ for any $l_1$. Hence $l=1$ or $l=2$. Without loss of generality, let $l=1$. In this case, $\{i^{\prime},j_2^{\prime}\}=\{j,i_{m_3}\}$ and $i^{\prime}=i_{m_3}$ and $j_2^{\prime}=j$. If $m_3\geq 3$, then $\{j,i_1\}=\{j^{\prime},i_{m_4}^{\prime}\}$ and $\{j,i_2\}=\{j^{\prime},i_{m_5}^{\prime}\}$ (otherwise $E=0$), which is not possible. Hence $m_3=1$ or $m_3=2$. Let $m_3=1$ (the argument for $m_3=2$ is the same). Then $\{j,i_2\}=\{j^{\prime},i_{m_6}^{\prime}\}$. If $m_6\geq3$, then $\{j^{\prime},i_{1}^{\prime}\}=\{i,j_{l_2}\}$ and $\{j^{\prime},i_{2}^{\prime}\}=\{i,j_{l_3}\}$ (otherwise $E=0$), which is not possible. Hence, $m_6=1$ or $m_6=2$. Let $m_6=1$ (the argument for $m_6=2$ is the same). Then $\{j^{\prime},i_{2}^{\prime}\}=\{i,j_{2}\}$, $i=i_{2}^{\prime}$ and $j_2=j^{\prime}$. There are at most $n^{2(r+v)-4}$ possible choices of the indices $i,j,i_1,\dots,i_v,j_1,\dots,j_r,i^{\prime},j^{\prime},i_1^{\prime},\dots,i_v^{\prime},j_1^{\prime},\dots,j_r^{\prime}$ satisfying these conditions. Then the sum in (\ref{tueeq7})  over these indices is bounded by
\begin{equation}\label{moneq1}
(np_n)^{2(r+v)-4}f^{(s,t)}(np_n,np_n)^2=o(\sigma_n^2).
\end{equation}

Suppose $r_0=2$ and $v_0=1$. Then $\{i,j_1\}=\{i^{\prime},j_l^{\prime}\}$ and $\{i,j_2\}=\{j^{\prime},i_m^{\prime}\}$ or $\{i,j_1\}=\{j^{\prime},i_m^{\prime}\}$ and $\{i,j_2\}=\{i^{\prime},j_l^{\prime}\}$ for some $l,m$. Otherwise $E=0$. Without loss of generality, let $\{i,j_1\}=\{i^{\prime},j_l^{\prime}\}$ and $\{i,j_2\}=\{j^{\prime},i_m^{\prime}\}$. By a similar argument as in the previous paragraph, $l=1$ or $l=2$. Let $l=1$. Then $\{i^{\prime},j_2^{\prime} \}=\{j,i_{m_3}\}$. If $m_3=1$ and $m=1$,  the sum in (\ref{tueeq7})  over these indices is bounded by
\begin{equation}\label{moneq2}
p_n(np_n)^{2(r+v)-2}f^{(s,t)}(np_n,np_n)^2=o(\sigma_n^2),
\end{equation}
If $m_3\geq2$ or $m\geq2$,  the sum in (\ref{tueeq7})  over these indices is bounded by
\begin{equation}\label{moneq3}
p_n^2(np_n)^{2(r+v)-4}f^{(s,t)}(np_n,np_n)^2=o(\sigma_n^2).
\end{equation}

Suppose $r_0=2$ and $v_0=0$. Then $\{i,j_1\}=\{i^{\prime},j_l^{\prime}\}$ and $\{i,j_2\}=\{j^{\prime},i_m^{\prime}\}$ or $\{i,j_1\}=\{j^{\prime},i_m^{\prime}\}$ and $\{i,j_2\}=\{i^{\prime},j_l^{\prime}\}$ for some $l,m$. Otherwise $E=0$. Without loss of generality, let $\{i,j_1\}=\{i^{\prime},j_l^{\prime}\}$ and $\{i,j_2\}=\{j^{\prime},i_m^{\prime}\}$. Then $l=1$ or $l=2$. Without loss of generality, let $l=1$. Then $\{i^{\prime},j_2^{\prime}\}=\{j,i_{m_1}\}$. The sum in (\ref{tueeq7})  over these indices is bounded by
\begin{equation}\label{moneq4}
p_n(np_n)^{2(r+v)-2}f^{(s,t)}(np_n,np_n)^2=o(\sigma_n^2).
\end{equation}

Suppose $r_0=1$ and $v_0=1$.  Then $\{i,j_1\}=\{i^{\prime},j_l^{\prime}\}$ or $\{i,j_1\}=\{j^{\prime},i_m^{\prime}\}$ for some $l,m$. Without loss of generality, let $\{i,j_1\}=\{i^{\prime},j_l^{\prime}\}$. Suppose $l=1$. Then $\{j,i_1\}=\{i^{\prime},j_{l_1}^{\prime}\}$ or $\{j,i_1\}=\{j^{\prime},i_{m_1}^{\prime}\}$. If $\{j,i_1\}=\{j^{\prime},i_{m_1}^{\prime}\}$ with $m_1=1$, the sum in (\ref{tueeq7})  over these indices is bounded by
\begin{equation}\label{moneq5}
(np_n)^{2(r+v)}f^{(s,t)}(np_n,np_n)^2=o(\sigma_n^2).
\end{equation}
If $m_1\geq2$, then $\{j^{\prime},i_{1}^{\prime}\}=\{i,j_{m_2}\}$ for $m_2\geq2$. Then the sum in (\ref{tueeq7})  over these indices is bounded by
\begin{equation}\label{moneq6}
p_n(np_n)^{2(r+v)-2}f^{(s,t)}(np_n,np_n)^2=o(\sigma_n^2).
\end{equation}
The case $\{j,i_1\}=\{i^{\prime},j_{l_1}^{\prime}\}$ can be similarly bounded as in (\ref{moneq6}).

Suppose $l\geq2$. Then $\{i^{\prime},j_1^{\prime}\}=\{j,i_{m_1}\}$. If $m_1=1$, then $\{j^{\prime},i_{1}^{\prime}\}=\{i,j_{l_1}\}$ for $l_1\geq2$ or $\{j^{\prime},i_{1}^{\prime}\}=\{j,i_{m_2}\}$ for $m_2\geq2$. In this case, the sum in (\ref{tueeq7})  over these indices is bounded by
\begin{equation}\label{moneq7}
p_n(np_n)^{2(r+v)-2}f^{(s,t)}(np_n,np_n)^2=o(\sigma_n^2).
\end{equation}
If $m_1\geq2$, then $\{j^{\prime},i_{1}^{\prime}\}=\{j,i_{1}\}$. The sum in (\ref{tueeq7})  over these indices is bounded by
\begin{equation}\label{moneq8}
p_n(np_n)^{2(r+v)-2}f^{(s,t)}(np_n,np_n)^2=o(\sigma_n^2).
\end{equation}

Suppose $r_0=1$ and $v_0=0$. Then $\{i,j_1\}=\{i^{\prime},j_{l_1}^{\prime}\}$ or $\{i,j_1\}=\{j^{\prime},i_{m_1}^{\prime}\}$. Suppose $\{i,j_1\}=\{j^{\prime},i_{m_1}^{\prime}\}$. Then $\{i^{\prime},j_1^{\prime}\}=\{j,i_{m_2}\}$. In this case,  the sum in (\ref{tueeq7})  over these indices is bounded by
\begin{equation}\label{moneq9}
(np_n)^{2(r+v)}f^{(s,t)}(np_n,np_n)^2=o(\sigma_n^2).
\end{equation}
Suppose $\{i,j_1\}=\{i^{\prime},j_{l_1}^{\prime}\}$. If $l_1=1$, then 
 the sum in (\ref{tueeq7})  over these indices is bounded by
\begin{equation}\label{moneq10}
n(np_n)^{2(r+v)+1}f^{(s,t)}(np_n,np_n)^2=o(\sigma_n^2).
\end{equation}
If $l_1\geq2$, then $\{i^{\prime},j_{1}^{\prime}\}=\{j,i_{m_3}\}$. In this case,
the sum in (\ref{tueeq7})  over these indices is bounded by
\begin{equation}\label{moneq11}(np_n)^{2(r+v)}f^{(s,t)}(np_n,np_n)^2=o(\sigma_n^2).
\end{equation}

\medskip

Now we consider the case (II): $\lambda_{r;l}\geq2$ for all $l=1,2,\dots, r$  and $\gamma_{v;m}\geq2$ for all $m=1,2,\dots, v$.  In this case, $r\leq\frac{s}{2}$ and $v\leq\frac{t}{2}$. The expectation of $V_{rv}$ is equal to
\begin{eqnarray}\nonumber
&&\mathbb{E}[V_{rv}]\\ \nonumber
&=&\sum_{i\neq j}\sum_{\substack{j_1,j_2,\dots,j_r\notin\{i,j\}\\ j_1\neq j_2\neq \dots\neq j_r\\
i_1,i_2,\dots,i_v\notin\{i,j\}\\i_1\neq i_2\neq \dots\neq i_v}}p_nw_{ij}f^{(s,t)}(w_{i(j)},w_{j(i)})\prod_{l=1}^r\mathbb{E}[(A_{ij_l}-p_nw_{ij_l})^{\lambda_{r;l}}]\\ \label{tueeq16}
&&\times\prod_{m=1}^v\mathbb{E}[(A_{ji_m}-p_nw_{ji_m})^{\gamma_{v;m}}].  
\end{eqnarray}
Next we bound the variance of $V_{rv}$, that is, $Var(V_{rv})=\mathbb{E}\left[\left(V_{rv}-\mathbb{E}[V_{rv}]\right)^2\right]$. 
Let $\eta_l\in\{0,1\}$ for $l=1,2,\dots,r$ and  $\xi_l\in\{0,1\}$ for $l=1,2,\dots,v$. Then we have
\begin{eqnarray}\nonumber
    &&\prod_{l=1}^r(A_{ij_l}-p_nw_{ij_l})^{\lambda_{r;l}}\prod_{m=1}^v(A_{ji_m}-p_nw_{ji_m})^{\gamma_{v;m}}\\ \nonumber
    &=&    \prod_{l=1}^r\Big[(A_{ij_l}-p_nw_{ij_l})^{\lambda_{r;l}}-\mathbb{E}\left[(A_{ij_l}-p_nw_{ij_l})^{\lambda_{r;l}}\right]\\ \nonumber
    &&+\mathbb{E}\left[(A_{ij_l}-p_nw_{ij_l})^{\lambda_{r;l}}\right]\Big]\\  \nonumber
    &&\times \prod_{l=1}^v\Big[(A_{ji_l}-p_nw_{ji_l})^{\gamma_{v;l}}-\mathbb{E}\left[(A_{ji_l}-p_nw_{ji_l})^{\gamma_{v;l}}\right]\\  \nonumber
    &&+\mathbb{E}\left[(A_{ji_l}-p_nw_{ji_l})^{\gamma_{v;l}}\right]\Big]\\ \nonumber    &=&\sum_{\substack{\eta_1,\dots,\eta_r\in\{0,1\}\\ \xi_1,\dots,\xi_v\in\{0,1\}}}\prod_{l=1}^r\Big[(A_{ij_l}-p_nw_{ij_l})^{\lambda_{r;l}}-\mathbb{E}\left[(A_{ij_l}-p_nw_{ij_l})^{\lambda_{r;l}}\right]\Big]^{\eta_l}\\\nonumber 
    &&\times\Big[\mathbb{E}\left[(A_{ij_l}-p_nw_{ij_l})^{\lambda_{r;l}}\right]\Big]^{1-\eta_l}\\ \nonumber
    &&\times \prod_{l=1}^r\Big[(A_{ji_l}-p_nw_{ji_l})^{\gamma_{v;l}}-\mathbb{E}\left[(A_{ji_l}-p_nw_{ji_l})^{\gamma_{v;l}}\right]\Big]^{\xi_l}\\ \label{tueeq14} 
&&\times\Big[\mathbb{E}\left[(A_{ji_l}-p_nw_{ji_l})^{\gamma_{v;l}}\right]\Big]^{1-\xi_l}. 
\end{eqnarray}
For convenience, denote 
\[X_{ij}(x)=(A_{ij}-p_nw_{ij})^{x}-\mathbb{E}[(A_{ij}-p_nw_{ij})^{x}]\]
and $Y_{ij}(x)=\mathbb{E}[(A_{ij}-p_nw_{ij})^{x}]$.

Since $i\neq j$, $j_l\notin\{i,j\}$ for $l=1,2,\dots,r$ and $i_l\notin\{i,j\}$ for $l=1,2,\dots,v$, then $\{i,j_l\}\neq \{j,i_m\}$ for any $l=1,2,\dots,r$ and $m=1,2,\dots,v$. Consequently, $A_{ij_l}$ ($l=1,2,\dots,r$) and $A_{ji_m}$ ($m=1,2,\dots,v$) are independent.

By (\ref{tueeq4}), (\ref{tueeq16}) and (\ref{tueeq14}), $V_{rv}-\mathbb{E}[V_{rv}]$ does not contain the term corresponding to
$\eta_l=0$ for all $l\in\{1,2,\dots,r\}$ and $\xi_m=0$ for 
 all $m\in\{1,2,\dots,v\}$. Hence we assume $\eta_1+\dots+\eta_r\geq 1$ or $\xi_1+\dots+\xi_v\geq1$. Without loss of generality, let $\eta_1=\eta_2=\dots=\eta_{l_0}=1$, $\eta_{l_0+1}=\dots=\eta_r=0$, $\xi_1=\xi_2=\dots=\xi_{m_0}=1$ and $\xi_{m_0+1}=\dots=\xi_v=0$, where $1\leq l_0+m_0\leq r+v$. 
 In this case,
\begin{eqnarray}\nonumber
\prod_{l=1}^r\Big[X_{ij_l}(\lambda_{r;l})\Big]^{\eta_l}\Big[Y_{ij_l}(\lambda_{r;l})\Big]^{1-\eta_l}
    =\left( \prod_{l=1}^{l_0} X_{ij_l}(\lambda_{r;l}) \right)\left(\prod_{l=l_0+1}^rY_{ij_l}(\lambda_{r;l})\right),
\end{eqnarray}
\begin{eqnarray}\nonumber
    \prod_{l=1}^v\Big[X_{ji_l}(\gamma_{v;l})\Big]^{\xi_l}\Big[Y_{ji_l}(\gamma_{v;l})\Big]^{1-\xi_l}=\left(\prod_{l=1}^{m_0}X_{ji_l}(\gamma_{v;l})\right)\left(\prod_{l=m_0+1}^vY_{ji_l}(\gamma_{v;l})\right).
\end{eqnarray}
Denote
\begin{eqnarray}\nonumber
&&V_{rv}(l_0,m_0) \\ \nonumber
   &=& \sum_{\substack{i\neq j\\ j_1,j_2,\dots,j_r\notin\{i,j\}\\ j_1\neq j_2\neq \dots\neq j_r\\
i_1,i_2,\dots,i_v\notin\{i,j\}\\i_1\neq i_2\neq \dots\neq i_v}}p_nw_{ij}f^{(s,t)}(w_{i(j)},w_{j(i)})\left( \prod_{l=1}^{l_0} X_{ij_l}(\lambda_{r;l}) \right) \\ \label{wedeq1}
&&\times\left(\prod_{l=l_0+1}^rY_{ij_l}(\lambda_{r;l})\right)\left(\prod_{l=1}^{m_0}X_{ji_l}(\gamma_{v;l})\right)\left(\prod_{l=m_0+1}^vY_{ji_l}(\gamma_{v;l})\right).
\end{eqnarray}
We only need to consider the variance of $V_{rv}(l_0,m_0)$. Since $\lambda_{r;l}\geq2$ and $w_{ij}\in[\beta,1]$, then
\begin{eqnarray*}
&&\mathbb{E}[(A_{ij_l}-p_nw_{ij_l})^{\lambda_{r;l}}]\\
&=&(1-p_nw_{ij_l})^{\lambda_{r;l}}p_nw_{ij_l}+(-p_nw_{ij_l})^{\lambda_{r;l}}(1-p_nw_{ij_l})\\
&=&\Theta(p_n).
\end{eqnarray*}
Hence one has
\[\prod_{l=l_0+1}^rY_{ij_l}(\lambda_{r;l})=\Theta(p_n^{r-l_0}),\ \ \ \ \ 
 \prod_{l=m_0+1}^vY_{ji_l}(\gamma_{v;l})=\Theta(p_n^{v-m_0}).\]
For convenience, denote 
\[a_{ij_{1},\dots,j_{l_0}}=\sum_{j_{l_0+1},\dots,j_r\notin\{j_1,\dots,j_{l_0}\}}\prod_{l=l_0+1}^rY_{ij_l}(\lambda_{r;l}),\]
\[b_{ji_{1},\dots,i_{m_0}}=\sum_{i_{m_0+1},\dots,i_v\notin\{i_1,\dots,i_{m_0}\}}\prod_{l=m_0+1}^vY_{ji_l}(\gamma_{v;l}).\]
Then
\begin{equation}\label{wedeq4}
a_{ij_{1},\dots,j_{l_0}}=\Theta((np_n)^{r-l_0}),\ \ \ \ \ \ b_{ji_{1},\dots,i_{m_0}}=\Theta((np_n)^{v-m_0}).
\end{equation}
In this case,  $V_{rv}(l_0,m_0)$ is written as
\begin{eqnarray}\nonumber
V_{rv}(l_0,m_0)
&=&\sum_{\substack{i\neq j\\ j_1,j_2,\dots,j_r\notin\{i,j\}\\ j_1\neq j_2\neq \dots\neq j_r\\
i_1,i_2,\dots,i_v\notin\{i,j\}\\i_1\neq i_2\neq \dots\neq i_v}}p_nw_{ij}f^{(s,t)}(w_{i(j)},w_{j(i)})a_{ij_{1},\dots,j_{l_0}}b_{ji_{1},\dots,i_{m_0}}\\ \label{wedeq2}
&&\times\left( \prod_{l=1}^{l_0} X_{ij_l}(\lambda_{r;l}) \right)\left(\prod_{l=1}^{m_0}X_{ji_l}(\gamma_{v;l})\right).
\end{eqnarray}
Next we bound the variance of $V_{rv}(l_0,m_0)$. Since $i\neq j$, $j_l\notin\{i,j\}$  and $i_m\notin\{i,j\}$ for all $l,m$, $X_{ij_l}(\lambda_{r;l}) $ and $X_{ji_l}(\gamma_{v;l})$ are independent. By definition, $\mathbb{E}[X_{ij_l}(\lambda_{r;l})]=\mathbb{E}[X_{ji_l}(\gamma_{v;l})]=0$. Then  $\mathbb{E}[V_{rv}(l_0,m_0)]=0$ and $Var\left(V_{rv}(l_0,m_0)\right)=\mathbb{E}[V_{rv}(l_0,m_0)^2]$, that is,
\begin{eqnarray}\nonumber
&&Var\left(V_{rv}(l_0,m_0)\right)\\ \nonumber
&=&\mathbb{E}\Bigg[\sum_{i,j,j_1,\dots,j_{l_0},i_1,\dots,i_{m_0}}p_nw_{ij}f^{(s,t)}(w_{i(j)},w_{j(i)})a_{ij_{1},\dots,j_{l_0}}b_{ji_{1},\dots,i_{m_0}}\\ \nonumber
&&\times\left( \prod_{l=1}^{l_0} X_{ij_l}(\lambda_{r;l}) \right)\left(\prod_{l=1}^{m_0}X_{ji_l}(\gamma_{v;l})\right)\Bigg]^2\\ \nonumber
&=&\sum_{\substack{i,j,j_1,\dots,j_{l_0},i_1,\dots,i_{m_0}\\ i^{\prime},j^{\prime},j_1^{\prime},\dots,j_{l_0}^{\prime},i_1^{\prime},\dots,i_{m_0}^{\prime}}}p_nw_{ij}f^{(s,t)}(w_{i(j)},w_{j(i)})p_nw_{i^{\prime}j^{\prime}}f^{(s,t)}(w_{i^{\prime}(j^{\prime})},w_{j^{\prime}(i^{\prime})})\\ \nonumber
&&\times a_{ij_{1},\dots,j_{l_0}}b_{ji_{1},\dots,i_{m_0}} a_{i^{\prime}j_{1}^{\prime},\dots,j_{l_0}^{\prime}}b_{j^{\prime}i_{1}^{\prime},\dots,i_{m_0}^{\prime}}\\ \label{wedeq3}
&&
\times\mathbb{E}\left( \prod_{l=1}^{l_0} X_{ij_l}(\lambda_{r;l}) X_{i^{\prime}j_l^{\prime}}(\lambda_{r;l}) \prod_{l=1}^{m_0}X_{ji_l}(\gamma_{v;l})X_{j^{\prime}i_l^{\prime}}(\gamma_{v;l})\right).
\end{eqnarray}

For each $j_{l_1}$ with $1\leq l_1\leq l_0$, if $\{i,j_{l_1}\}\neq \{i^{\prime},j_l^{\prime}\}$ or $\{j^{\prime},i_m^{\prime}\}$ for some $l\in\{1,2,\dots,l_0\}$ or $m\in\{1,2,\dots,m_0\}$, then
$X_{ij_{l_1}}(\lambda_{r;l_1})$ is independent of $X_{ji_l}(\lambda_{r;l})$, $X_{i^{\prime}j_l^{\prime}}(\lambda_{r;l})$ and $X_{j^{\prime}i_m^{\prime}}(\lambda_{r;m})$. In this case,
\begin{eqnarray}\nonumber
&&\mathbb{E}\left[ \prod_{l=1}^{l_0} X_{ij_l}(\lambda_{r;l}) X_{i^{\prime}j_l^{\prime}}(\lambda_{r;l}) \prod_{l=1}^{m_0}X_{ji_l}(\gamma_{v;l})X_{j^{\prime}i_l^{\prime}}(\gamma_{v;l})\right]\\ \nonumber
&=&\mathbb{E}\left[X_{ij_{l_1}}\right]\mathbb{E}\left[ \prod_{l=1,l\neq l_1}^{l_0} X_{ij_l}(\lambda_{r;l}) X_{i^{\prime}j_l^{\prime}}(\lambda_{r;l}) \prod_{l=1}^{m_0}X_{ji_l}(\gamma_{v;l})X_{j^{\prime}i_l^{\prime}}(\gamma_{v;l})\right]\\ \label{vkeq1}
&=&0.
\end{eqnarray}
Hence, 
for each $j_{l_1}$ with $1\leq l_1\leq l_0$, there exists $l$ or $m$ such that $\{i,j_{l_1}\}=\{i^{\prime},j_l^{\prime}\}$ or $\{i,j_{l_1}\}=\{j^{\prime},i_m^{\prime}\}$. Otherwise, (\ref{vkeq1}) holds. Moreover, if $\{i,j_{l_1}\}=\{i^{\prime},j_l^{\prime}\}$ and $\{i,j_{l_1}\}=\{j^{\prime},i_m^{\prime}\}$, then $\{i^{\prime},j_l^{\prime}\}=\{j^{\prime},i_m^{\prime}\}$. In this case, $i^{\prime}=i_m^{\prime}$ and $j^{\prime}=j_l^{\prime}$ (since $i^{\prime}\neq j^{\prime}$). This is not possible, due to the fact that $i_m^{\prime}\notin\{i^{\prime},j^{\prime}\}$ and $j_l^{\prime}\notin\{i^{\prime},j^{\prime}\}$ for all $m,l$.  Therefore, $\{i,j_{l_1}\}$ can only be equal to  one of  $\{i^{\prime},j_l^{\prime}\}$ and $\{j^{\prime},i_m^{\prime}\}$, but not both. Similarly, $\{j,i_m\}$ can only be equal to one of $\{i^{\prime},j_l^{\prime}\}$ and $\{j^{\prime},i_m^{\prime}\}$.

Let $m_0=0$. If $l_0\geq2$, then $i=i^{\prime}$, $\{j_1,j_2,\dots,j_{l_0}\}=\{j_1^{\prime},j_2^{\prime},\dots,j_{l_0}^{\prime}\}$. Then
\begin{eqnarray}\nonumber
&&\mathbb{E}\left[ \prod_{l=1}^{l_0} X_{ij_l}(\lambda_{r;l}) X_{i^{\prime}j_l^{\prime}}(\lambda_{r;l}) \prod_{l=1}^{m_0}X_{ji_l}(\gamma_{v;l})X_{j^{\prime}i_l^{\prime}}(\gamma_{v;l})\right]
\\ \label{wedeq21}
&=&\mathbb{E}\left[\prod_{l=1}^{l_0} X_{ij_l}(\lambda_{r;l})^2 \right] 
=O(p_n^{l_0}).
\end{eqnarray}
Then the sum over these indices in (\ref{wedeq3}) is bounded by
\begin{equation}\label{ssseq1}
O\left(n(np_n)^{2(r+v)-l_0+2}f^{(s,t)}(np_n,np_n)^2\right)=o(\sigma_n^2).
\end{equation}
If $l_0=1$, then $\{i,j_1\}=\{i^{\prime},j_1^{\prime}\}$. In this case, (\ref{ssseq1}) still holds.

Let $m_0=1$. If $l_0\geq3$, then $i=i^{\prime}$, $\{j_1,j_2,\dots,j_{l_0}\}=\{j_1^{\prime},j_2^{\prime},\dots,j_{l_0}^{\prime}\}$ and $\{j,i_1\}=\{j^{\prime},i_1^{\prime}\}$. Then the sum over these indices in (\ref{wedeq3}) is bounded by
\begin{equation}\label{ssseq2}
O\left((np_n)^{2(r+v)-l_0+1}f^{(s,t)}(np_n,np_n)^2\right)=o(\sigma_n^2).
\end{equation}
If $l_0=2$, there are two situations: (i) $i=i^{\prime}$, $\{j_1,j_2\}=\{j_1^{\prime},j_2^{\prime}\}$ and $\{j,i_1\}=\{j^{\prime},i_1^{\prime}\}$; (ii) $\{i,j_1\}=\{j_{l_1}^{\prime},i^{\prime}\}$, $\{i^{\prime},j_{l_2}^{\prime}\}=\{j,i_1\}$, $\{i,j_2\}=\{j^{\prime},i_1^{\prime}\}$ with $l_1\neq l_2$. For the case (i), the sum over these indices in (\ref{wedeq3}) is bounded by (\ref{ssseq2}) with $l_0=2$. For case (ii), the sum over these indices in (\ref{wedeq3}) is bounded by 
\begin{equation}\label{ssseq3}
O\left((np_n)^{2(r+v)-2}f^{(s,t)}(np_n,np_n)^2\right)=o(\sigma_n^2).
\end{equation}
If $l_0=1$, then $\{i,j_1\}=\{i^{\prime},j_1^{\prime}\}$ and $\{j,i_1\}=\{j^{\prime},i_1^{\prime}\}$ or $\{i,j_1\}=\{j^{\prime},i_1^{\prime}\}$ and $\{j,i_1\}=\{i^{\prime},j_1^{\prime}\}$. Either case, the sum over these indices in (\ref{wedeq3}) is bounded by 
\begin{equation}\label{ssseq4}
O\left((np_n)^{2(r+v)}f^{(s,t)}(np_n,np_n)^2\right)=o(\sigma_n^2).
\end{equation}
If $l_0=0$, then (\ref{ssseq1}) holds with $l_0$ replaced by $m_0=1$.

Let $m_0=2$. If $l_0\geq3$, then $i=i^{\prime}$, $\{j_1,j_2,\dots,j_{l_0}\}=\{j_1^{\prime},j_2^{\prime},\dots,j_{l_0}^{\prime}\}$, $j=j^{\prime}$ and $\{i_1,i_2\}=\{i_1^{\prime},i_2^{\prime}\}$. In this case, the sum over these indices in (\ref{wedeq3}) is bounded by 
\begin{equation}\label{ssseq5}
O\left((np_n)^{2(r+v)-l_0}f^{(s,t)}(np_n,np_n)^2\right)=o(\sigma_n^2).
\end{equation}
Let $l_0=2$. There are two cases (i) $i=i^{\prime}$, $\{j_1,j_2\}=\{j_1^{\prime},j_2^{\prime}\}$, $j=j^{\prime}$ and $\{i_1,i_2\}=\{i_1^{\prime},i_2^{\prime}\}$; (ii) $\{i,j_1\}=\{j_{l_1}^{\prime},i^{\prime}\}$, $\{j_{l_2}^{\prime},i^{\prime}\}=\{j,i_{m_1}\}$, $\{j,i_{m_2}\}=\{j^{\prime},i_{m_3}^{\prime}\}$, $\{j^{\prime},i_{m_4}^{\prime}\}=\{i,j_2\}$. For case (i), (\ref{ssseq5}) holds. For case (ii),  the sum over these indices in (\ref{wedeq3}) is bounded by 
\begin{equation}\label{ssseq6}
O\left((np_n)^{2(r+v)-2}f^{(s,t)}(np_n,np_n)^2\right)=o(\sigma_n^2).
\end{equation}
The case $l_0=0,1$ are similar to the case $m_0=0,1$ and $l_0=2$ discussed earlier.

The case $m_0\geq3$ is similar to the case $l_0\geq3$.

Now we study the second moment of (\ref{tueeq2}). Note that
\begin{eqnarray}\nonumber
   &&\sum_{i\neq j}f^{(s,t)}(w_{i(j)},w_{j(i)})(d_{i(j)}-w_{i(j)})^s(d_{j(i)}-w_{j(i)})^t(A_{ij}-p_nw_{ij})\\ \nonumber
   &=&\sum_{r=1}^s\sum_{v=1}^tU_{rv},
\end{eqnarray}
where
\begin{eqnarray}\nonumber
U_{rv}
&=&\sum_{i\neq j}\sum_{\substack{j_1,j_2,\dots,j_r\notin\{i,j\}\\ j_1\neq j_2\neq \dots\neq j_r\\
i_1,i_2,\dots,i_v\notin\{i,j\}\\i_1\neq i_2\neq \dots\neq i_v}}f^{(s,t)}(w_{i(j)},w_{j(i)})(A_{ij}-p_nw_{ij})\\ \nonumber
&&\times\prod_{l=1}^r(A_{ij_l}-p_nw_{ij_l})^{\lambda_{r;l}}\prod_{m=1}^v(A_{ji_m}-p_nw_{ji_m})^{\gamma_{v;m}}.
\end{eqnarray}
Since $s,t$ are fixed finite integers less than $k_0$, we only need to bound the variance of $U_{rv}$. Obviously, $\mathbb{E}[U_{rv}]=0$. Then $Var(U_{rv})$  is equal to

\begin{eqnarray}\nonumber
 &&   \mathbb{E}[U_{rv}^2]\\ \nonumber
 &=&\sum_{\substack{i\neq j\\ i^{\prime}\neq j^{\prime}}}\sum_{\substack{j_1,j_2,\dots,j_r\notin\{i,j\}\\ j_1\neq j_2\neq \dots\neq j_r\\
i_1,i_2,\dots,i_v\notin\{i,j\}\\i_1\neq i_2\neq \dots\neq i_v}}\sum_{\substack{i_1^{\prime},i_2^{\prime},\dots,i_r^{\prime}\notin\{i^{\prime},j^{\prime}\}\\ i_1^{\prime}\neq i_2^{\prime}\neq \dots\neq i_r^{\prime}\\
j_1^{\prime},j_2^{\prime},\dots,j_v^{\prime}\notin\{i^{\prime},j^{\prime}\}\\j_1^{\prime}\neq j_2^{\prime}\neq \dots\neq j_v^{\prime}}}f^{(s,t)}(w_{i(j)},w_{i(j)})\\ \nonumber
&&\times f^{(s,t)}(w_{i^{\prime}(j^{\prime})},w_{j^{\prime}(i^{\prime})})\mathbb{E}\Bigg[(A_{ij}-p_nw_{ij})(A_{i^{\prime}j^{\prime}}-p_nw_{i^{\prime}j^{\prime}})\left(\prod_{l=1}^r\bar{A}_{ij_l}^{\lambda_{r;l}}\right)\\ \label{wedeq24}
&&\times\left(\prod_{l=1}^r\bar{A}_{i^{\prime}j_l^{\prime}}^{\lambda_{r;l}}\right)\left(\prod_{m=1}^{v}\bar{A}_{ji_m}^{\gamma_{v;m}}\right)\left(\prod_{m=1}^{v}\bar{A}_{j^{\prime}i_m^{\prime}}^{\gamma_{v;m}}\right)\Bigg].
\end{eqnarray}
Denote the expectation in (\ref{wedeq24}) as $E_1$. Fix $r\in\{1,2,\dots,s\}$ and $v\in\{1,2,\dots,t\}$.
There are two cases: (a) there are some indices $\lambda_{r;l}$ or $\gamma_{v;m}$ are equal to one; (b)  $\lambda_{r;l}\geq2$ and $\gamma_{v;m}\geq2$ for all $l\in\{1,2,\dots,r\}$ and $m\in\{1,2,\dots,v\}$.

Consider case (a) first. Suppose there are some indices $\lambda_{r;l}$ or $\gamma_{v;m}$ which are equal to one. Without loss of generality, let $\lambda_{r;1}=\lambda_{r;2}=\dots=\lambda_{r;r_0}=1$ and $\lambda_{r;l}\geq2$ for $l\in\{r_0+1,\dots,r\}$. Let $\gamma_{v;1}=\gamma_{v;2}=\dots=\gamma_{v;v_0}=1$ and $\gamma_{v;l}\geq2$ for $l\in\{v_0+1,\dots,v\}$. Here $r_0+v_0\geq1$. Then
\begin{eqnarray}\nonumber
 E_1 
 &=&\mathbb{E}\Bigg[(A_{ij}-p_nw_{ij})(A_{i^{\prime}j^{\prime}}-p_nw_{i^{\prime}j^{\prime}})\left(\prod_{l=1}^{r_0}\bar{A}_{ij_l}\right)\left(\prod_{l=r_0+1}^r\bar{A}_{ij_l}^{\lambda_{r;l}}\right)\\ \nonumber
&&\times\left(\prod_{l=1}^{r_0}\bar{A}_{i^{\prime}j_l^{\prime}}\right)\left(\prod_{l=r_0+1}^r\bar{A}_{i^{\prime}j_l^{\prime}}^{\lambda_{r;l}}\right)\left(\prod_{m=1}^{v_0}\bar{A}_{ji_m}\right)\left(\prod_{m=v_0+1}^{v}\bar{A}_{ji_m}^{\gamma_{v;m}}\right)\\ \nonumber
&&\times\left(\prod_{m=1}^{v_0}\bar{A}_{j^{\prime}i_m^{\prime}}\right)\left(\prod_{m=v_0+1}^{v}\bar{A}_{j^{\prime}i_m^{\prime}}^{\gamma_{v;m}}\right)\Bigg].
\end{eqnarray}

 We split the sum of (\ref{wedeq24}) into two cases:   $\{i,j\}\neq\{i^{\prime},j^{\prime}\}$ and  $\{i,j\}=\{i^{\prime},j^{\prime}\}$.

Suppose $\{i,j\}\neq\{i^{\prime},j^{\prime}\}$. Let $v_0\geq2$ and $r_0\geq2$. For $1\leq m_1\leq v_0$, if $\{j,i_{m_1}\}\neq\{i^{\prime},j_l^{\prime}\}$ for any $l$, $\{j,i_{m_1}\}\neq\{j^{\prime},i_m^{\prime}\}$ for any $m$ and $\{j,i_{m_1}\}\neq\{i^{\prime},j^{\prime}\}$, then $A_{ji_{m_1}}$ is independent of $A_{i^{\prime}j_l^{\prime}}$ for all $l$, $A_{j^{\prime}i_m^{\prime}}$ for all $m$ and $A_{i^{\prime}j^{\prime}}$. Then
\begin{eqnarray}\nonumber
 E_1 
 &=&\mathbb{E}\Bigg[(A_{ij}-p_nw_{ij})(A_{i^{\prime}j^{\prime}}-p_nw_{i^{\prime}j^{\prime}})\left(\prod_{l=1}^{r_0}\bar{A}_{ij_l}\right)\left(\prod_{l=r_0+1}^r\bar{A}_{ij_l}^{\lambda_{r;l}}\right)\\ \nonumber
&&\times\left(\prod_{l=1}^{r_0}\bar{A}_{i^{\prime}j_l^{\prime}}\right)\left(\prod_{l=r_0+1}^r\bar{A}_{i^{\prime}j_l^{\prime}}^{\lambda_{r;l}}\right)\left(\prod_{m=1,m\neq m_1}^{v_0}\bar{A}_{ji_m}\right)\left(\prod_{m=v_0+1}^{v}\bar{A}_{ji_m}^{\gamma_{v;m}}\right)\\ \label{TTteq2}
&&\times\left(\prod_{m=1}^{v_0}\bar{A}_{j^{\prime}i_m^{\prime}}\right)
\left(\prod_{m=v_0+1}^{v}\bar{A}_{j^{\prime}i_m^{\prime}}^{\gamma_{v;m}}\right)\Bigg]\mathbb{E}\left[\bar{A}_{ji_{m_1}}\right]=0.
\end{eqnarray}
Hence, $\{j,i_{m_1}\}=\{i^{\prime},j_l^{\prime}\}$ for some $l$ or $\{j,i_{m_1}\}=\{j^{\prime},i_m^{\prime}\}$ for some $m$ or $\{j,i_{m_1}\}=\{i^{\prime},j^{\prime}\}$. 
Similar results hold for $\{i,j\}$, $\{i,j_l\}$ with $1\leq l\leq r_0$, $\{j^{\prime},i_m^{\prime}\}$ with $1\leq m\leq v_0$, $\{i^{\prime},j^{\prime}\}$, and $\{i^{\prime},j_l^{\prime}\}$ with $1\leq l\leq r_0$. (i)
 Suppose $\{j,i_1\}=\{i^{\prime},j_{l_1}^{\prime}\}$. Then $\{j,i_2\}=\{i^{\prime},j_{l_2}^{\prime}\}$ or $\{j,i_2\}=\{i^{\prime},j^{\prime}\}$. If $\{j,i_2\}=\{i^{\prime},j_{l_2}^{\prime}\}$, then $\{i,j\}=\{i^{\prime},j_{l_3}^{\prime}\}$. Since $j_l^{\prime}\neq j^{\prime}$ for all $l$, then either $\{j^{\prime},i_1^{\prime}\}\neq \{j_{l_3}^{\prime},j_l\}$ for all $l$ or $\{j^{\prime},i_2^{\prime}\}\neq \{j_{l_3}^{\prime},j_l\}$ for all $l$. By a similar argument as in (\ref{TTteq2}), $E_1=0$. The same result hods if $\{j,i_2\}=\{i^{\prime},j^{\prime}\}$. (ii) Suppose $\{j,i_1\}=\{j^{\prime},i_{m_1}^{\prime}\}$ or $\{j,i_1\}=\{i^{\prime},j^{\prime}\}$. By a similar argument as in the case $\{j,i_1\}=\{i^{\prime},j_{l_1}^{\prime}\}$, it is easy to get $E_1=0$.

Let $v_0\geq2,r_0\leq1$. Given $1\leq m\leq v_0$, if $\{j,i_m\}=\{i^{\prime},j_l^{\prime}\}$ for some $l$, then $E_1=0$. Hence $\{j,i_m\}=\{j^{\prime},i_l^{\prime}\}$ for some $l$. Then $j=j^{\prime}$ and 
\begin{equation}\label{TTeq2}
\{\{j,i\},\{j,i_1\},\dots,\{j,i_{v_0}\}\}\subset\{\{j^{\prime},i^{\prime}\},\{j^{\prime},i_1^{\prime}\},\dots,\{j^{\prime},i_{v}^{\prime}\}\}. 
\end{equation}
Similarly,
\begin{equation}\label{TTeq3}
\{\{j^{\prime},i^{\prime}\},\{j^{\prime},i_1^{\prime}\},\dots,\{j^{\prime},i_{v_0}^{\prime}\}\}\subset\{\{j,i\},\{j,i_1\},\dots,\{j,i_{v}\}\}.
\end{equation}
Without loss of generality, let $i=i_1^{\prime}$, $i_1=i^{\prime}$, $i_2=i_2^{\prime}$, $\dots$, $i_{v_1}^{\prime}=i_{v_1}^{\prime}$ for $v_0\leq v_1\leq v$.
There are at most $n^{2(r+v)+2-v_1}$ choices of these indices. If $r_0=0$, $E_1$ is bounded by $O\left(p_n^{2(r+v)+1-v_1}\right)$. Then the sum of (\ref{wedeq24}) over these indices is bounded by
\begin{equation}\label{TTeq3}
n(np_n)^{2(r+v)+1-v_1}f^{(s,t)}(np_n,np_n)^2=o(\sigma_n^2).
\end{equation}
If $r_0=1$, then $\{i,j_1\}=\{i^{\prime},j_1^{\prime}\}$. In this case, the sum of (\ref{wedeq24}) over these indices is bounded by
\begin{equation}\label{TTeq4}
(np_n)^{2(r+v)-v_1}f^{(s,t)}(np_n,np_n)^2=o(\sigma_n^2).
\end{equation}

Let $v_0=r_0=1$. In this case, $\{j,i_1\}=\{i^{\prime},j_l^{\prime}\}$ for some $l$ or $\{j,i_1\}=\{j^{\prime},i_m^{\prime}\}$ for some $m$ or $\{j,i_1\}=\{i^{\prime},j^{\prime}\}$. Let  $\{j,i_1\}=\{i^{\prime},j_{l_1}^{\prime}\}$. If $j=i^{\prime}$, then $\{i,j\}=\{i^{\prime},j_{l_2}^{\prime}\}$. In this case, $\{i^{\prime},j^{\prime}\}=\{j,i_{m_1}\}$ for some $m_1\geq2$, $\{i,j_1\}=\{j^{\prime},i_1^{\prime}\}$ and $l_1=1$. Otherwise $E_1=0$. There are at most $n^{2(r+v)-2}$ choices of these indices. In this case, the sum of (\ref{wedeq24}) over these indices is bounded by
\begin{equation}\label{TTeq5}
(np_n)^{2(r+v)-2}f^{(s,t)}(np_n,np_n)^2=o(\sigma_n^2).
\end{equation}
If $j=j_1^{\prime}$, then $\{i,j\}=\{j^{\prime},i_1^{\prime}\}$,  $\{i,j_1\}=\{i^{\prime},j^{\prime}\}$ and $l_1=1$. Otherwise $E_1=0$. There are at most $n^{2(r+v)-1}$ choices of these indices. In this case, the sum of (\ref{wedeq24}) over these indices is bounded by
\begin{equation}\label{TTeq5}
(np_n)^{2(r+v)-1}f^{(s,t)}(np_n,np_n)^2=o(\sigma_n^2).
\end{equation}
The cases $\{j,i_1\}=\{j^{\prime},i_m^{\prime}\}$ for some $m$ and $\{j,i_1\}=\{i^{\prime},j^{\prime}\}$ can be similarly bounded as in (\ref{TTeq4}) and (\ref{TTeq5}).

Let $v_0=1$ and $r_0=0$.  In this case, $\{j,i_1\}=\{i^{\prime},j_l^{\prime}\}$ for some $l$ or $\{j,i_1\}=\{j^{\prime},i_m^{\prime}\}$ for some $m$ or $\{j,i_1\}=\{i^{\prime},j^{\prime}\}$. Let  $\{j,i_1\}=\{j^{\prime},i_{m_1}^{\prime}\}$. If $m_1=1$, then $\{i,j\}=\{j^{\prime},i_{m_2}^{\prime}\}$ and $\{i^{\prime},j^{\prime}\}=\{j,i_{m_3}\}$. There are at most $n^{2(r+v)-1}$ choices of these indices. Then the sum of (\ref{wedeq24}) over these indices is bounded by (\ref{TTeq5}). If $m_1\geq2$ and $m_3=1$, then the sum of (\ref{wedeq24}) over these indices is bounded by (\ref{TTeq5}). If $m_1\geq2$ and $m_3\geq2$, then $\{j^{\prime},i_1^{\prime}\}=\{j,i_{m_3}\}$. There are at most $n^{2(r+v)-1}$ choices of these indices. Then the sum of (\ref{wedeq24}) over these indices is bounded by 
\begin{equation}\label{TTeq6}
n(np_n)^{2(r+v)-2}f^{(s,t)}(np_n,np_n)^2=o(\sigma_n^2).
\end{equation}
The case $\{j,i_1\}=\{i^{\prime},j_l^{\prime}\}$ for some $l$ and $\{j,i_1\}=\{i^{\prime},j^{\prime}\}$ can be similarly studied.

Suppose $\{i,j\}=\{i^{\prime},j^{\prime}\}$. Consider $i=i^{\prime}$ and $j=j^{\prime}$ first. For any $m$ with $1\leq m\leq v_0$, if $i_m\neq i_{m_1}^{\prime}$ for all $1\leq m_1\leq v$, then the expectation $E_1=0$. Hence, $\{i_{1},\dots,i_{v_0}\}\subset \{i_{1}^{\prime},\dots,i_{v}^{\prime}\}$. Similarly, $\{i_{1}^{\prime},\dots,i_{v_0}^{\prime}\}\subset \{i_{1},\dots,i_{v}\}$, $\{j_{1}^{\prime},\dots,j_{r_0}^{\prime}\}\subset \{j_{1},\dots,j_{r}\}$, $\{j_{1},\dots,j_{r_0}\}\subset \{j_{1}^{\prime},\dots,j_{r}^{\prime}\}$. Without loss of generality, let
\[\{i_{1}^{\prime},\dots,i_{v_1}^{\prime}\}=\{i_{1},\dots,i_{v_1}\}, \ \ \ \ \ \{i_{v_1+1}^{\prime},\dots,i_{v}^{\prime}\}\cap\{i_{v_1+1},\dots,i_{v}\}=\emptyset.\]
\[\{j_{1}^{\prime},\dots,j_{r_1}^{\prime}\}=\{j_{1},\dots,j_{r_1}\}, \ \ \ \ \ \{j_{r_1+1}^{\prime},\dots,j_{r}^{\prime}\}\cap\{j_{r_1+1},\dots,j_{r}\}=\emptyset,\]
where $v_0\leq v_1\leq v$, $r_0\leq r_1\leq r$. There are at most $n^{2+r_1+2(r-r_1)+v_1+2(v-v_1)}$ indices.
Let $\sigma_1$ be a one-to-one map from $\{i_{1},\dots,i_{v_1}\}$ to $\{i_{1}^{\prime},\dots,i_{v_1}^{\prime}\}$ and $\sigma_2$ be a one-to-one map from $\{j_{1},\dots,j_{r_1}\}$ to $\{j_{1}^{\prime},\dots,j_{r_1}^{\prime}\}$. Then
\begin{eqnarray}\nonumber
E_1
&=& \mathbb{E}\Bigg[\bar{A}_{ij}^2\left(\prod_{l=1}^{r_1}\bar{A}_{ij_l}^{\lambda_{r;l}+\lambda_{r;\sigma_2(l)}}\right)\left(\prod_{l=r_1+1}^{r}\bar{A}_{ij_l}^{\lambda_{r;l}}\bar{A}_{ij_l^{\prime}}^{\lambda_{r;l}}\right)\\ \nonumber
&&\times\left(\prod_{m=1}^{v_1}\bar{A}_{ji_m}^{\gamma_{v;m}+\gamma_{v;\sigma_1(m)}}\right)\left(\prod_{m=v_1+1}^{v}\bar{A}_{ji_m}^{\gamma_{v;m}}\bar{A}_{ji_m^{\prime}}^{\gamma_{v;m}}\right)\Bigg]\\ \label{thueq2}
&=&\Theta\left(p_n^{1+r_1+2(r-r_1)+v_1+2(v-v_1)}\right).
\end{eqnarray}
Then the sum over indices $\{i,j\}=\{i^{\prime},j^{\prime}\}$ in (\ref{wedeq24}) is bounded by
\begin{equation}\label{TTeq7}
n(np_n)^{2(r+v)-(r_1+v_1)+1}f^{(s,t)}(np_n,np_n)^2=o(\sigma_n^2).
\end{equation}

Similarly, (\ref{TTeq7}) holds for the case $i=j^{\prime}$ and $j=i^{\prime}$ with $\max\{r_0,v_0\}\leq r_1,v_1\leq \min\{r,v\}$.

Now we consider case (b). Suppose $\lambda_{r;l}\geq2$ for $l\in\{1,\dots,r\}$ and $\gamma_{v;m}\geq2$ for $m\in\{1,\dots,v\}$. In this case $r\leq\frac{s}{2}$ and $v\leq \frac{t}{2}$. If $\{i,j\}\neq\{i^{\prime},j^{\prime}\}$, $\{i,j\}\neq\{i^{\prime},j_l^{\prime}\}$ for all $l$, $\{i,j\}\neq\{j^{\prime},i_m^{\prime}\}$ for all $m$, then $E_1=0$. 

Suppose $\{i,j\}=\{i^{\prime},j^{\prime}\}$. There are two cases: $i=i^{\prime}$ and $j=j^{\prime}$ or $i=j^{\prime}$ and $j=i^{\prime}$. Let $i=i^{\prime}$ and $j=j^{\prime}$. Suppose $|\{j_1,\dots,j_r\}\cap\{j_1^{\prime},\dots,j_r^{\prime}\}|=r_1$ and $|\{i_1,\dots,i_v\}\cap\{i_1^{\prime},\dots,i_v^{\prime}\}|=v_1$. There are at most $n^{2+r_1+2(r-r_1)+v_1+2(v-v_1)}$ possible indices. Without loss of generality, assume $j_l=j_l^{\prime}$ for $1\leq l\leq r_1$ and $i_l=i_l^{\prime}$ for $1\leq l\leq v_1$. Then
\begin{eqnarray}\nonumber
E_1
&=&\mathbb{E}\Bigg[(A_{ij}-p_nw_{ij})^2\left(\prod_{l=1}^{r_1}\bar{A}_{ij_l}^{2\lambda_{r;l}}\right)\left(\prod_{l=r_1+1}^r\bar{A}_{ij_l}^{\lambda_{r;l}}\right)\left(\prod_{l=r_1+1}^r\bar{A}_{ij_l^{\prime}}^{\lambda_{r;l}}\right)\\ \nonumber
&&\times\left(\prod_{m=1}^{v_1}\bar{A}_{ji_m}^{2\gamma_{v;m}}\right)\left(\prod_{m=v_1+1}^{v}\bar{A}_{ji_m}^{\gamma_{v;m}}\right)\left(\prod_{m=v_1+1}^{v}\bar{A}_{ji_m^{\prime}}^{\gamma_{v;m}}\right)\Bigg]\\
&=&O\left(p_n^{1+r_1+2(r-r_1)+v_1+2(v-v_1)}\right).
\end{eqnarray}
Then the sum over indices $\{i,j\}=\{i^{\prime},j^{\prime}\}$ with $i=i^{\prime}$ and $j=j^{\prime}$ in (\ref{wedeq24}) is bounded by
\[n(np_n)^{2(r+v)-r_1-v_1+1}f^{(s,t)}(np_n,np_n)^2=o(\sigma_n^2).\]

Let $i=j^{\prime}$ and $j=i^{\prime}$. Suppose $|\{j_1,\dots,j_r\}\cap\{i_1^{\prime},\dots,i_v^{\prime}\}|=r_1$ and $|\{i_1,\dots,i_v\}\cap\{j_1^{\prime},\dots,j_r^{\prime}\}|=v_1$, where $0\leq r_1\leq \min\{r,v\}$ and $0\leq v_1\leq \min\{r,v\}$. There are at most $n^{2+r_1+(r-r_1)+(r-v_1)+v_1+(v-v_1)+(v-r_1)}$ such indices. Without loss of generality, let $j_l=i_l^{\prime}$ for $l\leq r_1$ and $i_l=j_l^{\prime}$ for $l\leq v_1$. Then
\begin{eqnarray}\nonumber
E_1
&=&\mathbb{E}\Bigg[(A_{ij}-p_nw_{ij})^2\left(\prod_{l=1}^{r_1}\bar{A}_{ij_l}^{\lambda_{r;l}+\gamma_{v;l}}\right)\left(\prod_{l=r_1+1}^r\bar{A}_{ij_l}^{\lambda_{r;l}}\right)\left(\prod_{l=v_1+1}^r\bar{A}_{jj_l^{\prime}}^{\lambda_{r;l}}\right)\\ \nonumber
&&\times\left(\prod_{m=1}^{v_1}\bar{A}_{ji_m}^{\gamma_{v;m}+\lambda_{r;m}}\right)\left(\prod_{m=v_1+1}^{v}\bar{A}_{ji_m}^{\gamma_{v;m}}\right)\left(\prod_{m=r_1+1}^{v}\bar{A}_{ii_m^{\prime}}^{\gamma_{v;m}}\right)\Bigg]\\
&=&\Theta\left(p_n^{1+r_1+(r-r_1)+(r-v_1)+v_1+(v-v_1)+(v-r_1)}\right).
\end{eqnarray}
 Then the sum over indices $\{i,j\}=\{i^{\prime},j^{\prime}\}$ with $i=j^{\prime}$ and $j=i^{\prime}$ in (\ref{wedeq24}) is bounded by
\[O\left(n(np_n)^{1+2(r+v)-(r_1+v_1)}f^{(s,t)}(np_n,np_n)^2\right)=o(\sigma_n^2).\]

Suppose $\{i,j\}=\{i^{\prime},j_{l_1}^{\prime}\}$ for some $l_1$. If $i=j_{l_1}^{\prime}$ and $j=i^{\prime}$, then $j^{\prime}=i_{m_1}$ for some $m_1$, otherwise $E_1=0$. There are at most $n^{3+r+(r-1)+v+(v-1)}$ possible nodes. In this case,
\begin{eqnarray}\nonumber
E_1
&=&\mathbb{E}\Bigg[(A_{jj_{l_1}^{\prime}}-p_nw_{jj_{l_1}^{\prime}})^{1+\lambda_{r;l_1}}\left(\prod_{l=1}^{r}\bar{A}_{ij_l}^{\lambda_{r;l}}\right)\left(\prod_{l=1,l\neq l_1}^r\bar{A}_{jj_l^{\prime}}^{\lambda_{r;l}}\right)\\ \nonumber
&&\times(A_{ji_{m_1}}-p_nw_{ji_{m_1}})^{1+\gamma_{v;m_1}}\left(\prod_{m=1,m\neq m_1}^{v}\bar{A}_{ji_m}^{\gamma_{v;m}}\right)\left(\prod_{m=1}^{v}\bar{A}_{i_{m_1}i_m^{\prime}}^{\gamma_{v;m}}\right)\Bigg]\\
&=&\Theta\left(p_n^{2+r+(r-1)+v+(v-1)}\right).
\end{eqnarray}
 Then the sum over indices $\{i,j\}=\{i^{\prime},j_{l_1}^{\prime}\}$ in (\ref{wedeq24}) is bounded by
\begin{equation}\label{TTeq8}
O\left(n(np_n)^{2(r+v)}f^{(s,t)}(np_n,np_n)^2\right)=o(\sigma_n^2).
\end{equation}
Similarly, (\ref{TTeq8}) holds for $i=i^{\prime}$ and $j=j_{l_1}^{\prime}$ or $\{i,j\}=\{j^{\prime},i_{m_1}^{\prime}\}$ for some $m_1$.

\subsubsection{Bound the last term of (\ref{iiineq1})}

Now we prove the last term of (\ref{iiineq1}) converges in probability to zero. To this end, we will show that 
\begin{equation}\label{apeq82}
\mathbb{E}\left[\left|\sum_{ i\neq j}R_{ij}A_{ij}\right|\right]=o\left(\sigma_n\right). 
\end{equation}

Let $s,t$ satisfy $s+t=k_0$. By $(C3)$ of Assumption \ref{assump}, $|f^{(s,t)}(x,y)|$ is monotone in $x$ and $y$. There are four cases: $|f^{(s,t)}(x,y)|$ is decreasing in $x$ and $y$,  $|f^{(s,t)}(x,y)|$ is increasing in $x$ and $y$, $|f^{(s,t)}(x,y)|$ is increasing in $x$ and decreasing in $y$, $|f^{(s,t)}(x,y)|$ is decreasing in $x$ and increasing in $y$.

Suppose $|f^{(s,t)}(x,y)|$ is decreasing in $x$ and $y$.
Let $\delta_n=\left[\log (np_n)\right]^{-2}$. Then
\begin{eqnarray}\nonumber
&&\mathbb{E}[|f^{(s,t)}(X_{i(j)},X_{j(i)})||d_{i(j)}-w_{i(j)}|^s|d_{j(i)}-w_{j(i)}|^t] \\ \nonumber
&=&\mathbb{E}\Bigg[|f^{(s,t)}(X_{i(j)},X_{j(i)})||d_{i(j)}-w_{i(j)}|^s|d_{j(i)}-w_{j(i)}|^t\\ \nonumber
&&\times I[X_{i(j)}\geq \delta_nw_{i(j)},X_{j(i)}\geq \delta_nw_{j(i)}]\Bigg]\\ \nonumber
&&+\mathbb{E}\Bigg[|f^{(s,t)}(X_{i(j)},X_{j(i)})||d_{i(j)}-w_{i(j)}|^s|d_{j(i)}-w_{j(i)}|^t\\ \nonumber
&&\times I[X_{i(j)}\geq \delta_nw_{i(j)},X_{j(i)}< \delta_nw_{j(i)}]\Bigg]\\ \nonumber 
&&+\mathbb{E}\Bigg[|f^{(s,t)}(X_{i(j)},X_{j(i)})||d_{i(j)}-w_{i(j)}|^s|d_{j(i)}-w_{j(i)}|^t\\ \nonumber
&&\times I[X_{i(j)}< \delta_nw_{i(j)},X_{j(i)}\geq \delta_nw_{j(i)}]\Bigg]\\  \nonumber
&&+\mathbb{E}\Bigg[|f^{(s,t)}(X_{i(j)},X_{j(i)})||d_{i(j)}-w_{i(j)}|^s|d_{j(i)}-w_{j(i)}|^t\\\label{eqn2}
&&\times I[X_{i(j)}< \delta_nw_{i(j)},X_{j(i)}< \delta_nw_{j(i)}]\Bigg]. 
\end{eqnarray}
By the Cauchy–Schwarz inequality and (\ref{mmmeq1}), we have
\begin{eqnarray*}  \mathbb{E}\left[|d_{i(j)}-w_{i(j)}|^s\right]&\leq& \sqrt{\mathbb{E}\left[(d_{i(j)}-w_{i(j)})^{2s}\right]}\\
&=&\sqrt{\sum_{r=1}^{2s}\sum_{\substack{j_1,j_2,\dots,j_r\notin\{i,j\}\\ j_1\neq j_2\neq \dots\neq j_r}}\prod_{l=1}^r\mathbb{E}\left[(A_{ij_l}-p_nw_{ij_l})^{\lambda_{r;l}}\right]}\\
&=&O\left(\sqrt{\sum_{r=1}^s(np_n)^r}\right)=O\left((np_n)^{\frac{s}{2}}\right).
\end{eqnarray*}
Then
\begin{eqnarray}\nonumber
&&\mathbb{E}\Bigg[|f^{(s,t)}(X_{i(j)},X_{j(i)})||d_{i(j)}-w_{i(j)}|^s|d_{j(i)}-w_{j(i)}|^t\\ \nonumber
&&\times I[X_{i(j)}\geq \delta_nw_{i(j)},X_{j(i)}\geq \delta_nw_{j(i)}]\Bigg]\\ \nonumber
&\leq& \mathbb{E}\Bigg[|f^{(s,t)}(\delta_nw_{i(j)},\delta_nw_{j(i)})|d_{i(j)}-w_{i(j)}|^s|d_{j(i)}-w_{j(i)}|^t\\ \nonumber
&&\times I[X_{i(j)}\geq \delta_nw_{i(j)},X_{j(i)}\geq \delta_nw_{j(i)}]\Bigg]\\ \nonumber
&\leq&|f^{(s,t)}(\delta_nw_{i(j)},\delta_nw_{j(i)})|\mathbb{E}\left[|d_{i(j)}-w_{i(j)}|^s|d_{j(i)}-w_{j(i)}|^t\right]\\ \nonumber
&=&|f^{(s,t)}(\delta_nw_{i(j)},\delta_nw_{j(i)})|\mathbb{E}\left[|d_{i(j)}-w_{i(j)}|^s\right]\mathbb{E}\left[|d_{j(i)}-w_{j(i)}|^t\right]\\ \label{eqn3}
&\leq&(np_n)^{\frac{k_0}{2}}|f^{(s,t)}(\delta_nw_{i(j)},\delta_nw_{j(i)})|.
\end{eqnarray}

On the event $\{X_{i(j)}< \delta_nw_{i(j)}\}$, if $X_{i(j)}<d_{i(j)}$, then $X_{i(j)}$ cannot be between $d_{i(j)}$ and $w_{i(j)}$. Hence $X_{i(j)}< \delta_nw_{i(j)}$ implies $d_{i(j)}\leq X_{i(j)}< \delta_nw_{i(j)}$. Similar result holds for $X_{j(i)}$. By definition, $d_{i(j)}$ and $d_{j(i)}$ are independent if $i\neq j$. By Lemma \ref{lem1} and $(C2)$ of Assumption \ref{assump}, one has
\begin{eqnarray}\nonumber
&&\mathbb{E}\Big[|f^{(s,t)}(X_{i(j)},X_{j(i)})||d_{i(j)}-w_{i(j)}|^s|d_{j(i)}-w_{j(i)}|^t\\ \nonumber
&&\times I[X_{i(j)}< \delta_nw_{i(j)},X_{j(i)}< \delta_nw_{j(i)}]\Big]\\ \nonumber
&\leq& \mathbb{E}\Big[|f^{(s,t)}(X_{i(j)},X_{j(i)})||d_{i(j)}-w_{i(j)}|^s|d_{j(i)}-w_{j(i)}|^t\\ \nonumber
&&\times I[d_{i(j)}\leq X_{i(j)}< \delta_nw_{i(j)},d_{j(i)}\leq X_{j(i)}< \delta_nw_{j(i)}]\Big]\\ \nonumber
&\leq&\mathbb{E}\Big[|f^{(s,t)}(d_{i(j)},d_{j(i)})||d_{i(j)}-w_{i(j)}|^s|d_{j(i)}-w_{j(i)}|^t\\ \nonumber
&&\times I[d_{i(j)}< \delta_nw_{i(j)},d_{j(i)}< \delta_nw_{j(i)}]\Big]\\ \nonumber
&=&\sum_{k=1}^{\delta_nw_{i(j)}}\sum_{l=1}^{\delta_nw_{j(i)}}|f^{(s,t)}(k,l)||k-w_{i(j)}|^s|l-w_{j(i)}|^t\mathbb{P}(d_{i(j)}=k)\mathbb{P}(d_{j(i)}=l)\\ \nonumber
&=&O\left((\delta_nnp_n)^M\exp(-2np_n\beta(1+o(1))\right)\\ \label{eqn4}
&=&\exp(-2np_n\beta(1+o(1)),
\end{eqnarray}
where $M$ is a positive constant.

Similarly, the second term of (\ref{eqn2}) is bounded as follows.
\begin{eqnarray}\nonumber
&&\mathbb{E}\Big[|f^{(s,t)}(X_{i(j)},X_{j(i)})||d_{i(j)}-w_{i(j)}|^s|d_{j(i)}-w_{j(i)}|^t\\ \nonumber
&&\times I[X_{i(j)}\geq \delta_nw_{i(j)},X_{j(i)}< \delta_nw_{j(i)}]\Big]\\ \nonumber
&\leq&\mathbb{E}\Big[|f^{(s,t)}(\delta_nw_{i(j)},X_{j(i)})||d_{i(j)}-w_{i(j)}|^s|d_{j(i)}-w_{j(i)}|^t\\ \nonumber
&&\times I[X_{i(j)}\geq \delta_nw_{i(j)},d_{j(i)}\leq X_{j(i)}< \delta_nw_{j(i)}]\Big]\\ \nonumber
&\leq& \mathbb{E}\Big[|f^{(s,t)}(\delta_nw_{i(j)},d_{j(i)})||d_{i(j)}-w_{i(j)}|^s|d_{j(i)}-w_{j(i)}|^t\\ \nonumber
&&\times I[d_{j(i)}< \delta_nw_{j(i)}]\Big]\\ \nonumber
&=&\sum_{k=1}^{\delta_nw_{j(i)}}|f^{(s,t)}(\delta_nw_{i(j)},k)|(k-w_{j(i)})^t\mathbb{P}(d_{j(i)}=k)\mathbb{E}[|d_{i(j)}-w_{i(j)}|^s]\\ \nonumber
&=&O\left((\delta_nnp_n)^M\exp(-np_n\beta(1+o(1))\right)\\ \label{eqn5}
&=&\exp(-np_n\beta(1+o(1)).
\end{eqnarray}
The third term of (\ref{eqn2}) can be similarly bounded.

Combining (\ref{eqn2})-(\ref{eqn5}) yields
\begin{eqnarray*} 
&&\mathbb{E}\left[|f^{(s,t)}(X_{i(j)},X_{j(i)})||d_{i(j)}-w_{i(j)}|^s|d_{j(i)}-w_{j(i)}|^t\right] \\
&=&O\left((np_n)^{\frac{k_0}{2}}|f^{(s,t)}(\delta_nw_{i(j)},\delta_nw_{j(i)})|+\exp(-np_n\beta(1+o(1))\right).
\end{eqnarray*}
By $(C4)$ of Assumption \ref{assump}, we have
\begin{equation}\nonumber
\mathbb{E}\left[\left|\sum_{ i\neq j}R_{ij}A_{ij}\right|\right]=O\left(n(np_n)^{\frac{k_0}{2}+1}|f^{(s,t)}(\delta_nnp_n,\delta_nnp_n)|\right)=o\left(\sigma_n\right). 
\end{equation}
Then (\ref{apeq82}) holds.

Suppose $|f^{(s,t)}(x,y)|$ is increasing in $x$ and $y$.
Let $M$ be a large positive constant. Then
\begin{eqnarray}\nonumber
&&\mathbb{E}\left[|f^{(s,t)}(X_{i(j)},X_{j(i)})||d_{i(j)}-w_{i(j)}|^s|d_{j(i)}-w_{j(i)}|^t\right] \\ \nonumber
&=&\mathbb{E}\Bigg[|f^{(s,t)}(X_{i(j)},X_{j(i)})||d_{i(j)}-w_{i(j)}|^s|d_{j(i)}-w_{j(i)}|^t\\ \nonumber
&&\times I[X_{i(j)}\geq Mw_{i(j)},X_{j(i)}\geq Mw_{j(i)}]\Bigg]\\ \nonumber
&&+\mathbb{E}\Bigg[|f^{(s,t)}(X_{i(j)},X_{j(i)})||d_{i(j)}-w_{i(j)}|^s|d_{j(i)}-w_{j(i)}|^t\\ \nonumber
&&\times I[X_{i(j)}\geq Mw_{i(j)},X_{j(i)}< Mw_{j(i)}]\Bigg]\\ \nonumber 
&&+\mathbb{E}\Bigg[|f^{(s,t)}(X_{i(j)},X_{j(i)})||d_{i(j)}-w_{i(j)}|^s|d_{j(i)}-w_{j(i)}|^t\\ \nonumber
&&\times I[X_{i(j)}< Mw_{i(j)},X_{j(i)}\geq Mw_{j(i)}]\Bigg]\\ \nonumber
&&+\mathbb{E}\Bigg[|f^{(s,t)}(X_{i(j)},X_{j(i)})||d_{i(j)}-w_{i(j)}|^s|d_{j(i)}-w_{j(i)}|^t\\ \nonumber
&&\times I[X_{i(j)}< Mw_{i(j)},X_{j(i)}< Mw_{j(i)}]\Bigg].\\ \label{tueq4}
\end{eqnarray}

On the event $\{X_{i(j)}\geq Mw_{i(j)}\}$, if $X_{i(j)}>d_{i(j)}$, then $X_{i(j)}$ cannot be between $d_{i(j)}$ and $w_{i(j)}$. Hence $X_{i(j)}\geq  Mw_{i(j)}$ implies $Mw_{i(j)}\leq X_{i(j)}\leq d_{i(j)}$. Similar result holds for $X_{j(i)}$. By definition, $d_{i(j)}$ and $d_{j(i)}$ are independent if $i\neq j$. Suppose $np_n=\omega(\log n)$. By Lemma \ref{lem1} and $(C2)$ of Assumption \ref{assump}, one has
\begin{eqnarray}\nonumber
&&\mathbb{E}\Bigg[|f^{(s,t)}(X_{i(j)},X_{j(i)})||d_{i(j)}-w_{i(j)}|^s|d_{j(i)}-w_{j(i)}|^t\\ \nonumber
&&\times I[X_{i(j)}\geq Mw_{i(j)},X_{j(i)}\geq Mw_{j(i)}]\Bigg]\\ \nonumber
&\leq&\mathbb{E}\Bigg[|f^{(s,t)}(X_{i(j)},X_{j(i)})||d_{i(j)}-w_{i(j)}|^s|d_{j(i)}-w_{j(i)}|^t\\ \nonumber
&&\times I[d_{i(j)}\geq X_{i(j)}\geq Mw_{i(j)},d_{j(i)}\geq X_{j(i)}\geq Mw_{j(i)}]\Bigg]\\ \nonumber
&\leq&\mathbb{E}\Bigg[|f^{(s,t)}(d_{i(j)},d_{j(i)})||d_{i(j)}-w_{i(j)}|^s|d_{j(i)}-w_{j(i)}|^t\\ \nonumber
&&\times  I[d_{i(j)}-1\geq Mw_{i(j)}-1,d_{j(i)}-1\geq Mw_{j(i)}-1]\Bigg]\\  \nonumber
&=&\sum_{k=Mw_{i(j)}-1}^{n-2}\sum_{l=Mw_{j(i)}-1}^{n-2}|f^{(s,t)}(k+1,l+1)||k+1-w_{i(j)}|^s|l+1-w_{j(i)}|^t\\ \nonumber
&&\times \mathbb{P}(d_{i(j)}-1=k)\mathbb{P}(d_{i(j)}-1=l)\\ \nonumber
&=&O\left(n^{C_{s,t,f}}e^{-np_n\beta(1+o(1))}\right)\\ \label{tueq1}
&=&O\left(e^{-np_n\beta(1+o(1))}\right),
\end{eqnarray}
where $C_{s,t,f}$ is some constant dependent on $s,t$ and $f$. Similarly,
\begin{eqnarray}   \nonumber 
&&\mathbb{E}\Bigg[|f^{(s,t)}(X_{i(j)},X_{j(i)})||d_{i(j)}-w_{i(j)}|^s|d_{j(i)}-w_{j(i)}|^t\\ \nonumber
&&\times I[X_{i(j)}\geq Mw_{i(j)},X_{j(i)}< Mw_{j(i)}]\Bigg]\\ \nonumber
&\leq&\mathbb{E}\Bigg[|f^{(s,t)}(X_{i(j)},X_{j(i)})||d_{i(j)}-w_{i(j)}|^s|d_{j(i)}-w_{j(i)}|^t\\ \nonumber
&&\times I[d_{i(j)}\geq X_{i(j)}\geq Mw_{i(j)},X_{j(i)}< Mw_{j(i)}]\Bigg]\\ \nonumber
&\leq&\mathbb{E}\Big[|f^{(s,t)}(d_{i(j)},Mw_{j(i)})||d_{i(j)}-w_{i(j)}|^s|d_{j(i)}-w_{j(i)}|^t\\ \nonumber
&&\times I[d_{i(j)}\geq Mw_{i(j)}]\Big]\\ \nonumber
&=&\mathbb{E}|d_{j(i)}-w_{j(i)}|^t\sum_{k=Mw_{i(j)}-1}^{n-2}|f^{(s,t)}(k+1,Mw_{j(i)})||k+1-w_{i(j)}|^s\\ \nonumber
&&\times\mathbb{P}(d_{j(i)}-1=k)\\ \label{tueq2}
&=&O\left(e^{-np_n\beta(1+o(1))}\right),
\end{eqnarray}
and
\begin{eqnarray}\nonumber
&&\mathbb{E}\Bigg[|f^{(s,t)}(X_{i(j)},X_{j(i)})||d_{i(j)}-w_{i(j)}|^s|d_{j(i)}-w_{j(i)}|^t\\ \nonumber
&& \times I[X_{i(j)}< Mw_{i(j)},X_{j(i)}< Mw_{j(i)}]\Bigg]\\ \nonumber
&\leq&\mathbb{E}\left[|f^{(s,t)}(Mw_{i(j)},Mw_{j(i)})||d_{i(j)}-w_{i(j)}|^s|d_{j(i)}-w_{j(i)}|^t\right]\\ \nonumber
&=&\mathbb{E}\left[|f^{(s,t)}(Mw_{i(j)},Mw_{j(i)})|\right]\mathbb{E}\left[|d_{i(j)}-w_{i(j)}|^s\right]\mathbb{E}\left[|d_{j(i)}-w_{j(i)}|^t\right]\\ \label{tueq3}
&=&O\left((np_n)^{\frac{k_0}{2}}|f^{(s,t)}(Mnp_n,Mnp_n)|\right),
\end{eqnarray}

By (\ref{tueq4}), (\ref{tueq1}),(\ref{tueq2}) and (\ref{tueq3}), it follows that
\begin{equation*}\label{tueq5}
\mathbb{E}\left[\left|\sum_{ i\neq j}R_{ij}A_{ij}\right|\right]=O\left(n(np_n)^{\frac{k_0}{2}+1}|f^{(s,t)}(Mnp_n,Mnp_n)|\right)=o\left(\sigma_n\right). 
\end{equation*}
Then (\ref{apeq82}) holds.

Suppose $|f^{(s,t)}(x,y)|$ is decreasing in $x$ and increasing in $y$. Let $\delta_n=\left[\log (np_n)\right]^{-2}$ and $M$ be a large positive constant. Then
\begin{eqnarray}\nonumber
&&\mathbb{E}\left[|f^{(s,t)}(X_{i(j)},X_{j(i)})||d_{i(j)}-w_{i(j)}|^s|d_{j(i)}-w_{j(i)}|^t\right] \\ \nonumber
&=&\mathbb{E}\Bigg[|f^{(s,t)}(X_{i(j)},X_{j(i)})||d_{i(j)}-w_{i(j)}|^s|d_{j(i)}-w_{j(i)}|^t\\ \nonumber
&& \times I[X_{i(j)}\geq \delta_n w_{i(j)},X_{j(i)}\geq Mw_{j(i)}]\Bigg]\\ \nonumber
&&+\mathbb{E}\Bigg[|f^{(s,t)}(X_{i(j)},X_{j(i)})||d_{i(j)}-w_{i(j)}|^s|d_{j(i)}-w_{j(i)}|^t\\ \nonumber
&& \times I[X_{i(j)}\geq \delta_nw_{i(j)},X_{j(i)}< Mw_{j(i)}]\Bigg]\\ \nonumber 
&&+\mathbb{E}\Bigg[|f^{(s,t)}(X_{i(j)},X_{j(i)})||d_{i(j)}-w_{i(j)}|^s|d_{j(i)}-w_{j(i)}|^t\\ \nonumber
&& \times I[X_{i(j)}< \delta_nw_{i(j)},X_{j(i)}\geq Mw_{j(i)}]\Bigg]\\ \nonumber
&&+\mathbb{E}\Bigg[|f^{(s,t)}(X_{i(j)},X_{j(i)})||d_{i(j)}-w_{i(j)}|^s|d_{j(i)}-w_{j(i)}|^t\\ \nonumber
&& \times I[X_{i(j)}< \delta_nw_{i(j)},X_{j(i)}< Mw_{j(i)}]\Bigg].\\ \label{phtueq1}
\end{eqnarray}
The first term of 
(\ref{phtueq1}) can be bounded by
\begin{eqnarray}\nonumber
&&\mathbb{E}\Bigg[|f^{(s,t)}(X_{i(j)},X_{j(i)})||d_{i(j)}-w_{i(j)}|^s|d_{j(i)}-w_{j(i)}|^t\\ \nonumber
&&\times I[X_{i(j)}\geq \delta_n w_{i(j)},X_{j(i)}\geq Mw_{j(i)}]\Bigg]\\ \nonumber
&\leq&\mathbb{E}\Bigg[|f^{(s,t)}(\delta_n w_{i(j)},d_{j(i)})||d_{i(j)}-w_{i(j)}|^s|d_{j(i)}-w_{j(i)}|^t\\ \nonumber
&&\times I[X_{i(j)}\geq \delta_n w_{i(j)},d_{j(i)}\geq Mw_{j(i)}]\Bigg]\\ \nonumber
&=&(np_n)^{
\frac{s}{2}
}\sum_{k=Mw_{j(i)}}^{n-2}|f^{(s,t)}(\delta_n w_{i(j)},k)||k-w_{j(i)}|^t\mathbb{P}(d_{j(i)}=k)\\ \label{mmeq3}
&=&\exp(-np_n\beta(1+o(1)),
\end{eqnarray}
The second term of 
(\ref{phtueq1})  can be bounded by
\begin{eqnarray}\nonumber
&&\mathbb{E}\Bigg[|f^{(s,t)}(X_{i(j)},X_{j(i)})||d_{i(j)}-w_{i(j)}|^s|d_{j(i)}-w_{j(i)}|^t\\ \nonumber
&&\times I[X_{i(j)}\geq \delta_nw_{i(j)},X_{j(i)}< Mw_{j(i)}]\Bigg]\\ \label{mmeq2}
&\leq&|f^{(s,t)}(\delta_nw_{i(j)},Mw_{j(i)})|(np_n)^{\frac{k_0}{2}}.
\end{eqnarray}
The third term of 
(\ref{phtueq1})  can be bounded by
\begin{eqnarray}\nonumber
&&\mathbb{E}\Bigg[|f^{(s,t)}(X_{i(j)},X_{j(i)})||d_{i(j)}-w_{i(j)}|^s|d_{j(i)}-w_{j(i)}|^t\\ \nonumber
&&\times I[X_{i(j)}< \delta_nw_{i(j)},X_{j(i)}\geq Mw_{j(i)}]\Bigg]\\ \nonumber
&\leq&\mathbb{E}\Bigg[|f^{(s,t)}(d_{i(j)},d_{j(i)})||d_{i(j)}-w_{i(j)}|^s|d_{j(i)}-w_{j(i)}|^t\\ \nonumber
&&\times I[d_{i(j)}< \delta_nw_{i(j)},d_{j(i)}\geq Mw_{j(i)}]\Bigg]\\ \nonumber
&=&\sum_{k=0}^{\delta_nw_{i(j)}}\sum_{l=Mw_{j(i)}}^{n-2}|f^{(s,t)}(k,l)||k-w_{i(j)}|^s|l-w_{j(i)}|^t\\ \nonumber
&&\times\mathbb{P}(d_{i(j)}=k)\mathbb{P}(d_{j(i)}=l)\\ \label{mmeq1}
&=&\exp(-np_n\beta(1+o(1)).
\end{eqnarray}
By (\ref{phtueq1})-(\ref{mmeq1}),  (\ref{apeq82}) holds.

The case that $|f^{(s,t)}(x,y)|$ is increasing in $x$ and decreasing in $y$ can be similarly processed. We omit it. Then the proof is complete.

\qed

\subsection{\bf Proof of Theorem \ref{geRandic} in Subsection \ref{geR}} 
To prove Theorem \ref{geRandic} in Subsection \ref{geR}, we only need to derive the asymptotic distribution of the Randi\'{c} index of the Erd\H{o}s-R\'{e}nyi random graph $\mathcal{G}_n(\alpha)$, that is, $\tau=-\frac{1}{2}$.
Recall that for the Erd\H{o}s-R\'{e}nyi random graph $\mathcal{G}_n(\alpha)$, $w_{i(k)}=1+(n-2)p_n$. Then
\[f_x(w_{i(j)},w_{j(i)})=f_y(w_{i(j)},w_{j(i)})=-\frac{1}{2(1+(n-2)p_n))^2},\]
\[f_{xx}(w_{i(j)},w_{j(i)})=f_{yy}(w_{i(j)},w_{j(i)})=\frac{3}{4(1+(n-2)p_n))^3},\] 
\[ f_{xy}(w_{i(j)},w_{j(i)})=\frac{1}{4(1+(n-2)p_n))^3},\]
 \[ |f^{(s,t)}(np_n,np_n)|=O\left(\frac{1}{(np_n)^{1+(s+t)}}\right).\]
 Let $k_0=\max\left\{\lfloor2+\frac{1}{1-\alpha}\rfloor+1,3\right\}$. By (\ref{iineq1}) and the proof of Theorem \ref{mainthm}, we have
\begin{eqnarray*}\nonumber
\mathcal{I}_n-\mathbb{E}[\mathcal{I}_n]
&=&\frac{1}{2}\sum_{ i\neq j}(M_{ij}A_{ij}-\mathbb{E}[M_{ij}A_{ij}])+\frac{1}{2}\sum_{ i\neq j}(S_{ij}A_{ij}-\mathbb{E}[S_{ij}A_{ij}])\\ \label{ineq1}
&&+O_P\left(\sqrt{\frac{n}{(np_n)^3}}\right). \
\end{eqnarray*}

By the calculations in Section \ref{geR}, $a_{ij}=\Theta\left(\frac{1}{(np_n)^{2}}\right)$ for $\tau=-\frac{1}{2}$. Then
\[Var\left(\sum_{ i\neq j}a_{ij}(A_{ij}-p_nw_{ij})\right)=O\left(\frac{n}{(np_n)^{3}}\right).\]

By equations (\ref{apeq2})- (\ref{fyineq3}), we have
\begin{eqnarray}\nonumber
&&\frac{1}{2}\sum_{ i\neq j}\left(M_{ij}A_{ij}-\mathbb{E}[M_{ij}A_{ij}]\right)\\ \nonumber
&=&-\frac{\sum_{i\neq j\neq l}(A_{il}-p_n)(A_{ij}-p_n)}{2(1+(n-2)p_n))^2} +O_P\left(\sqrt{\frac{n}{(np_n)^3}}\right).
\end{eqnarray}
By equations (\ref{apeq10})-(\ref{caneq3}), we have
\begin{eqnarray}\nonumber
&&\frac{1}{2}\sum_{ i\neq j}\left(S_{ij}A_{ij}-\mathbb{E}[S_{ij}A_{ij}]\right)\\ \nonumber
&=&\frac{1}{4}\sum_{\substack{i\neq j,
s\neq t\\s,t\notin\{i,j\}}}p_n\left(f_{xx}(w_{i(j)},w_{j(i)})+f_{yy}(w_{i(j)},w_{j(i)})\right)\\ \nonumber
&&\times (A_{is}-p_n)(A_{it}-p_n)+O_P\left(\sqrt{\frac{n}{(np_n)^3}}\right)
\\ \nonumber
&=&\frac{3(n-3)p_n}{8(1+(n-2)p_n))^3}\sum_{\substack{s\neq t\neq i}}(A_{is}-p_n)(A_{it}-p_n)+O_P\left(\sqrt{\frac{n}{(np_n)^3}}\right).
\end{eqnarray}
Hence, we get
\begin{eqnarray}\label{ineq1}
\mathcal{I}_n -\mathbb{E}[\mathcal{I}_n ]
&=&\mathcal{X}_n+O_P\left(\sqrt{\frac{n}{(np_n)^3}}\right), 
\end{eqnarray}
where
\[\mathcal{X}_n=-\frac{1}{8(1+(n-2)p_n))^2}\sum_{\substack{s\neq t\neq i}}(A_{is}-p_n)(A_{it}-p_n).\]
Note that
\begin{eqnarray*}
-\mathcal{X}_n
&=&\sum_{\substack{i<j<k}}\frac{(A_{ij}-p_n)(A_{ik}-p_n)}{4(1+(n-2)p_n))^2}+\sum_{\substack{i<j<k}}\frac{(A_{ji}-p_n)(A_{jk}-p_n)}{4(1+(n-2)p_n))^2}\\
&&+\sum_{\substack{i<j<k}}\frac{(A_{ki}-p_n)(A_{kj}-p_n)}{4(1+(n-2)p_n))^2},
\end{eqnarray*}
and the variance $\sigma_n^2$ of $\mathcal{X}_n$ is equal to
\[\sigma_n^2=\frac{n(n-1)(n-2)p_n^2(1-p_n)^2}{32(1+(n-2)p_n))^4}=\Theta\left(\frac{n}{(np_n)^2}\right).\]
By Theorem 6.1 in \cite{GL17b}, we have
\[\frac{\mathcal{X}_n}{\sigma_n}\Rightarrow \mathcal{N}(0,1),\]
from which and (\ref{ineq1}) it follows that
\[\frac{\mathcal{I}_n -\mathbb{E}[\mathcal{I}_n]}{\sigma_n}\Rightarrow \mathcal{N}(0,1).\]
Then the proof is complete.

\qed

\acknowledgment{The author is grateful to the anonymous reviewers for valuable comments.}

\singlespacing


\begin{thebibliography}{9}

\bibitem{AZG19}
A. Ali, L. Zhong, I. Gutman, 
Harmonic index and its generalizations:
extremal results and bounds,
\textit{MATCH Commun. Math. Comput. Chem.} {\bf 81} (2019) 249-311.


\bibitem{BE98}
B. Bollob\'{a}s, P. Erd\H{o}s, Graphs of extremal weights, \textit{Ars Comb.} {\bf 50} (1998) 225-233.


\bibitem{BES99}
B. Bollobás, P. Erd\H{o}s, A. Sarkar,
Extremal graphs for weights,
\textit{Discrete Math.} {\bf 200} (1999) 5-19.




\bibitem{BBGG00}
S. C. Basak, A. T. Balaban, G. Grunwald, B. D. Gute,
Topological Indices: Their Nature and Mutual Relatedness, \textit{J. Chem. Inf. Comput. Sci.} {\bf 40} (2000) 891-898.




\bibitem{CFK10}
M. Cavers, S.  Kirkland, On the normalized Laplacian energy and general
Randi\'{c} index $r_1$ of graphs, \textit{Lin. Algebra Appl.}  {\bf433} (2010) 172-190.



\bibitem{DSG17}
K.C. Das, S. Sun, I. Gutman,
Normalized Laplacian eigenvalues and Randi\'c energy of graphs,
\textit{MATCH Commun. Math. Comput. Chem.} {\bf 77} (2017) 45-59






\bibitem{DHIR20}
T. Do\v{s}li\'c, M. A. Hosseinzadeh, S. Hossein-Zadeh, A. Iranmanesh, F. Rezakhanlou,
On generalized Zagreb indices
of random graphs, \textit{MATCH Commun. Math. Comput. Chem. } {\bf 84} (2020) 499-511.


\bibitem{DMRSV18}
P. De Meo, F. Messina, D. Rosaci, G. M. L. Sarne, A. V. Vasilakos,
Estimating graph robustness
through the Randi\'{c} index,
\textit{IEEE Trans Cybern. } {\bf 48} (2018) 3232-3242.





\bibitem{DMMRSF21}
S. Dattola, N. Mammone, F. C. Morabito, D. Rosaci, G.M.L. Sarne, F. L. Foresta, Testing graph robustness indexes for EEG analysis in
Alzheimer’s disease diagnosis, \textit{Electronics} {\bf 10} (2021) 1440.







\bibitem{E10}
E. Estrada, 
Quantifying network heterogeneity, \textit{Physical Review E }  {\bf 82} (2010) 066102.




\bibitem{FT13}
D. Fourches, A. Tropsha,  Using graph indices for the analysis and comparison of chemical datasets, \textit{Molecular Informatics}  {\bf 32} (2013)  827-842.








\bibitem{G13}
I. Gutman,
Degree-based topological indices, \textit{Croatica Chemica Acta}  {\bf 86} (2013) 351-361.



\bibitem{GL17b} 
C. Gao, J. Lafferty,  
Testing network structure using relations between small subgraph probabilities,
\url{https://arxiv.org/pdf/1704.06742.pdf}





\bibitem{LSG21}
S. Li, L. Shi, W. Gao, 
Two modified Zagreb indices for random
structures, \textit{Main Group Met. Chem.}  {\bf 44} (2021) 150–156.



\bibitem{LS08}
X. Li, Y.  Shi,  A survey on the Randi\'c index,
\textit{MATCH Commun. Math. Comput. Chem.}  {\bf 59} (2008) 127-156.










\bibitem{MCSGDF18}
Y. Ma, S. Cao, Y. Shi, I. Gutman, M. Dehmer, B. Furtula,
From the connectivity index to
various Randi\'c-type descriptors, \textit{MATCH Commun. Math. Comput. Chem.} {\bf 80} (2018) 85-106.


\bibitem{MMRS20}
C. T. Mart\'inez-Mart\'inez, J. A. M\'endez-Berm\'udez,   J. M. Rodriguez, J. Sigarreta, Computational and analytical studies of the Randi\'{c} index in Erd\H{o}-R\'{e}nyi models,
\textit{Applied Mathematics and Computation} {\bf 377} (2020) 125137.


\bibitem{MMRS21}
C. T. Mart\'inez-Mart\'inez, J. A. M\'endez-Berm\'udez,   J. M. Rodriguez, J. Sigarreta, 
Computational and analytical studies of
the harmonic index on Erd\H{o}s-R\'{e}nyi models,
\textit{MATCH Commun. Math. Comput. Chem.} {\bf 85} (2021) 395-426.



\bibitem{P18}
K. Pattabiraman, Inverse sum indeg index of graphs, \textit{AKCE International Journal of Graphs and Combinatorics} {\bf 15} (2018) 155-167.





\bibitem{R75}
 M. Randi\'{c}, Characterization of molecular branching, \textit{J. Am. Chem. Soc.} {\bf 97} (1975) 6609–6615.
 
\bibitem{R08}
M. Randi\'{c},
On history of the Randi\'{c} index and emerging hostility toward chemical graph theory,
\textit{MATCH Commun. Math. Comput. Chem.} {\bf 59} (2008) 5-124.

\bibitem{RNP16}
M. Randi\'{c}, M. Novi\'c, D. Plav\v si\'c,
\textit{Solved and unsolved problems in structural chemistry}, CRC Press, Boca Raton, 2016.








\bibitem{SV15}
J. Sedlar, D. Stevanovi\'c, A. Vasilyev,
On the inverse sum indeg index, 
\textit{Discrete Appl. Math.}  {\bf 184} (2015), 202-212.



\bibitem{VG10}
D. Vuki\v cevi\'c, M. Ga\v sperov,
Bond aditive modeling 1. Ariatic indices,
\textit{Croat. Chem. Acta}  {\bf 83} (2010) 243-260.



\bibitem{Y22}
M. Yuan, On the Randi\'{c} index and its variants of network data, submitted to \textit{Test}, 2023. 



\bibitem{ZT09}
B. Zhou,  N. Trinajsti\'c, On a novel connectivity index, \textit{J. Math. Chem.} {\bf 46} (2009)
1252-1270.

\bibitem{ZT10}
B. Zhou,  N. Trinajsti\'c,  On general sum-connectivity index, \textit{J. Math. Chem.} {\bf 47} (2010)
210-218.

\bibitem{Z12}
L. Zhong, The harmonic index for graphs, \textit{Appl. Math. Lett.} {\bf25} (2012) 561-566.


\end{thebibliography}
\end{document}